\documentclass{amsart}

\usepackage{amssymb}
\usepackage{amsmath}
\usepackage{mathtools}
\usepackage{graphicx}
\usepackage{bbm}
\usepackage{array}
\usepackage{eepic}
\usepackage{color}
\usepackage[linktocpage=true]{hyperref}
\usepackage[all]{xy}

\allowdisplaybreaks

\numberwithin{equation}{section}

\newtheorem{theorem}{Theorem}[section]
\newtheorem{lemma}[theorem]{Lemma}
\newtheorem{proposition}[theorem]{Proposition}
\newtheorem{corollary}[theorem]{Corollary}
\newtheorem{remark}[theorem]{Remark}

\newtheorem{definition}[theorem]{Definition}

\renewenvironment{proof}{{\bf \noindent Proof.}}{\qed}
\newcommand{\dem}{\noindent {\bf Proof. }}

\theoremstyle{theorem}
\newtheorem{ltheorem}{Theorem}

\newcommand{\demMax}{\noindent {\bf Proof of Theorem \ref{BoundMultiplierMaximal}. }}
\newcommand{\demA}{\noindent {\bf Proof of Theorem A. }}
\newcommand{\demB}{\noindent {\bf Proof of Theorem B i). }}
\newcommand{\demBB}{\noindent {\bf Proof of Theorem B ii). }}
\newcommand{\demC}{\noindent {\bf Proof of Theorem C. }}
\newcommand{\fin}{\hspace*{\fill} $\square$ \vskip0.2cm}

\newcommand{\N}{\mathbb{N}}

\newcommand{\R}{\mathbb{R}}
\newcommand{\C}{\mathbb{C}}

\newcommand{\summ}{\sum\nolimits}

\DeclareMathOperator*\sotlim{\mathrm{SOT-lim}}

\def\G{\mathrm{G}}
\def\1{\mathbf{1}}
\def\Id{\mathrm{Id}}
\def\H{\mathcal{H}}

\def\M{\mathcal{M}}

\def\B{\mathcal{B}}
\def\CB{\mathcal{C B}}
\def\Ha{\mathcal{H}^{\infty}}
\def\S{\mathcal{S}}

\def\V{\mathcal{L} G }

\def\L{\mathcal{L}}
\def\op{\mathrm{op}}
\def\Aut{\mathrm{Aut}}

\def\mintensor{\otimes_{\min}}

\def\weaktensor{ \, \overline{\otimes} \, }
\def\Fubinitensor{\otimes_{\mathcal{F}}}
\def\osprojtensor{\, \widehat{\otimes} \,}

\def\esssup{\mathop{\mathrm{ess \, sup}}}

\addtocounter{tocdepth}{0}

\begin{document}

\null

\vskip-30pt

\null

\title[Smooth multipliers via Sobolev dimension]{
    Smooth Fourier multipliers \\
    in group algebras via Sobolev dimension}

\author[Gonz\'alez-P\'erez, Junge, Parcet]{
    Adri\'an M. Gonz\'alez-P\'erez,
    Marius Junge, 
    Javier Parcet}

\maketitle

\null

\vskip-30pt

\null

\begin{abstract}
  We investigate Fourier multipliers with smooth symbols defined over locally
  compact Hausdorff groups. Our main results in this paper establish new 
  H\"ormander-Mikhlin criteria for spectral and non-spectral multipliers.
  The key novelties which shape our approach are three. First, we control
  a broad class of Fourier multipliers by certain maximal operators in 
  noncommutative $L_p$ spaces. This general principle ---exploited in 
  Euclidean harmonic analysis during the last 40 years--- is of independent
  interest and might admit further applications. Second, we replace the
  formerly used cocycle dimension by the Sobolev dimension. This is based
  on a noncommutative form of the Sobolev embedding theory for Markov
  semigroups initiated by Varopoulos, and yields more flexibility to measure
  the smoothness of the symbol. Third, we introduce a dual notion of polynomial
  growth to further exploit our maximal principle for non-spectral Fourier
  multipliers. The combination of these ingredients yields new $L_p$ estimates
  for smooth Fourier multipliers in group algebras.     
\end{abstract}


\addtolength{\parskip}{+0.8ex}

\null

\vskip-30pt

\null

\section*{{\bf Introduction \label{Sect0}}}

The aim of this paper is to study Fourier multipliers on group von Neumann algebras for locally compact Hausdorff groups. More precisely, the relation between the smoothness of their symbols and $L_p$-boundedness. This is a central topic in Euclidean harmonic analysis. In the context of  nilpotent groups, it has also been intensively studied in the works of Cowling, M\"uller, Ricci, Stein and others. In this paper we will consider the dual problem, placing our nonabelian groups in the frequency side. Today it is well
understood that the dual of a nonabelian group can only be described as a quantum group, its underlying algebra being the group von Neumann algebra. The interest of Fourier multipliers over such group algebras was recognized early in the pioneering work of Haagerup \cite{Haa1978}, as well as in the research carried out thereafter. It was made clear how Fourier multipliers on these algebras can help in their classification, through the use of certain approximation properties which become invariants of the algebra. Unfortunately, the literature on this topic does not involve the $L_p$-theory ---with a few exemptions like \cite{JunRuan2003} and the very recent paper of Lafforgue and de la Salle \cite{LaffSall2011}--- as it is mandatory from a harmonic analysis viewpoint. In this respect, our work is a continuation of \cite{JunMeiPar2014,JunMeiPar2014NoncommRiesz} where $1$-cocycles into finite-dimensional Hilbert spaces were used to lift multipliers from the group into a more Euclidean space. This yields H\"ormander-Mikhlin type results depending of the dimension of the Hilbert space involved. Here, we shall follow a different approach through the introduction of new notions of dimension allowing more room for the admissible class of multipliers. These notions rely on noncommutative forms of the Sobolev embedding theory for Markov semigroups, which carrie an \lq encoded geometry\rq${}$ in the commutative setting. Prior to that, we need to investigate new maximal bounds whose Euclidean analogues are central in harmonic analysis. In this
paper we shall limit ourselves to unimodular groups to avoid technical issues concerning modular theory.  

\renewcommand{\theequation}{$\mathrm{W}L_2$}
\addtocounter{equation}{-1}

This text is divided into three blocks which are respectively
devoted to maximal bounds, Sobolev dimension and polynomial co-growth.
Let us first put in context our maximal estimates for Fourier multipliers.
Given a symbol $m: \R^n \to \C$ with corresponding Fourier multiplier $T_m$,
there is a long tradition in identifying maximal operators $\mathcal{M}$
which satisfy the weighted $L_2$-norm inequality below for all admissible
input functions $f$ and weights $w$
\begin{equation}
  \label{EuclidenMaxWeight}
  \int_{\R^n} |T_mf|^2 w \, \lesssim \, \int_{\R^n} |f|^2 \mathcal{M}w.
\end{equation}
It goes back to the work of C\'ordoba and Fefferman in the 70's. This
general principle has deep connections with Bochner-Riesz 
multipliers and also with $A_p$ weight theory. The Introduction of \cite{Benn2014}
gives a very nice historical summary and new results in this direction. The
main purpose of this estimate is that elementary duality arguments yield for
$p>2$ that $$\big\| T_m: L_p(\R^n) \to L_p(\R^n) \big\| \, \lesssim \, \big\| \mathcal{M}: L_{(p/2)'}(\R^n) \to L_{(p/2)'}(\R^n) \big\|^{\frac12}.$$ The most general noncommutative form of this inequality would require too
much terminology for this Introduction. Instead, let us just introduce the
basic concepts to give a reasonable but weaker statement. Stronger results
will be given in the body of the paper. Let $G$ be a locally compact 
Hausdorff group. If we write $\mu$ for the left Haar measure of $G$ and 
$\lambda$ for the left regular representation $\lambda: G \to \B(L_2 G)$, the group 
von \hskip-1pt Neumann algebra \hskip-1pt $\V$ is the weak operator closure
in $\B(L_2 G)$ of $\lambda(L_1(G))$. We refer to Section \ref{Sect1} for a
construction of the Plancherel weight $\tau$ on $\V$, a noncommutative
substitute of the Haar measure. Note that $\tau$ is tracial iff $G$ is
unimodular ---which we assume--- and it coincides with the finite trace
given by $\tau(x) = \langle \delta_e, x \delta_e \rangle$ when $G$ is 
discrete. In the unimodular case, $(\V, \tau)$ is a semifinite von Neumann
algebra with a trace and it is easier to construct the noncommutative
$L_p$-spaces $L_p(\V, \tau)$ with norm $\|x\|_p=\tau ( |x|^p)^{1/p}$,
where $|x|^p = (x^* \hskip-2pt x)^{p/2}$ by functional calculus on the
(unbounded) operator $x^* \hskip-2pt x$. Given a bounded symbol 
$m: G \to \C$, the corresponding Fourier multiplier is densely defined by
$T_m \lambda(f) = \lambda(m f)$. Alternatively, it will be useful to 
understand these operators as convolution maps in the following way
\[
  T_m (x) = \lambda(m) \star x = (\tau \otimes \Id) \big( \delta \lambda (m) \, ( \sigma x \otimes \1) \big),
\]
where $\delta : \V \to \V \weaktensor \V$ is determined by
$\delta( \lambda_g ) = \lambda_g \otimes \lambda_g$ and $\sigma : \V \to \V$
is the anti-automorphism given by linear extension of
$\sigma(\lambda_g) = \lambda_{g^{-1}}$. The first map is called the
comultiplication map for $\V$, whereas $\sigma$ is the corresponding
coinvolution. Our next ingredient is the $L_p$-norm of maximal operators.
Given a family of noncommuting operators $(x_\omega)_\omega$ affiliated
to a semifinite von Neumann algebra $\M$, their supremum is not well-defined.
We may however consider their $L_p$-norms through
\[
  \Big\| \, {\sup_{\omega \in \Omega}}^+ x_\omega \Big\|_{L_p (\M)}
  =
  \big\| (x_\omega)_{\omega \in \Omega } \big\|_{L_p(\M; L_\infty(\Omega))},
\]
where the mixed-norm $L_p(L_\infty)$-space has a nontrivial definition obtained by Pisier for hyperfinite $\M$ in 
\cite{Pi1998} and later generalized in \cite{Jun20024,JunPar2010}. This definition recovers
the norm in $L_p(\Sigma; L_\infty(\Omega))$ for abelian 
$\M = L_\infty(\Sigma)$, further details in Section \ref{Sect1}. Finally,
conditionally negative lengths $\psi:G \to \R_+$ are symmetric functions
vanishing at the identity $e$ which satisfy
$\sum_{g,h} \overline{a_g} a_h \psi(g^{-1}h) \le 0$ for any family of
coefficients with $\sum_g a_g = 0$. Due to its one-to-one relation to
Markov convolution semigroups, they will play a crucial role in this paper.
In the classical multiplier theorems, the symbols $m$ are cut
out with functions $\eta(|\xi|)$ for some compactly supported 
$\eta \in C^\infty(\R_+)$. Our techniques do not allow us to
use compactly supported functions in $\R_+$. Instead, we are going to use
analytic functions decaying fast near $0$ and near $\infty$. We will call 
such $\eta$ an $\H^\infty_0$-cut-off, see Section \ref{Sect1} for the precise
definitions. The archetype of such functions will be $\eta(z) = z \, e^{-z}$. 

\renewcommand{\theequation}{\arabic{equation}}
\numberwithin{equation}{section}

\begin{ltheorem}
  \label{TheoremA}
  Let $G$ be a unimodular group equipped with any conditionally negative length
  $\psi : G \rightarrow \R_+$. Let $\eta$ be an $\H_0^\infty$-cut-off and
  $m: G \to \C$ an essentially bounded symbol constant on
  $G_0 = \{g \in G : \psi(g)=0\}$. Assume $B_t = \lambda ( m \eta ( t \psi))$
  admits a decomposition $B_t = \Sigma_t M_t$ with $M_t$ positive and
  satisfying $M_t = \sigma M_t$, and consider the convolution map 
  $\mathcal{R}(x) = ( |M_t|^2 \star \, x )_{t \ge 0}$. Then the following
  inequality holds for $2 < p < \infty$
  \[
    \| T_m \|_{\B( L_p(\V))} \lesssim_{(p)} \, \Big( \sup_{t \geq 0} \| \Sigma_t \|_2 \Big) \Big\| \mathcal{R}: L_{(p/2)'}(\V) \to L_{(p/2)'}(\V; L_\infty(\R_+)) \Big\|^\frac{1}{2}.
  \]
\end{ltheorem}

By duality, a similar stamens holds for $1 < p < 2$. Moreover, a stronger
result holds in terms of noncommutative Hardy spaces which allows more general
symbols and decompositions. Theorem \ref{TheoremA} combines in a very neatly way
noncommutative generalizations of \eqref{EuclidenMaxWeight} with square
function estimates. In the particular case of H\"ormander-Mikhlin symbols ---as we
shall see along this paper--- the decomposition splits the 
assumptions. Namely, the $L_2$-norm of $\Sigma_t$ is bounded using
the smoothness condition while the maximal $\mathcal{R}$ is bounded
through the geometrical assumptions regarding the 
dimensional behaviour of $\psi$. Apart from the direct
consequences given in the present paper, this result is of independent interest
and admits potential applications in other directions to be explored in a
forthcoming publication.

Given a conditionally negative length $\psi: G \to \R_+$, the infinitesimal 
generator of the semigroup $\lambda_g \mapsto \exp(-t \psi(g)) \lambda_g$ is
the map determined by $A(\lambda_g) = \psi(g) \lambda_g$. In particular,
$\psi$-radial Fourier multipliers fall in the category of spectral operators
of the form $m(A)$. These maps are known as spectral multipliers and play a
central role in the theory. Our aim in this second block is to find smoothness
criteria on $m$ which implies $L_p$-boundedness of the spectral multiplier
$T_{m \circ \psi}$.

It is well understood ---specially after \cite{Cow1983,Ste1970}--- that if we want to obtain $L_p$ boundedness for
$m(A)$ from the smoothness of $m$, for every semigroup, we need to impose 
analyticity on $m$. To obtain an smoothness condition with a finite number
of derivatives our space needs to be finite-dimensional. Indeed,
it is known that the optimal smoothness 
order may growth with the dimension. This indicates the
necessity of defining a notion of dimension in the
non-commutative setting. We will take as dimension the value $d>0$ for which a
Sobolev type embedding holds for $A$. Recall that there is a Sobolev
embedding theory for Markov semigroups introduced by Varopoulos \cite{Va1985}.
More precisely, given a measure space $(\Omega,\mu)$ and certain elliptic
operator $A$ generating the Markov process $S_t = \exp(-tA)$, one can
introduce the Sobolev dimension $d$ for which the equivalence below holds
\[
  \|f\|_{L_{\frac{2d}{d-2}}(\Omega)} \, \lesssim \, \big\| A^{\frac12} f \big\|_{L_2(\Omega)}
  \ \Longleftrightarrow \
  \big\| S_t f \big\|_{L_\infty(\Omega)} \, \lesssim \, t^{-\frac{d}{2}} \|f\|_{L_1(\Omega)}.
\]
The property of the right hand side is known as ultracontractivity. When it
holds for the semigroup generated by an invariant Laplacian on a Lie  group,
it forces $\mu(B_t(e)) \sim t^d$. Thus, we can understand ultracontractivity
as a way of describing the growth of balls. With that motivation we
introduce general ultracontractivity properties
\[
  \big\| S_t : L_1(\M) \to \M \big\|_\mathrm{cb} \lesssim \frac{1}{\Phi(\sqrt{t})}.
\]
where $\mathrm{cb}$ stands for completely bounded. The function $\Phi$ will
measure the ``growth of the balls''. Since doubling measure spaces are widely
recognized as a natural setting for performing harmonic analysis, we will
impose $\Phi$ to be doubling, i.e.:
\[
  \sup_{t > 0} \left\{ \frac{\Phi(2 t)}{\Phi(t)} \right\} < \infty,
\]
and our doubling dimension will be given by
\[
  D_\Phi = \log_2 \, \sup_{t > 0} \left\{ \frac{\Phi(2 t)}{\Phi(t)} \right\}.
\]
In the classical abelian setting, apart from the ultracontractivity
---or on-diagonal behaviour of $S_t$--- we need to impose off-diagonal decay on
$S_t$, typically Gaussian bounds. Let $(G,\psi,X)$ be a triple formed by a
locally compact Hausdorff unimodular group $G$, a conditionally negative length
$\psi: G \to \R_+$ and an element $X$ in the extended positive cone
$\V_+^\wedge$, see \cite{Haa1979,Haa1979ii} for precise definitions. We will
say that the triple satisfies the \emph{standard assumptions} when:
\begin{itemize}
  \item[i)] \textbf{Doublingness} 
  \vskip-10pt 
  \null \hskip20pt $\displaystyle \Phi_X(s) = \tau (\chi_{[0,s)}(X))$ is doubling.
  \vskip5pt
  \item[ii)] \textbf{$L_2$ Gaussian upper bounds} 
  \vskip-10pt 
  \null \hskip20pt $\displaystyle \tau \Big( \chi_{[r,\infty)} (X) \big| \lambda(e^{-s \psi}) \big|^2 \Big) \, \lesssim \, \frac{e^{-\beta \frac{r^2}{s}}}{\Phi_X(\sqrt{s})} \quad \mbox{for some } \beta > 0.$
  \vskip5pt
  \item[iii)] \textbf{Hardy-Littlewood maximality} 
  \vskip-8pt
  \null \hskip20pt $\displaystyle \left\| \, {\sup_{s \ge 0}}^+ \frac{\chi_{[0,s)}(X)}{\Phi_X(s)} \star \hskip1pt x \right\|_{L_p(\V)} \hskip-1pt \lesssim_{(p)} \hskip-1pt \|x\|_{L_p(\V)} \quad \mbox{for every } 1 < p < \infty.$
\end{itemize}
We will also require the inequality $\mathrm{iii)}$ to hold
uniformly for matrix amplifications. As we shall see, inequality $\mathrm{ii)}$ implies ultracontractivity with $\Phi_X$
as growth function. We will omit the dependency of $X$ from $\Phi_X$ 
when it can be understood from the context. It is also interesting to point
out that, in the classical case, Gaussian bounds can be deduced from the
ultracontractivity in the presence of geometrical assumptions like locality
or finite speed of propagation for the wave equation, see \cite{Si1996, Si2004}
and \cite[Section 3]{Sa2009}. Generalizing such results to the noncommutative
setting will be the object of forthcoming research. The connection of
standard assumptions with smooth $\psi$-radial Fourier multipliers is nearly
optimal.

\begin{ltheorem}
  \label{TheoremB}
  Let $(G, \psi, X)$ be any triple satisfying the standard assumptions which we considered
  above. Given an $\H_0^\infty$-cut-off function $\eta$ and a symbol $m: \R_+ \to \C$, the following inequalities hold for $1 < p < \infty \hskip-1pt :$
  \begin{itemize}
    \item[i)] If $s > (D_\Phi+1)/2$
    \[
      \|T_{m \circ \psi}\|_{\CB(L_p(\V))} \, \lesssim \, \sup_{t \ge 0} \big\| m(t \cdot) \eta(\cdot) \big\|_{W^{2,s}(\R_+)}.
    \]
    \item[ii)] If $s > D_\Phi/2$ and $\psi \in \mathrm{CBPlan}_q^\Phi$
    for some $q \ge 2$
    \[
      \|T_{m \circ \psi}\|_{\CB(L_p(\V))} \, \lesssim \, \sup_{t \ge 0} \big\| m(t \cdot) \eta(\cdot) \big\|_{W^{q,s}(\R_+)}.
    \]
  \end{itemize}
  The last inequality holds with $q = \infty$ under the sole assumption of
  $s > D_\Phi/2$.
\end{ltheorem}

The term $\CB$ also stands for ``complete bounded'' and the property
$\mathrm{CBPlan}_q^\Phi$ plays the role of the $q$-Plancherel
property introduced by Duong-Ouhabaz-Sikora \cite{DuOuSi2002}, see the body
of the paper for concrete definitions. The proof of Theorem B is the most 
technical in this paper. It will explain the 
decoupling nature of Theorem \ref{TheoremA}. The $\Sigma_t$ are controlled
using the Sobolev smoothness (via the Phragmen-Lindel\"of theorem)
for any degree $s > 0$, whereas the maximal bound determines optimal
restrictions in terms of the Sobolev dimension $D_\Phi$.

Theorem B should be illustrated with interesting examples. The existence
of natural triples satisfying the standard assumptions for nonabelian groups
is the subject of current research, which will appear elsewhere. In this
paper we shall construct such triples out of finite-dimensional cocycles.
This permits to recover the results in
\cite{JunMeiPar2014,JunMeiPar2014NoncommRiesz} for $\psi$-radial multipliers.
In fact, we should emphasize at this point that the notion of dimension in
the previous approach was limited to the Hilbert space dimension of the
cocycle determined by the length $\psi$. Working with finite-dimensional cocycles is an unfortunate
limitation which we could remove for noncommutative Riesz transforms in
\cite{JunMeiPar2014NoncommRiesz}. Theorem B allows to go even further
for smooth radial multipliers.

In our third and last block of this paper, we consider general (non-spectral)
Fourier multipliers. Apart for the semigroup over $\V$ generated by $\psi$ we
will endow $G$ with two semigroups $\S_1/\S_2 : L_\infty(G) \to L_\infty(G)$
of left/right invariant operators. The intuition here is that $\S_j$ will
describe the geometry of $G$ while the semigroup generated by $\psi$ will
describe the geometry of its dual. If $A$ denotes the infinitesimal generator
of a semigroup over $L_\infty(G)$, we use the standard notation for its
nonhomogeneous Sobolev spaces
\[
  \|f\|_{W_A^{p,s}(G)} = \big\| ( \mathbf{1} + A )^{\frac{s}{2}} f \big\|_{L_p(G)}.
\]
When $\mathcal{S}$ is left invariant there exists a positive densely defined
operator $\widehat{A}$ affiliated to $\V$ such that 
$\lambda(Af) = \lambda(f) \widehat{A}$ for all $f \in \mathrm{dom}_2(A)$. In 
a similar way we obtain $\lambda(Af) = \widehat{A} \lambda(f)$ when $\S$ is
right invariant, see Proposition \ref{multSymbols} for the proof. Then we define
the \emph{polynomial co-growth of $\widehat{A}$} as follows
\[
  \mathrm{cogrowth}(\widehat{A})
  \, = \,
  \inf \left\{ r>0 : \big( \mathbf{1} + \widehat{A} \hskip2pt \big)^{-\frac{r}{2}} \in L_1(\V) \right\}.
\]
Our choice for the term ``polynomial co-growth'' sits on the intuition that
$\widehat{A}$ behaves like $|\xi|^2$ in the case of the Laplacian $\Delta$ on
$\R^D$ and therefore $\mathrm{cogrowth}(\widehat{\Delta}) = D$ follows from
the fact that large balls grow like $r^D$. Further in Section \ref{Sect3} we
will characterize polynomial co-growth by relating the behavior of small
balls in $G$ with ``large balls'' in $\V$, see Remark
\ref{LocalvsAssymptotic} for further explanations. It is also worth mentioning the close relation between polynomial growth and Sobolev dimension as it will be analyzed in the body of the paper. Our main result in this direction is the following criterium for non-spectral multipliers.

\begin{ltheorem}
  \label{TheoremC}
  Let $G$ be a unimodular group equipped with a conditionally negative length
  $\psi$.  Let $\S_1/\S_2$ be respectively left/right invariant submarkovian
  semigroups on $L_\infty(G)$ whose generators $A_j$ satisfy
  $\mathrm{cogrowth}(\widehat{A}_j) = D_j$ for $j=1,2$. Consider an
  $\H_0^\infty$-cut-off function $\eta$ and a symbol $m: G \to \C$ which is 
  constant in the subgroup $G_0 = \{g \in G : \psi(g)=0 \}$. Then, if $s_j > D_j/2$ for
  $j=1,2$, the following inequality holds for $1 < p < \infty$
  \[
    \| T_m \|_{\CB(L_p(\V))}
    \, \lesssim_{(p)} \,
    \sup_{t \ge 0} \, \max \Big\{ \big\| \eta(t \psi) \, m \big\|_{W_{A_1}^{2,s_1}(G)}, \big\| \eta(t \psi) \, m \big\|_{W_{A_2}^{2,s_2}(G)} \Big\}.
  \]
\end{ltheorem}

Theorem \ref{TheoremC} establishes a link between the, a priori unrelated,
geometries which determine $\psi$ and $\S_j$. Indeed, we use the
length $\psi$ to cut $m$ ---determining the size of the support--- and use $A_j$ to measure the smoothness of $m$. It is interesting to
note that passing to the dual requires a size condition on
$\widehat{A}$, reinforcing the intuition that duality 
switches size and smoothness. The main difference with Theorem 
\ref{TheoremB} is that in this general context we 
have been forced to place the dilation in the cut-off function $\eta$
instead of the multiplier $m$. We conclude the paper illustrating
Theorem \ref{TheoremC} for Lie groups of polynomial growth by means
of the subriemannian metrics determined by sublaplacians, see
Corollary \ref{LieGroup}.

\section{{\bf Maximal bounds \label{Sect1}}}

\subsection{{\hskip3pt Preliminaries}}

Although the material here exposed is probably well-known to experts, let us
review some notions and results in the interface between harmonic analysis 
and operator algebra that we will need throughout this section. We will start 
with a brief exposition of noncommutative integration theory, including the 
construction of noncommutative $L_p$ spaces. Our main example will be the group
von Neumann algebra of an unimodular Lie group equipped with its canonical 
Plancherel trace. Then we will review some basics of operator space theory 
as well as the construction of certain mixed-norm spaces. Finally
we will consider Markov semigroups with an special emphasis on semigroups of 
convolution type. We will revisit Hardy spaces and square function estimates
associated with a semigroup.

\subsubsection{{\bf \hskip3pt Noncommutative $L_p$ spaces}}

Part of von Neumann algebra theory has evolved as the noncommutative form of 
measure theory and integration. A von Neumann algebra $\M$
\cite{KaRing1997,TaI,TaII}, is a unital weak-operator closed 
$\mathrm{C}^*$-subalgebra of $\B(\H)$, the algebra of bounded linear operators on
a Hilbert space $\H$. We will write $\1_\M$, or simply $\1$, for the unit. 
The positive cone $\M_+$ is the set of positive operators in $\M$ and a trace
$\tau: \M_+ \to [0,\infty]$ is a linear map satisfying $\tau(x^* x) = \tau(x x^*)$. 
Such map is said to be normal if 
$\sup_\alpha \tau(x_\alpha) = \tau(\sup_\alpha x_\alpha)$ for 
bounded increasing nets $(x_\alpha)$; it is semifinite if for 
$x \in \M_+ \setminus \{0\}$ there exists $0 < x' \le x$ with 
$\tau(x') < \infty$; and it is faithful if 
$\tau(x) = 0$ implies $x = 0$. The trace $\tau$ plays the role of the 
integral in the classical case. A von Neumann algebra $\M$ is semifinite 
when it admits a normal semifinite faithful (n.s.f. in short) trace
$\tau$. Any $x \in \M$ is a linear combination $x_1 - x_2 + i x_3 - i x_4$ of
four positive operators. Thus, $\tau$ extends as an unbounded operator to 
nonpositive elements and the tracial property takes the familiar form 
$\tau(x y) = \tau(y x)$. The pairs $(\M,\tau)$ composed by a von Neumann 
algebra and a n.s.f. trace will be called
\emph{noncommutative measure spaces}. Note that commutative von Neumann 
algebras correspond to classical measurable spaces.

By the GNS construction, the noncommutative analogue of measurable sets 
(characteristic functions) are orthogonal projections. Given $x \in \M_+$, 
its support is the least projection $q$ in $\M$ such that $q x = x = x q$ 
and is denoted by $\mbox{supp} \hskip1pt x$. Let $\mathcal{S}_\M^+$ be the 
set of all $f \in \M_+$ such that $\tau(\mbox{supp} \hskip1pt f) < \infty$ 
and set $\mathcal{S}_\M$ to be the linear span of $\mathcal{S}_\M^+$. If we 
write $|x|=\sqrt{x^*x}$, we can use the spectral measure 
$d E$ of $|x|$ to observe that
\[
  x \in \mathcal{S}_\M \Rightarrow |x|^p = \int_{\R_+} s^p \, d E(s) \in \mathcal{S}_\M^+ \Rightarrow \tau(|x|^p) < \infty.
\]
If we set $\|x\|_p = \tau( |x|^p )^{\frac1p}$, we obtain a norm in 
$\mathcal{S}_\M$ for $1 \le p < \infty$. By the strong density of 
$\mathcal{S}_\M$ in $\M$, the \emph{noncommutative $L_p$ space} $L_p(\M)$
is the corresponding completion for $p < \infty$ and $L_\infty(\M) = \M$. 
Many basic properties of classical $L_p$ spaces like duality, real and 
complex interpolation, H\"older inequalities, etc hold in this setting. 
Elements of $L_p(\M)$ can be described as measurable operators 
affiliated to $(\M,\tau)$, we refer to Pisier/Xu's survey \cite{PiXu2003} 
for more information and historical references. Note that classical 
$L_p$ spaces $L_p(\Omega,\mu)$ are denoted in this terminology as 
$L_p(\M)$ where $\M$ is the commutative von Neumann algebra 
$L_\infty(\Omega,\mu)$.

\subsubsection{{\bf \hskip3pt Group algebras and comultiplication formulae}}

Our main example of noncommutative measure space in this paper is that of group von Neumann algebra. Let $G$ be a locally compact and Hausdorff group (LCH group in short) equipped with its left Haar measure $\mu$. Let $\lambda : G \to \B(L_2G)$ be the left regular representation. We will also use $\lambda$ to denote the linear extension of $\lambda$ to the space $L_1(G)$. We will denote by $C_{\lambda}^{\ast} G$ the norm closure of $\lambda(L_1(G))$ and by $\V$ the closure of $C_{\lambda}^{\ast} G$ in the weak operator topology. $\V$ is usually referred to as the \emph{group von Neumann algebra} associated to $G$. There is a distinguished normal faithful weight $\tau : \V_+ \to \R_+$ such that $\lambda: L_1(G) \cap L_2(G) \to \V$ extends to an isometry from $L_2(G)$ to $L_2(\V, \tau)$, the GNS construction associated to $\tau$. Such weight is unique and it is called the Plancherel weight. When the function $f$ belongs to the dense class $C_c(G) \ast C_c(G)$ we have $\tau(\lambda(f)) = f(e)$. The Placherel weight is tracial if and only if $G$ is unimodular. In this case it is called the Placherel trace. From now on we will focus on unimodular groups. We will often work with the spaces $L_p(\V,\tau)$ although the dependency on $\tau$ will be dropped in our terminology.


$\V$ has a natural comultiplication given by linear extension of $\delta( \lambda_g ) = \lambda_g \otimes \lambda_g$ which extends to a $\ast$-homomorphism $\delta : C_{\lambda}^{\ast} G \to C_{\lambda}^{\ast} G \mintensor C_{\lambda}^{\ast} G$.  
There is a unique normal extension $\delta : \V \to \V \weaktensor \V$. This is a consequence of the fact that if $\delta$ is normal then $\delta_{\ast} : \V_{\ast} \osprojtensor \V_{\ast} \to \V_{\ast}$. Here $\mintensor$ and $\osprojtensor$ represent respectively the minimal and projective o.s. tensor products \cite{Pi2003} and $\weaktensor$ denotes the weak operator closure of the algebraic tensor product. Identifying $\L (G \times G)_{\ast}$ with $\V_{\ast} \osprojtensor \V_{\ast}$ we have
\[
  \delta_{\ast} \Big( \int_{G \times G} f(g_1,g_2) \lambda_{(g_1,g_2)} d \mu(g_1) d\mu(g_2) \Big) = \int_{G} f(g,g) \lambda_g d \mu(g),
\]
for every $f \in C_c(G) \ast C_c(G)$. The boundedness of $\delta_{\ast}$ is 
then a consequence of the Herz restriction theorem \cite{Herz1972}. It is 
interesting to note that the Plancherel weight can be characterized as 
the unique normal, nontrivial and $\delta$-invariant weight, where 
$\delta$-invariant means that
\[
  ( \tau \otimes \Id ) \delta x = \tau(x) \1.
\]
Analogously, Fourier multipliers are characterized as $\delta$-equivariant maps
\[
  \delta T = (T \otimes \Id) \delta = (\Id \otimes T) \delta.
\]
We will denote by $\sigma : \V \to \V$ the anti-automorphism 
given by linear extension of $\sigma(\lambda_g) = \lambda_{g^{-1}}$. 
The \emph{quantized convolution} of two elements $x,y$ affiliated to $\V$ is defined by
\[
  x \star y = (\tau \otimes \Id) \big( \delta x \, ( \sigma y \otimes \1) \big).
\]
Observe that given $m \in L_\infty(G)$, the corresponding Fourier multiplier has the form
\[
  T_m(x) = \lambda(m) \star x = (\tau \otimes \mathrm{Id}) \big( \delta \lambda (m) \, ( \sigma x \otimes \1) \big).
\]

\subsubsection{{\bf \hskip3pt Operator space background}}

The theory of operator spaces is regarded  as a noncommutative or quantized form of 
Banach space theory. An \emph{operator space} $E$ is a closed subspace of
$\B(\H)$. Let $M_m(E)$ be the space of 
$m \times m$ matrices with entries in $E$ and impose on it 
the norm inherited from $M_m(E) \subset \B(\H^m)$. 
The morphisms in this category are the \emph{completely bounded}
linear maps (c.b. in short) $u: E \to F$, i.e. those 
satisfying
\[
  \| u \|_{\CB(E, F)} = \big\| u: E \to F \big\|_{\mathrm{cb}} = \sup_{m \ge 1} \big\| \Id_{M_m} \otimes u \big\|_{\B(M_m(E), M_m(F))} < \infty.
\]
Similarly, given $\mathrm{C}^*$-algebras $A$ and $B$, a linear map 
$u: A \to B$ is called completely positive (c.p. in short)
when $\Id_{M_m} \otimes u$ is positivity preserving
for $m \ge 1$. When a c.p. map $u: A \to B$ is contractive 
(resp. unital) we will say it is a c.c.p. (resp. u.c.p.) map. 
The Kadison-Schwartz inequality for a c.c.p. map $u: \M \to \M$ claims that
\[
  u(x)^* u(x) \le u(x^*x) \quad \mbox{for all} \quad x \in \M.
\]
Ruan's axioms describe axiomatically those sequences of matrix norms 
which can occur from an isometric embedding into $\B(\H)$. Admissible 
sequences of matrix norms are called operator space structures 
(o.s.s. in short) and become crucial in the theory. Given a 
Banach space $X$ and an isometric embedding $\rho : X \to \B(\H)$ 
we will denote by $X^ {\rho}$ the corresponding operator space. Central 
branches from the theory of Banach spaces like duality, tensor 
norms or complex interpolation have been successfully extended 
to the category of operator spaces. Rather complete expositions 
are given in {\cite{EfRu2000,Pau2002,Pi2003}}. Two particular aspects of
operator space theory which are relevant in this paper are the following:

\noindent \textbf{A. Vector-valued Schatten classes.} We will denote by $S_p$ the Schatten $p$-class given by 
$S_p = L_p (\B( \ell_2), \mathrm{Tr})$ with $\mathrm{Tr}$ 
the standard trace in $\mathcal{B}(\ell_2)$. Similarly, $S_p^m$ 
stands for the same space over $m \times m$ matrices. Vector-valued 
forms of these spaces can be defined as long as we define an o.s.s. 
over the space where we take values. Given an operator space $E$, 
we may define the \emph{$E$-valued Schatten classes} $S_p^m[E]$ as the 
operator spaces given by interpolation 
\[
  S_p^m[E] \, := \, \big[ S_\infty^m[ E], S_1^m[E] \big]_{\frac1p} \, := \, \big[ S_{\infty}^m \mintensor E, S_1^m  \osprojtensor E \big]_{\frac1p}.
\]
These classes provide a useful characterization of complete boundedness $$\| u \|_{\CB(E,F)} = \sup_{m \ge 1} \big\| \Id_{M_m} \otimes u \big\|_{\B(S_p^m(E),S_p^m(F))} \quad \mbox{for} \quad 1 \le p \le \infty.$$ For a general hyperfinite von Neumann algebra $\M$ the construction of $L_p(\M;E)$ is carried out by direct limits of $E$-valued Schatten classes. We refer to Pisier's book \cite{Pi1998} for more on vector-valued noncommutative $L_p$ spaces. The space $L_p(\M;E)$ for nonhyperfinite $\M$ cannot be constructed without losing fundamental properties like projectivity/injectivity of the functor $E \mapsto L_p (\M; E)$. Fortunately, this drawback is solvable for the vector-valued $L_p$ space we shall be working with.

\noindent \textbf{B. Operator space structure of $L_p$.} Given an operator space $E$, its opposite $E_{\op}$ is 
the operator space which comes equipped with the operator space 
structure determined by the o.s.s. of $E$ as follows
\[
  \Big\| \sum_{j,k=1}^m a_{jk} \otimes e_{jk} \Big\|_{M_m(E_{\op})} \, = \, \Big\| \sum_{j,k=1}^m a_{kj} \otimes e_{jk} \Big\|_{M_m(E)},
\]
where $e_{jk}$ stand for the matrix units in $M_m$. Alternatively, if
$E \subset \B(\H)$, then $E_\op = E^{\top} \subset \B(\H)$,
where $\top$ is the transpose. The $\op$ construction plays a 
role in the construction of a \lq\lq natural'' 
o.s.s. for noncommutative $L_p$ spaces. If $\M$ is a von Neumann 
algebra we will denote by $\M_{\op}$ it opposite algebra, 
the original algebra with the multiplication reversed. It is a 
well-known result that $\M_{\op}$ and $\M$ need not be 
isomorphic \cite{Co1975}. For every operator space $E$ the natural 
inclusion $j : E \rightarrow E^{\ast \ast}$ is a complete isometry. 
This allows us to build an operator space structure in the predual 
$\M_{\ast}$ as the only operator space structure that makes the 
inclusion $j : \M_{\ast} \rightarrow \M^{\ast}$ completely isometric. 
The operator space structure of $L_p (\M)$ is given by operator space 
complex interpolation between $L_1 (\M) = (\M_{\op})_{\ast}$ 
and $\M$. In particular, it turns out that 
\[
  L_p(\M)^* \simeq L_{p'}(\M_{\op})
\]
is a complete isometry for $1 \le p < \infty$, see 
\cite[pp. 120-121]{Pi2003} for further details. 

\subsubsection{{\bf \hskip3pt $L_\infty$-valued $L_p$ spaces}}

Maximal inequalities are a cornerstone in harmonic analysis. Unfortunately, the supremun of a family of noncommuting operators is not well-defined, so that we do not have a proper noncommutative analogue of maximal functions. Nevertheless, this difficulty can be overcome if all we want is to bound is the maximal function in noncommutative $L_p$, as usually happens in harmonic analysis for commutative spaces. In that case we exploit the fact that the $p$-norm of a maximal function can always be written as a 
mixed $L_p(L_\infty)$-norm of the corresponding entries. This reduces the problem to construct the vector-valued spaces $L_p(\M; L_\infty(\Omega))$. This construction can be carried out without requiring $\M$ to be hyperfinite, relaying in the commutativity of 
$L_\infty(\Omega)$. $L_p(\M;L_\infty(\Omega))$ is defined as the subspace of functions $x \in L_\infty(\Omega; L_p(\M))$ which admit a factorization of the form $x_{\omega} = \alpha \, y_{\omega} \, \beta$ with $\alpha, \beta \in L_{2p}(\M)$ and $y \in L_\infty(\Omega;\M)$. The norm in such space is then given by
\[
  \big\| (x_\omega)_{\omega \in \Omega} \big\|_{L_p (\M ; L_{\infty} (\Omega))} \, = \, \inf \Big\{ \| \alpha \|_{2p}  \Big( \esssup_{\omega \in \Omega} \| y_\omega \|_\M \Big)  \| \beta \|_{2 p} \, : \, x = \alpha y \beta \Big\}.
\]
When $x_\omega \ge 0$ the norm coincides with 
\begin{equation} \label{PostiveMax}
\big\| (x_\omega)_{\omega \in \Omega} \big\|_{L_p ( \M ; L_{\infty} (\Omega))} \, = \, \inf \Big\{ \|y\|_{L_p(\M)} \, : \, x_\omega \leq y \mbox{ for a.e. $\omega \in \Omega$} \Big\}.
\end{equation}
Its operator space structure satisfies
\[
  S_p^m \big[ L_p ( \M; L_{\infty} ( \Omega)) \big] \, = \, L_p \big( M_m \otimes \M ; L_{\infty} (\Omega) \big).
\]
It is standard to use the following notation for the noncommutative $L_p(L_\infty)$-norm
\[
  \Big\| \, {\sup_{\omega \in \Omega}}^{+} \, x_\omega \, \Big\|_{L_p(\M)} = \big\| ( x_\omega )_{\omega \in \Omega} \big\|_{L_p(\M; L_\infty(\Omega))},
\]
where the $\sup$ is just a symbolic notation without an intrinsic 
meaning. In the proof of Theorem \ref{TheoremB} we will use the
fact that if $(\mu_{\omega_2})_{\omega_2 \in \Omega_2}$ is a family of finite
positive measures in $\Omega_1$ and $(R_{\omega_1})_{\omega_1 \in \Omega_1}$
is a family of positivity preserving operators, then the following bound 
holds for $x \in L_p(\M)_+$
\begin{equation}
  \Big\| \, {\mathop{\mathrm{sup}^+}_{\omega_2 \in \Omega_2}} \Big\{ \int_{\Omega_1} R_{\omega_1}(x) d \mu_{\omega_2}( \omega_1 ) \Big\} \Big\|_p
  \le
  \Big( \sup_{\omega_2 \in \Omega_2} \| \mu_{\omega_2} \|_{M(\Omega)} \Big) \Big\| \, {\mathop{\mathrm{sup}^+}_{\omega_1 \in \Omega_1}} R_{\omega_1}(x) \Big\|_p.
  \label{averageMax}
\end{equation}

When $\M$ is hyperfinite, this definition of $L_p(\M ; L_\infty(\Omega))$ 
coincides with the corresponding vector-valued space as defined by 
Pisier \cite{Pi1998}. This approach to handle 
maximal inequalities in von Neumann algebras has been successfully used in 
\cite{Jun20024} to find noncommutative forms of Doob's maximal inequality for 
martingales and the maximal ergodic inequalities for Markov 
semigroups \cite{JunXu2007}. The predual can be explicitly described as the $L_1$-valued space $L_{p'}(\M;L_1(\Omega))$. Indeed, let $S_p(\Omega)$ be the Schatten
class associated to the Hilbert space $L_2(\Omega)$. Note that there is a
hermitian form $q: L_{2 \, p}(\M) \otimes S_2^c(\Omega)
\times L_{2 \, p}(\M) \otimes S_2^c(\Omega) \to L_p(\M) \otimes L_1(\Omega)$
given by
\[
  q( x \otimes m, y \otimes n ) = x^* y \otimes \mathrm{diag}(m^* n),
\]
where $\mathrm{diag} : S_1(\Omega) \to L_1(\Omega)$ is the restriction
to the diagonal. Define 
\[
  \| x \|_{L_p(\M;L_1(\Omega))} = \inf \Big\{ \| a \|_{L_{2p}(\M;S_2^c(\Omega))} \, \| b \|_{L_{2p}(\M;S_2^c(\Omega))} : q(a,b) = x \Big\}.
\]
This space satisfies that $L_p(\M;L_1(\Omega))^{\ast} = 
L_{p'}(\M_{\mathrm{op}};L_\infty(\Omega))$ for $1 \leq p < \infty$. 

\subsubsection{{\bf \hskip3pt Hilbert-valued $L_p$ spaces}}

For certain operator spaces whose underlying Banach space is a Hilbert space 
we can define vector-valued noncommutative $L_p$ spaces for general von Neumann 
algebras. Indeed, let $\H$ be a Hilbert space and and 
$P_{e} \xi = \langle e, \xi \rangle e$ for some $e \in \H$ of unit 
norm. We define the following two Hilbert-valued forms of $L_p (\M)$
\begin{eqnarray*}
  L_p (\M ; \H^c) & = & L_p(\M \weaktensor \B (\H)) (\Id_\M \otimes P_e), \\
  L_p (\M ; \H^r) & = & (\Id_\M \otimes P_e) L_p ( \M \weaktensor \B (\H)), 
\end{eqnarray*}
called the $L_p$ spaces with $\H$-column (resp. $\H$-row) values. Their 
o.s.s. are the ones inherited from $L_p( \M \weaktensor \B(\H))$. 
If $\H = \ell_2^n$, then we can identify $L_p(\B(\H) \weaktensor \M)$ with 
$L_p (\M)$-valued $n \times n$ matrices. Under that identification 
$L_p (\M ;\H^c)$ (resp. $L_p(\M;\H^r)$) corresponds to the matrices 
with zero entries outside the first column (resp. row) and we have that
\begin{eqnarray*}
  \Big\| \sum_{j=1}^n x_j \otimes e_{j1} \Big\|_{L_p ( \M \otimes \B ( \ell_2^n))} & = & \Big\| \Big( \sum_{j = 1}^n x_j^{\ast} x_j \Big)^{\frac{1}{2}}  \Big\|_{L_p (\M)},\\
  \Big\| \sum_{j=1}^n x_j \otimes e_{1j} \Big\|_{L_p ( \M \otimes \B ( \ell_2^n))} & = & \Big\| \Big( \sum_{j = 1}^n x_j x_j^{\ast} \Big)^{\frac{1}{2}}  \Big\|_{L_p (\M)} .
\end{eqnarray*}
The same formulas hold after replacing the finite sums by infinite ones of 
even by integrals.
For every $1 \leq p \leq \infty$ we can embed $\H$ isometrically
in $S_p$ by sending $c_p(e_j) = e_{1 \, j}$ or $r_p(e_j) = e_{j \, 1}$, where
$\{e_j\}$ is an orthonormal basis of $\H$. Such maps are called
the $p$-column/$p$-row embedings. These isometries endow $\H$ with several o.s.
structures. Observe that, as an o.s, $L_p(\M;\H^c)$ (resp. $L_p(\M;\H^r)$)
coincides with Pisier's vector-valued $L_p$-space $L_p(\M;\H^{c_p})$
(resp. $L_p(\M;\H^{r_p})$) for $\M$ hyperfinite. For $1 \leq p < \infty$
the duals are given by $L_p(\M;\H^c)^* = L_{p'}(\M_\op; \H^c)$ and
$L_p(\M;\H^c)^* = L_{p'}(\M_\op, \H^c)$. The duality pairing can be express as
\[
  \Big\langle \summ_{j} x_j \otimes e_j, \summ_{k} y_k \otimes e_k \Big\rangle = \summ_{j} \tau(x_j^*  y_j).
\]
The spaces $L_p(\M; \H^r)$ and 
$L_p(\M; \H^c)$ form complex interpolation scales for $p \ge 1$
$$	\begin{array}{rcl}
	  \left[ L_{\infty}(\M; \H^r), L_p(\M ;\H^r) \right]_{\theta} & = & L_{\frac{p}{\theta}}(\M ;\H^r),\\
	  \left[ L_{\infty}(\M; \H^c), L_p(\M ;\H^c) \right]_{\theta} & = & L_{\frac{p}{\theta}}(\M ;\H^c).
	\end{array}
$$

In order to treat square functions and Hardy spaces we 
will need to control sums and intersections of these 
Hilbert valued noncommutative $L_p$ spaces. The algebraic tensor product 
$L_p(\M) \otimes \H$ embeds in 
$L_p(\M \weaktensor \B(\H))$ by $\Id \otimes r$ and $\Id \otimes c$.
Taking direct sums we obtain an embedding in 
$X = L_p(\M \weaktensor \B(\H)) \oplus L_p(\M \weaktensor \B(\H))$. The space
$L_p(\M; \H^{r \cap c})$ is defined as the norm closure (or weak-$\ast$ closure
if $p = \infty$) of $L_p(\M) \otimes \H$ inside $X$. Such space 
comes equipped with the norm given by
\[
  \| x \|_{L_p (\M; \H^{r \cap c})} = \max \Big\{ \| x \|_{L_p ( \M ; \H^r)}, \| x \|_{L_p ( \M ; \H^c)} \Big\}.
\]
The embedding also gives $L_p(\M;\H^{r \cap c})$ an o.s.s. 
We will denote the dual spaces by 
$L_p(\M;\H^{r + c}) = L_{p'}(\M_{\op}; \H^{r \cap c})^*$ for
$1 < p \leq \infty$. The space $L_1(\M;\H^{r+c})$ is defined as the 
subset of weak-$\ast$ continuous functionals in 
$L_\infty(\M_{\op}; \H^{r \cap c})^*$. The sum notation 
comes from the fact that
\[
  \| x \|_{L_p (\M ; \H^{r + c})} = \inf \Big\{ \| y \|_{L_p (\M ; \H^r)} + \| z \|_{L_p ( \M ; \H^c)} : x = y + z \Big\}.
\]
We will denote by $L_p(\M; \H^{rc})$ the spaces given by
\[
   L_p(\M ; \H^{r c})
   =
  \begin{cases}
	  L_p(\M ;\H^{r + c})    & \mbox{when  } 1 \leq p < 2\\
	  L_p(\M ;\H^{r \cap c}) & \mbox{when  } 2 \leq p \leq \infty .
  \end{cases}
\]
The spaces $L_p (\M;\H^{rc})$ are crucial for the noncommutative Khintchine 
inequalities \cite{Lust1986, LustPi1991}, the noncommutative Burkholder-Gundy 
inequalities \cite{JunXu2003}, noncommutative Littlewood-Paley estimates {\cite{JunMerXu2006}} and other related results.

\subsubsection{{\bf \hskip3pt Markovian semigroups and length functions}}

\hskip-3pt A semigroup $\S=(S_t)_{t \geq 0}$ over a Banach space $X$ is a family of 
operators $S_t : X \to X$ such that $S_0 = \Id$ and 
$S_t S_s = S_{t + s}$. Let $(\M,\tau)$ be a noncommutative measure 
space, we will say that a semigroup $\S$ over $\M$ is 
\emph{submarkovian} iff:
\begin{itemize}
  \item[i)]   Each $S_t$ is a weak-$\ast$ continuous and c.c.p. map.
  \item[ii)]  Each $S_t$ is a self-adjoint, ie:
  $\tau(S_t x^* y) = \tau(x^* S_t y)$.
  \item[iii)] The map $t \mapsto S_t$ is pointwise weak-$\ast$ continuous.
\end{itemize}
$\S$ is \emph{Markovian} if each $S_t$ is a u.c.p. map, 
ie $S_t(\1) = \1$. Markovian operators satisfy $\tau \circ S_t = \tau$ 
while submarkovian ones satisfy $\tau \circ S_t \le \tau$. Sometimes these semigroups
are called symmetric and Markovian, where symmetric is synonym with 
self-adjoint. All the semigroups in this paper will be symmetric, so we will
drop the adjective. Using the first two properties it is easy to see that $S_t$
extends to a c.c.p. map on $L_1(\M)$. By the Riesz-Thorin theorem $S_t$ is a 
complete contraction over $L_p(\M)$ for $1 \leq p \leq \infty$. The third 
property implies that $t \mapsto S_t$ is SOT continuous in $L_1(\M)$. By 
interpolation between the pointwise norm continuity on $L_1(\M)$ and the 
pointwise weak-$\ast$ continuity on $\M$ we obtain that $t \mapsto S_t$ 
is SOT continuous on $L_p(\M)$ for $1 \leq p < \infty$. For every 
$1 \leq p < \infty$ there is a densely defined and closable operator 
$A$ whose closed domain is given by
\[
  \mathrm{dom}_{p}(A) = \Big\{ x \in L_p(\M) : \exists \lim_{t \to 0^{+}} \frac{x - S_t x}{t} \mbox{ in the norm topology }  \Big\}.
\]
When $p=2$ we have that $S_t = e^{-t A}$ and 
$S_t[L_p(\M)] \subset \mathrm{dom}_p(A)$ for 
$1 \leq p < \infty$. In the case $p=\infty$ we have that $A$ is densely
defined and closable with respect to the weak-$\ast$ topology with domain given
by those $x \in \M$ such that $\lim_{t \to 0^{+}} ( x - S_t x ) / t$ 
exists in the weak-$\ast$ topology. We will call $A$ the infinitesimal
generator of $\S$.

We are interested in those (sub)markovian semigroups 
over $\M=\V$ which are of convolute type. In other words, each $S_t$ is a Fourier multiplier. It can be proved
that $S_t = T_{e^{-t \psi}}$ for some function $\psi$. Let us recall a characterization 
of these functions. First, recall some definitions. 
A continuous function $\psi : G \to \C$ is said to be 
conditionally negative (c.n. in short) iff $\psi(e) = 0$
and for every finite subset $\{ g_1, ..., g_m \} \subset G$ and vector 
$(v_1, ..,v_m) \in \C^n$  we have
\[
  \sum_{i=1}^m v_i = 0 \ \ \Rightarrow \ \sum_{i,j=1}^{m} \bar{v}_i \psi( g_{i}^{-1} g_j ) v_j \leq 0.
\]
When $\psi : G \to \R_+$ is symmetric ($\psi(g) = \psi(g^{-1})$) and c.n. we will say that $\psi$ is a \emph{conditionally negative length}. Let  $\H$ be a real Hilbert space. Given an orthogonal representation $\alpha : G \to O(\H)$ we say that a continuous map $b : G \to \H$ is a \emph{$1$-cocycle} (with respect to $\alpha$) iff it satisfies the $1$-cocycle law $$b(g h) = \alpha(g) b(h)  + b(g).$$ The following characterization is proved in \cite[Appendix C]{BeHarVal2008} or \cite[Chapter 1]{CheCowJoJulgVal2001}.

\begin{theorem}
  \label{genCN}
  Let $\S=(S_t)_{t \geq 1}$ be a 
  semigroup of convolution type over the group algebra $\V$. Then, the following statements 
  are equivalent:
  \begin{itemize}
    \item[i)]   $\S$ is a Markovian semigroup.
    \item[ii)]  There is a c.n. length $\psi : G \to \R_+$ such that 
    $S_t = T_{e^{- t \psi}}$.
    \item[iii)] There is a real Hilbert space $\H$, an orthogonal 
    representation $\alpha : G \to O(\H)$ and a $1$-cocycle 
    $b:G \to \R_+$, such that $\psi(g) = \| b(g) \|_{\H}^2$ and
    $S_t = T_{e^{- t \psi}}$
  \end{itemize}
\end{theorem}

\subsubsection{{\bf \hskip3pt Holomorphic calculus and noncommutative Hardy spaces}}

We now introduce the Hardy spaces associated with a submarkovian semigroup on $(\M,\tau)$ as well as the corresponding $\mathcal{H}^{\infty}$-functional calculus. Both tools were introduced in the noncommutative setting 
in \cite{JunMerXu2006}. If $\S$ is a submarkovian semigroup, the 
fixed point subspace $F_p = \{ x \in L_p(\M) : S_t(x) = x \ \forall \, t \geq 0 \}$
coincides with $\ker A \subset \mathrm{dom}_p(A)$ and it is a subalgebra 
when $p = \infty$. It is also easily seen to be a complemented 
subspace with projection given by $Q_p(x) = \lim_{t \to \infty} S_t x$ where
the limit converges in the norm topology of $L_p$, for $p < \infty$
and in the weak-$\ast$ topology when $p = \infty$. We will denote by 
$L^{\circ}_p(\M) = L_p(\M) / F_p$ which is also a complemented subspace with
projection given by $P_p = \Id - Q_p$. Note that 
$L_p(\M) \simeq L^{\circ}_p(\M) \oplus_{p} F_p$. When $S_t$ are Fourier 
multipliers over $\M = \V$ with symbol $e^{-t \psi}$ we define 
$G_0 = \{ g \in G : \psi(g) = 0 \}$. In that case
\[
  F_p = \overline{\Big\{ x \in L_p(\M) : \, x = \lambda(f) \mbox{ with } \mathrm{supp}(f) \subset G_0 \Big\}}
\]
and in a similar way we find that $\lambda(f) = L^{\circ}_p(\M)$ if and only if $f_{|_{G_0}} = 0$.

For any given $x \in \M$ we define the function $Tx : (0,\infty) \to L_p(\M)$
given by $t \mapsto t \, \partial_t S_t x$. We can see $x \mapsto T x$ as a
map from certain domain $D \subset \M$ into $L_p(\M;\H^r)$, 
$L_p(\M;\H^c)$ or $L_p(\M;\H^{r c})$, where $\H = L_2(\R_+, dt / t)$.
The induced seminorms on $D \subset \M$ are called the 
row Hardy space, column Hardy space or Hardy space seminorms. Observe that
the map $T$ has as kernel those elements fixed by $\S$. Quotient out
the nulspace and taking the completion with respect to any of those norms
when $p < \infty$ (resp. the weak-$\ast$ topology for 
$p=\infty$) gives the Hardy spaces $H^{r}_p(\M;\S)$, 
$H^{c}_p(\M;\S)$ or $H_p(\M;\S)$. We can represent such norms as follows
\begin{eqnarray*}
    \| x \|_{H_p^c(\M;\S)} & = & \Big\| \Big( \int_{\R_+} \Big( t \frac{d}{dt} S_t x \Big)^{\ast} \Big( t \frac{d}{dt} S_t x \Big) \frac{dt}{t} \Big)^{\frac{1}{2}} \Big\|_{L_p(\M)}, \\
    \| x \|_{H_p^r(\M;\S)} & = & \Big\| \Big( \int_{\R_+} \Big( t \frac{d}{dt} S_t x \Big) \Big( t \frac{d}{dt} S_t x \Big)^{\ast} \frac{dt}{t} \Big)^{\frac{1}{2}} \Big\|_{L_p(\M)}. 
\end{eqnarray*}
We will drop the dependency on the semigroup and write $H_p^c(\M)$ whenever 
it can be understood from the context. These spaces inherit their o.s.s. from 
that of $L_p(\M;\H^r)$ or $L_p(\M;\H^c)$. Therefore we have the following 
identities
\begin{eqnarray*}
    S_p^n[H^{c}_p(\M;\S)] & = & H^{c}_p(\M \weaktensor \B(\ell^n_2);\S \otimes \Id),\\
    S_p^n[H^{r}_p(\M;\S)] & = & H^{r}_p(\M \weaktensor \B(\ell^n_2);\S \otimes \Id).
\end{eqnarray*} 
The duality is obtained from that of $L_p(\M;\H^c)$ or $L_p(\M;\H^r)$, resulting in the cb-isometries $H^r_p(\M;\S)^* = H^r_{p'}(\M_{\op};\S)$ for 
$1 \leq p < \infty$. The same holds for the column case. Finally let 
us recall that by \cite[Chapters 7 and 10]{JunMerXu2006} we have that if 
$1 < p < \infty$ then
\begin{equation}
  H_p(\M;\S) \simeq L^{\circ}_p(\M),
  \label{EqHardyLp}
\end{equation}
with the equivalence as operator spaces depending on the constant $p$. 
The result fails for $p=1,\infty$ and $H_1(\M;\S)$ is smaller in general 
than $L_1^{\circ}(\M)$.
Observe that $t \partial_t S_t x = \eta( t A ) x$ where $\eta(z) = z e^{-z}$. 
Due to the results in \cite{JunMerXu2006} we can change $\eta$  by other 
analytic functions in certain class obtaining equivalent norms. We will 
say that a holomorphic function $\rho$ defined over the sector 
$\Sigma_{\theta} = \{ z \in \C : | \arg(z) | < \theta \}$ is in 
$\H^ \infty (\Sigma_{\theta})$ iff it is bounded and 
we will say that it is in 
$\H_0^\infty (\Sigma_{\theta}) \subset \mathcal{H}^\infty (\Sigma_{\theta})$ 
iff there is an $s > 0$ such that
\[
  | \rho(z) | \lesssim \frac{ |z|^s}{(1 + |z|)^{2s}}.
\]
We will denote by $\H^{\infty}$ or $\mathcal{H}_0^{\infty}$ the 
spaces $\bigcap_{0 <\theta < \pi/2} \mathcal{H}^{\infty}(\Sigma_\theta)$ or 
$\bigcap_{0 <\theta < \pi/2} \mathcal{H}^{\infty}_0(\Sigma_\theta)$ 
respectively. If needed, we will equip these spaces with their natural inverse
limit topologies. We have that for any $\rho \in \mathcal{H}^{\infty}_0$ the
following holds
\begin{equation} \label{columnSQF}
\begin{array}{rcl}  \| x \|_{H_p^c(\M)} & \sim_{(p)} & \displaystyle \Big\| \Big( \int_{\R_+} \big( \rho( t A ) x \big)^{\ast} \, \rho( t A ) x \frac{dt}{t} \Big)^{\frac{1}{2}} \Big\|_{L_p(\M)}, \\
  \| x \|_{H_p^r(\M)} & \sim_{(p)} & \displaystyle \Big\| \Big( \int_{\R_+} \rho( t A ) x \,  \big( \rho( t A ) x \big)^{\ast} \frac{dt}{t} \Big)^{\frac{1}{2}} \Big\|_{L_p(\M)}. \end{array}
\end{equation}
The equivalence also holds after matrix amplifications. This type of identities also hold for wider classes of unbounded operators $A$ satisfying certain resolvent estimates, see \cite{JunMerXu2006} for further details.

\subsection{{\hskip3pt The general principle}} 

We are now ready to prove our maximal bounds in Theorem A. In fact, we shall obtain a more general principle in Theorem \ref{BoundMultiplierMaximal} which decouples in terms of row and column Hardy spaces. 

\begin{definition}
  Let $(B_t)_{t \geqslant 0}$ be a family of operators affiliated to
  $\V$. We say that $(B_t)_{t \geq 0}$ has an $L_p$-\emph{square-max 
  decomposition} when there is a decomposition $B_t = \Sigma_t
  M_t$ such that$\, :$
  \begin{equation} \label{SM} \tag{$\mathrm{SM}_p$}
  \begin{array}{rll} 
  \displaystyle \sup_{t \geq 0} \| \Sigma_t \|_2 & < & \infty, \\ [5pt]
  \displaystyle \Big\| {\sup_{t > 0}}^+ \sigma | M_t |^2 \star u \Big\|_p & \lesssim_{(p)} & \| u \|_p.
 \end{array}
 \end{equation}     
Similarly, $(B_t)_{t \geq 0}$ has an 
  $L_{p}$-\emph{max-square decomposition} when $B_t = M_t \Sigma_t$ with$\, :$
  \begin{equation} \label{MS} \tag{$\mathrm{MS}_p$}
  \begin{array}{rll} 
  \displaystyle \sup_{t \geq 0} \| \Sigma_t \|_2 & < & \infty, \\ [5pt]
  \displaystyle \Big\| {\sup_{t > 0}}^+ \sigma | M^{\ast}_t |^2 \star u \Big\|_p & \lesssim_{(p)} & \| u \|_p.
 \end{array}
 \end{equation}     
\end{definition}
When we say that $(B_t)_{t \geq 0}$ has a max-square (resp. square-max) decomposition we mean that it has an $L_p$-max-square (resp. $L_p$-square-max) decomposition for every $1 < p <\infty$.

\begin{theorem}
\label{BoundMultiplierMaximal}
Let $G$ be a LCH group equipped with a conditionally negative length $\psi : G \rightarrow \R_+$. Let $\S= (S_t)_{t \geq 0}$ be the convolution semigroup generated by $\psi$ and pick any $\eta \in \H_0^{\infty}$. If $m : G \rightarrow \C$ is a bounded function satisfying that $B_t = \lambda(m \eta (t \psi))$ has an $L_{(p/2)'}$-square-max decomposition $B_t = \Sigma_t M_t$ for some $2 < p <\infty$, then $T_m : H^c_p(\V) \to H^c_p (\V)$ and $$\| T_m \|_{ \B( H^c_p ) } \lesssim_{(p)} \Big( \sup_{t \geq 0} \| \Sigma_t \|_2 \Big) \, \Big\| (R_t^c)_{t \geq 0} : L_{(p/2)'}(\V) \to L_{(p/2)'}(\V; L_\infty) \Big\|^\frac{1}{2}$$ where $R_t^c(x) = \sigma | M_t |^2 \star x$. Similarly, when $(B_t)_{t \geq 0}$ admits an $L_{(p/2)'}$-max-square decomposition $B_t = M_t \Sigma_t$ for some $2 < p < \infty$, we get $T_m : H^r_p(\V) \rightarrow H^r_p(\V)$ and the following estimate holds $$\| T_m \|_{ \B( H^r_p )} \lesssim_{(p)} \Big( \sup_{t \geq 0} \| \Sigma_t \|_2 \Big) \, \Big\| (R_t^r)_{t \geq 0} : L_{(p/2)'}(\V) \to L_{(p/2)'}(\V;L_\infty) \Big\|^\frac{1}{2}$$ where $R_t^r(x) = \sigma | M_t^{\ast} |^2 \star x$. By duality, similar identities also hold for $1 < p < 2$.
\end{theorem}

\begin{corollary}
\label{ColandRowMultiplierMax}
If $G$, $\psi$, $\eta$ and $m$ are as above and $B_t = \lambda(m \eta( t \psi))$ admits both a $L_{(p/2)'}$-max-square and a $L_{(p/2)'}$-square-max decomposition, then it turns out that $T_m : L_p^{\circ}(\V) \to L_p^{\circ}(\V)$ boundedly. Furthermore, if $m \equiv c$ in $G_0 = \{g \in G : \psi(g)\}$ then $T_m$ is a bounded map on $L_p(\V)$. 
\end{corollary}

\dem
The first assertion follows trivially from \eqref{EqHardyLp}. For the second we use that  $L^{\circ}_p(\V)$ is a complemented subspace, and so 
\begin{eqnarray*}
\| T_m x \|_{p} & \le & \| P_p T_m x \|_p + \| Q_p T_m x \|_{p} \\ & = & \| T_m P_p x \|_{L^{\circ}_p(\V)} + \| T_{m_{|_{G_0}}} Q_p x \|_{p} \ \lesssim_{(p)} \ \Big( \| T_m \|_{\B(L^{\circ}_p(\V))} + c \Big) \| x \|_p. \hskip5pt \square
\end{eqnarray*}

\demMax 
Assume that $B_t = \lambda( m \eta ( t \psi))$ has an $L_{(p/2)'}$-square-max decomposition. According to \eqref{columnSQF} with $\rho(z) = \eta(z) \varrho(z)$ for some $\varrho \in \mathcal{H}_0^{\infty}$, and using that $T_m$ commutes  with the spectral calculus of $A$ (the generator of $\S$) we obtain
\begin{eqnarray*}
\| T_m (x) \|_{H^c_p} & \sim_{(p)} & \big\| \big( \eta(tA) \varrho(tA) T_m x \big)_{t \geq 0} \big\|_{L_p(\V; L^c_2)} \\ & = & \big\| \big( \eta(tA) T_m \varrho(tA) x \big)_{t \geq 0} \big\|_{L_p(\V; L_2^c)} \\ & = & \big\| \big( T_{m_t} ( x_t ) \big)_{t \geq 0} \big\|_{L_p(\V; L^c_2)}, 
\end{eqnarray*}
where $m_t(g) = m(g) \eta(t \psi(g))$ and $x_t = T_{\varrho( t \psi)} x$. Recall also that the $L_2$-space involved is $L_2\left(\R_+, {d t}/t\right)$. Now we may express the term on the right hand side as follows   
\begin{eqnarray}
\hskip15pt \big\| \big( T_{m_t} ( x_t ) \big)_{t \geq 0} \big\|_{L_p(\V; L^c_2)}^2 & \!\!\! = \!\!\! & \Big\| \int_{\R_+} | T_{m_t} x_t |^2 \frac{dt}{t} \Big\|_{\frac{p}{2}} \label{Lp2Norm} \\ & \!\!\! = \!\!\! & \tau \Big( u \int_{\R_+} | T_{m_t} x_t |^2 \frac{dt}{t} \Big) = \int_{\R_+} \tau \big( u | T_{m_t} x_t |^2 \big) \frac{dt}{t}, \nonumber
\end{eqnarray}
where $u \in L_{( p / 2)'}(\V)_+$ is the unique element realizing the $L_{p / 2}$-norm, which exists by the weak-$\ast$ compactness of the unit ball of $L_{(p / 2)'}(\V)$. Now we have to estimate the term inside the integral. As $u \geq 0$, we may write $u = w^{\ast} w$ for some $w \in L_{2 ( p / 2)'}$ and $$\big\langle u, | T_{m_t} (x_t) |^2 \big\rangle = \tau (w | T_{m_t} (x_t) |^2 w^{\ast}) = \tau \Big( w \underbrace{\big| ( \tau \otimes \mathrm{Id}) \big( \delta B_t  (\sigma x_t \otimes \1) \big) \big|^2}_{L_t} w^{\ast} \Big).$$ As $L_t \mapsto w L_t w^{\ast}$ is order preserving, any bound of $L_t$ gives a bound of the above term. By the complete positivity of the canonical trace we can apply Proposition 1.1 in \cite{Lance1995}, i.e.
  \[
    \langle x, y \rangle^{\ast} \, \langle x, y \rangle \leq \|\langle x, x \rangle\| \langle y, y \rangle
  \]
to the operator-valued inner product $\langle x, y \rangle = ( \tau \otimes \Id) ( x^{\ast} y)$. This yields  
\begin{eqnarray*}
L_t & \!\! = \!\! & \big| ( \tau \otimes \mathrm{Id}) \big( \delta \Sigma_t \delta M_t  (\sigma x_t \otimes \1) \big) \big|^2  \\ [2pt]
      & \!\! \le \!\! & \left\| (\tau \otimes \Id) (\delta | \Sigma_t |^2 ) \right\|_{\V} \, ( \tau \otimes \Id) \Big( (\sigma x^{\ast}_t \otimes \1) \delta M_t^{\ast} \delta M_t  (\sigma x_t \otimes \1) \Big) \\
      & \!\! \le \!\! & \Big( \sup_{t > 0}  \| \Sigma_t \|^2_2 \Big) ( \tau \otimes \Id) \big( \delta | M_t |^2  ( \sigma ( x^{\ast}_t x_t) \otimes \1) \big) = \Big( \sup_{t > 0}  \| \Sigma_t \|^2_2 \Big) \big( | M_t |^2 \star x^{\ast}_t x_t \big).
\end{eqnarray*}
We have used the $\delta$-invariance of the trace in the second inequality and the definition of the noncommutative convolution in the last identity. Now, substituting inside the trace and using the identity for the adjoint of the noncommutative convolution operator gives $$\big\langle u, | T_{m_t} ( x_t) |^2 \big\rangle \le K^2 \tau \big( u ( | M_t |^2 \star x_t^{\ast} x_t ) \big) = K^2 \tau \big( ( \sigma | M_t |^2 \star u )x_t^{\ast} x_t \big),$$ where $K$ is the supremum of the $L_2$ norm of $\Sigma_t$. This gives rise to 
$$\begin{array}{rll}
\| T_m  ( x) \|^2_{H_p^c} & \lesssim_{(p)} & \displaystyle K^2 \, \int_{\R_+} \tau \big( (\sigma | M_t |^2 \star u) x_t^{\ast} x_t \big) \frac{dt}{t} \\
                                & \le & \displaystyle K^2 \inf_{\sigma | M_t^{} |^2 \star u \le A} \tau \Big( A \int_{\R_+} x_t^{\ast} x_t  \frac{dt}{t} \Big) \\
                                & \le & \displaystyle K^2 \inf_{\sigma | M_t |^2 \star u \le A}  \| A \|_{L_{(p/2)'}}  \Big\| \int_{\R_+} x_t^{\ast} x_t \frac{dt}{t} \Big\|_{p/2} \\ [10pt]
                                & \lesssim_{(p)} & \displaystyle K^2 \Big\| (R_t^c)_{t \geq 0} : L_{(p/2)'} \to L_{(p/2)'}(L_\infty) \Big\| \, \| x \|^2_{H_p^c} 
\end{array}$$  
by using Fatou's lemma in the second line and the definition of the $L_p(\V;L_\infty)$ norm for positive elements in the last inequality. Taking square roots gives the desired estimate. The calculations for the row case are entirely analogous. \fin

\begin{remark}
  \label{switch}
    \emph{
      Throughout this paper we construct max-square and square-max
      decompositions of $B_t = \lambda ( m \, \eta ( t \psi))$ by choosing an
      smoothing positive factor $M_t$ with $M_t = \sigma M_t = M_t^\ast$ and 
      satisfying the appropriate maximal inequalities. Then we extract $M_t$ from the
      left and from the right of $B_t$ as
      \begin{eqnarray*}
        B_t & = & ( B_t M_t^{- 1}) M_t, \\
        B_t & = & M_t ( M^{- 1}_t B_t) .
      \end{eqnarray*}
      If the family $\Sigma_t = B_t M_t^{- 1}$ is uniformly bounded
      in $L_2$ and $B_t$ is self-adjoint then the other is automatically 
      uniformly in $L_2$ by the traciality of $\tau$. Most of the times it will
      be enough to check one of the two decompositions.
  }
\end{remark}

\demA It easily follows from Corollary \ref{ColandRowMultiplierMax} and Remark \ref{switch}. \fin

\begin{remark} 
  \label{CBMultiplierMax}
  \emph{
    The technique employed here gives complete bounds assuming that the maximal
    inequalities are satisfied with complete bounds. In order to prove that 
    assertion, let us express the matrix extension $( T_m \otimes \Id_{M_n} )$
    as a matrix-valued multiplier whose symbol takes diagonal values. Indeed
    \[
      (T_m \otimes \Id_{M_n}) ([ x_{i j}]) = (\Id \otimes \tau \otimes \Id_{M_n}) \Big( \underbrace{\big( \delta \lambda(m) \otimes \1_{M_n} \big)}_{K} \big( \1 \otimes [\sigma x_{i j}] \big) \Big),
    \]
    where $K$ is the corresponding kernel affiliated with $\V \weaktensor \V \weaktensor \C \1_{M_n}$. Clearly, any square-max decomposition $B_t = \Sigma_t M_t$ of 
    $B_t = \lambda ( m \, \eta ( t \psi))$ yields a diagonal
    decomposition $(\delta \Sigma_t \otimes \1_{M_n}) (\delta M_t \otimes \1_{M_n})$ of 
    $K_t = \delta B_t \otimes \1_{M_n}$. On the other hand recall that $T_m : H_p^{c} \to H_p^{c}$ is c.b. iff 
    $T_m \otimes \Id_{M_n} : S_p^n[H_p^{c}] \to S_p^n[H_p^{c}]$ is uniformly 
    bounded for $n \ge 1$ and that
    $S_p^n[H_p^{c}(\V;\S)] = H_p^{c}( M_n \otimes \V; \Id \otimes \S)$. That
    allows us to write the norm of $S_p^n[H_p^{c}(\V;\S)]$ as an $L_{p/2}$-norm
    like in (\ref{Lp2Norm}). Then, using \cite[Proposition 1.1]{Lance1995} for  
    $\langle x, y \rangle = (\Id \otimes \tau \otimes \Id_{M_n})( x^* y )$
    as in the proof of Theorem \ref{BoundMultiplierMaximal}, gives for
    $2 < p < \infty$ $$\| T_m \|_{\CB(H^c_p)} \lesssim_{(p)} \Big( \sup_{t \geq 0} \| \Sigma_t \|_2 \Big) \, \Big\| (R_t^c)_{t \geq 0} : L_{(p/2)'}(\V) \to L_{(p/2)'}(\V; L_\infty) \Big\|_{\mathrm{cb}}^\frac{1}{2}.$$
The row case is similar. The discussion of 
    Corollary \ref{ColandRowMultiplierMax} generalizes to c.b. norms.
  }
\end{remark}

\section{{\bf Spectral multipliers \label{Sect2}}}

\subsection{{ \hskip3pt Ultracontractivity}}

Let $(\M,\tau)$ be a noncommutative measure space and consider a Markov 
semigroup $\S = (S_t)_{t \ge 0}$ defined on it. Given a positive function 
$\Phi: \R_+ \to \R_+$ and $1 \le p < q \le \infty$, we say that $\S$ 
satisfies the $\mathrm{R}_{\Phi}^{p,q}$ ultracontractivity property when
\[
  \begin{array}{>{\displaystyle}r>{\displaystyle}l>{\displaystyle}l>{\displaystyle}r}
    \big\| S_t: L_p(\M) \to L_q(\M) \big\| & \lesssim & \frac{1}{\Phi(\sqrt{t})^{\frac1p - \frac1q}} & \forall \, t > 0.
    \label{ultracontractivity}
    \tag{$\mathrm{R}_{\Phi}^{p,q}$}
  \end{array}
\]
Similarly, $\S$ has the 
$\mathrm{CBR}_{\Phi}^{p,q}$ property when the above estimate 
holds for the c.b. norm of $S_t: L_p(\M) \to L_q(\M)$. These inequalities have 
been extensively studied for commutative measure spaces 
\cite[Chapter 1]{VaSaCou1992}. In the theory of Lie groups with an invariant 
Riemannian metric (equipped with the heat semigroup generated by the
invariant Laplacian) ultracontractivity holds for the function 
$\Phi(t) = \mu(B_t(e))$ which assigns the volume of a ball for a 
given radius. Influenced by that, we will interpret the above-defined 
properties as a way of describing the \lq\lq growth of the balls" in the 
noncommutative geometry determined by $\S = (S_t)_{t \ge 0}$. For that reason, 
we will work with \emph{doubling functions} $\Phi$. Doubling functions are
increasing functions $\Phi: \R_+ \to \R_+$ with $\Phi(0) = 0$ and 
satisfying
\[
  \sup_{t > 0} \left\{ \frac{\Phi(2t)}{\Phi(t)} \right\} \, < \, \infty.
\]
The doubling condition for $\Phi$ is a natural requirement since metric measure 
spaces $(\Omega, \mu, d)$ with $\Phi_x(t) = \mu(B_x(t))$ uniformly 
doubling in $x$ constitute an adequate setting for performing harmonic analysis 
in commutative measure spaces. Given a Markov semigroup $\S = (S_t)_{t \ge 0}$ 
over a noncommutative measure space $(\M,\tau)$, let us recall the following:
\begin{itemize}
  \item[i)] If $\S$ satisfies $\mathrm{R}_{\Phi}^{p_0,q_0}$, it satisfies 
  $\mathrm{R}^{p,q}_{\Phi}$ for $1 \le p_0 \le p < q \le q_0 \le \infty$.
  \item[ii)] If $\Phi$ is doubling and $\S$ satisfies 
  $\mathrm{R}^{p_0,q_0}_{\Phi}$ for some $1 \le p_0 < q_0 \le \infty$, then it satisfies 
  $\mathrm{R}^{p,q}_\Phi$ for $1 \le p \leq q \le \infty$ 
\end{itemize}
The same holds for the $\mathrm{CBR}^{p_0,q_0}_\Phi$ ultracontractivity property. The proof follows the 
same lines than \cite[Theorem II.1.3]{VaSaCou1992}. In the noncommutative 
setting a similar result is stated in  \cite[Lemma 1.1.2 ]{JunMe2010}
for $\Phi(t)=t^D$. As a consequence, all the ultracontractivity properties 
$\mathrm{R}^{p,q}_\Phi$ are equivalent for doubling $\Phi$. We shall denote 
them simply by $\mathrm{R}_\Phi$ and similarly $\mathrm{CBR}_\Phi$. As a 
corollary, we obtain that if $\M$ is an abelian von Neumann algebra 
$\mbox{CBR}_{\Phi}^{p,q}$ and  $\mbox{R}_{\Phi}^{p,q}$ are equivalent for 
doubling $\Phi$ since $\mbox{R}_{\Phi}^{p,q}$ is equivalent to 
$\mbox{R}_{\Phi}^{p,\infty}$ and any bounded map into an abelian 
$C^{\ast}$-algebra is completely bounded. For any doubling function
$\Phi$ we may define its doubling dimension $D_\Phi$ as
\[
  D_{\Phi} = \log_2  \sup_{t > 0} \left\{ \frac{\Phi(2t)}{\Phi(t)} \right\}.
\]
It is quite simple to show that any doubling $\Phi: \R_+ \to \R_+$ admits 
upper/lower polynomial bounds for large/small values of $t>0$. More precisely, 
we have the bounds
\begin{equation} \label{eq:lower}
\begin{array}{rcl}
    \Phi(t) & \lesssim_{(D_\Phi)} & t^{D_\Phi} \, \Phi(1) \quad \mbox{when} \quad t>1, \\
    \Phi(t) & \gtrsim_{(D_\Phi)}  & t^{D_\Phi} \, \Phi(1) \quad \mbox{when} \quad t \le 1. 
\end{array}
\end{equation}
Of course, the converse of this assertion is false. Whenever a Markovian 
semigroup $\S$ satisfies $\mbox{R}_\Phi$ (resp. $\mbox{CBR}_\Phi$) for 
doubling $\Phi$ we will call $D_\Phi$ the \emph{Sobolev dimension} (resp. 
\emph{c.b. Sobolev dimension}) of $(\M,\tau)$ with respect to $\S$. The 
reason for this name is based on the well-known relation between 
ultracontractivity estimates for a Markov semigroup and Sobolev embedding 
estimates for its infinitesimal generator. One of the first contributions 
to that relation is in the work of Varopoulos, who proved in \cite{Va1985} 
that when $\Phi(t)=t^D$ the property $\mathrm{R}_\Phi$ is equivalent to a whole 
range of Sobolev type estimates for the infinitesimal generator of the 
semigroup. See also \cite{VaSaCou1992} for more on that topic. Whenever 
$\Phi(t) = t^D$ we will denote the ultracontractivity properties by $\mbox{R}_D$ or $\mbox{CBR}_D$. By adding a zero, like 
$\mbox{R}_{\Phi}(0)$, we will mean that the inequality 
\ref{ultracontractivity} is satisfied for $t \leq 1$. This notation
is borrowed from \cite[II.5]{VaSaCou1992}. Recall that if $\S$ satisfies 
$\mathrm{R}_\Phi$ (resp. $\mbox{CBR}_\Phi$) for some doubling function 
$\Phi$ then, by the polynomial bounds in (\ref{eq:lower}), we have 
$\mbox{R}_{D_\Phi}(0)$ (resp. $\mbox{CBR}_{D_\Phi}(0)$).

Our characterization of co-polynomial growth in Section
\ref{Sect3} bellow requires the following equivalence for Sobolev-type
inequalities in term of the ultracontractivity properties $\mathrm{R}_D(0)$.
We did not find the proposition below in the literature, but it could be well-known to experts. We include a sketch of the proof. 

\begin{proposition}
\label{RdeqSob}
Let $\S$ be a submarkovian semigroup acting on a noncommutative measure space $(\M,\tau)$. Let $A$ denote its infinitesimal generator. Then, the following properties are equivalent$\, :$
  \begin{itemize}
    \item[i)] For every $\varepsilon > 0$, $\S$ satisfies the 
    $\mathrm{R}_{D+\varepsilon}(0)$ property.
    \item[ii)] For every $\varepsilon > 0$, we have that
    \[
      \| ( \1 + A )^{-D/4-\varepsilon} : L_2(\M) \rightarrow \M\| \lesssim_{(\varepsilon)}1.
    \]
  \end{itemize}
Similarly, $\S \in \mathrm{CBR}_{D+\varepsilon}(0)$ for all $\varepsilon > 0$ iff $(\1+A)^{-s} : L_2(\M) \stackrel{\mathrm{cb}}{\longrightarrow} \M$ for all $\varepsilon > 0$.
\end{proposition}

\dem
The implication $\mathrm{i}) \Rightarrow \mathrm{ii})$ follows from the identity
\begin{eqnarray*}
(\1 + A)^{- s} (x) & = & \frac{1}{\Gamma (s)} \Big( \int_{\R_+} t^{s} e^{- t} S_t (x) \frac{dt}{t} \Big).
\end{eqnarray*}
The integral in $[0,1]$ may be estimated applying the $\mathrm{R}_{D}(0)$ property, whereas the integral for $t>1$ is easily estimated using the semigroup law. This gives the desired implication. For the converse, we now take $s = D / 4 + \varepsilon$ and
use that $\|f(A)\|_{\B(L_2)} \le \|f\|_\infty$
\begin{eqnarray*}
\lefteqn{\hskip-20pt \big\| S_t : L_2(\M) \to \M \big\| = \left\| ( \1 + A)^{-\frac{s}{2}} ( \1 + A)^{\frac{s}{2}} S_t \right\|_{\B(L_2(\M),\M)}} \\ \hskip25pt & \!\!\! \le \!\!\! & \left\| ( \1 + A)^{-\frac{s}{2}} \right\|_{\B(L_2(\M),\M)}  \left\| ( \1 + A )^{\frac{s}{2}} S_t \right\|_{\B(L_2)}  \lesssim_{(\varepsilon,s)} \left( \frac{s}{2} \right)^{\left( \frac{s}{2} \right)} e^{- \frac{s}{2}} \frac{e^t}{t^{\frac{s}{2}}}. \hskip15pt \square
\end{eqnarray*}

\begin{remark}
  \emph{
    Observe that if $\mbox{R}_{D}(0)$ is satisfied then ii) also 
    holds. Nevertheless the converse is not true since the norm 
    $\| S_t : L_1(\M) \rightarrow \M \|$ could be comparable to, say, 
    $t^D (1+\log(t))$ for $0 \leq t \leq 1$. The original result proved by 
    Varopoulos \cite{Va1985} established a equivalence 
    between $\mbox{R}_D(0)$ and the bounds 
    $$(\1+A)^{-s} : L_p(\M) \rightarrow L_{\frac{p n}{n - s p}}(\M)$$ for every 
    $0 \leq s < n/p$. When $s > n/p$ the image space of $L_p(\M)$ is certainly 
    much smaller than $L_\infty(\M)$, for example in $\R^n$ with the usual 
    Laplacian the image space lies inside spaces of H\"older functions. Therefore,
    by describing the behavior of $(\1+A)^{-s}$ in $L_\infty(\M)$ we lose 
    information and we can no longer recover $\mbox{R}_D(0)$.
  }
\end{remark}

We will denote by $W_{A}^{p,s}(\M)$, or simply 
$W^{p,s}(\M)$ when the semigroup $S_t = e^{- t A}$ can be understood 
from the context, the closed domain in $L_p(\M)$ of the unbounded 
operator $(\1 + A)^{s/2}$, with norm given by
\[
	\| x \|_{W_A^{p,s}} = \big\| ( \1 + A )^{s/2} f \big\|_{p}.
\]
These are called the \emph{fractional Sobolev spaces} associated with 
$\S$. They satisfy the natural interpolation identities. Namely, if we set $1/p_3 = (1-\theta)/p_1 + \theta/p_2$ we get
\begin{eqnarray*}
\big[ W_{A}^{p_1,s}(\M), W_{A}^{p_2,s}(\M) \big]_\theta & \simeq & W_{A}^{p_3,s}(\M), \\
\big[ W_{A}^{p,s_1}(\M), W_{A}^{p,s_2}(\M) \big]_\theta & \simeq & W_{A}^{p,s_1 \theta + s_2 (1-\theta)}(\M), 
\end{eqnarray*}
Point ii) in Proposition \ref{RdeqSob} may be rephrased as $W_{A}^{2,s}(\M) \subset \M$ for every $s > D/2$.

\subsubsection{{\bf \hskip3pt $L_2$ bounds for $\CB(L_2(\V),\V)$ multipliers.}}

We shall work extensively with Markovian convolution semigroups over $\V$ with the $\mathrm{CBR}_\Phi$ ultracontractivity property for doubling $\Phi$. In general, determining the c.b. norm of a multiplier between 
general $L_p$ spaces is a problem that nobody expects to be solvable with a 
closed formula. Despite that, we can obtain characterizations in some 
particular cases. One of these cases is that of the c.b. multipliers 
$T_m : L_2(\V) \to \V$. That will allow us to express the
$\mathrm{CBR}^{2,\infty}_\Phi$ property of $\S = (T_{e^{- t \psi}})_{t \geq 0}$
as a condition over $\psi$. The next theorem is probably known to experts.
Since we could not find it in the literature, we include it here for the sake of
completeness.
\begin{theorem}
\label{CBL2Multiplier}
If $T$ denotes the map $m \mapsto T_m \hskip-1pt :$
\begin{itemize}
    \item[i)]  $T: L^r_2(G) \to \CB(L^{c}_2(\V), \V)$ is a complete isometry.
    \item[ii)] $T: L^c_2(G) \to \CB(L^{r}_2(\V), \V)$ is a complete isometry.
\end{itemize}
The image of $T$ is the set of multipliers $T_m : L_2^\dag(\V) \stackrel{\mathrm{cb}}{\longrightarrow} \V$ for $\dag \in \{c,r\}$ resp.
\end{theorem}

\begin{proof}
Let $V$ and $W$ be operator spaces and pick $x \otimes y \in V^* \otimes W^*$. According to \cite[Theorem 4.1]{Pi2003} the map $\mathcal{I}_{x \otimes y}(w) = x \langle y , w \rangle$ extends linearly to an isomorphism $\mathcal{I} : (V \osprojtensor W)^* \to \CB(W,V^*)$. Using the pairying $\langle \, , \rangle : L_2^r(\V) \times L_2^c(\V) \to \C$ given by $\langle y,w \rangle = \tau(y \, \sigma w)$ we obtain as a consequence that $$\mathcal{I}_{\delta z}(w) = ( \Id \otimes \tau ) \big( \delta z \, (\1 \otimes \sigma w) \big) = z \star w,$$ where $\delta z$ denotes the comultiplication map acting on $z$. This yields $$\big\| T_m: L_2^\dag (\V) \to \V \big\|_{\mathrm{cb}} = \big\| \delta \lambda(m) \big\|_{(\V_* \widehat{\otimes} L_2^\dag(\V))^*}$$ where $\dag \in \{r,c\}$ is either the row or the column o.s.s. We now claim that the natural map $$\iota : L_\infty(\V;L_2^{\dag^{\op}}(\V)) \hookrightarrow \big(\V_{\ast} \osprojtensor L_2^{\dag}(\V)\big)^*$$ is a complete isometry with $\dag^{\op} = r$ for $\dag = c$ and viceversa. This is all what is needed to complete the argument since we have the following commutative diagram of complete isometries
\begin{equation*}
  \xymatrix{
    L_2^{\dagger}(G) \ar[rr]^{T} \ar[d]^{\lambda} & & \CB(L_2^{\dagger^\op}(\V), \V)\\
    L_2^{\dagger}(\V) \ar[dr]^{\delta} & & (\V_\ast \osprojtensor L^{\dagger^\op}_2(\V))^* \ar[u]^{\mathcal{I}}\\
    & L_\infty(\V;L_2^{\dagger}(\V)) \ar@{^{(}->}[ur]^{\iota} &
  }.
\end{equation*}  
Let us therefore justify our claim. According to \cite{EfRu2000} $$(\V_{\ast} \osprojtensor L_2^{\dag}(\V))^* \simeq \V \Fubinitensor L_2^{\dag^\op}(\V)$$ where $\Fubinitensor$ stands for the Fubini tensor product of dual operator spaces. Bear in mind that
  if $V^*$ and $W^*$ are dual operator spaces, there are weak-$\ast$ continuous
  embeddings $V^* \subset \B(\H_1)$ and $W^* \subset \B(\H_2)$ and we can 
  define the weak-$*$ spatial tensor product $V^* \weaktensor W^*$ as
  \[
    V^* \weaktensor W^* = \overline{(V^* \otimes W^*)^{\mathrm{w}^\ast}}.
  \] 
  Such construction is representation independent and
  $V^* \weaktensor W^*$ embeds completely isometricaly in 
  $V^* \Fubinitensor W^*$. Since the column and row
  embeddings of $L_2(\V)$ into $\B(L_2(\V))$ are weak-$\ast$ continuous,
  $L_\infty(\V;L_2^{\dag^\op}(\V)) = \V  \weaktensor L_2^{\dag^\op}(\V)$. This proves that $\iota$ is a complete isometry and so is the map $m \mapsto T_m = \mathcal{I}_{\iota \delta \lambda(m)}$.  
  \end{proof}

\begin{remark}
  \label{CBL2Norm}
  \emph{
    Since $L^r_2(\V)$ and $L^c_2(\V)$ are isometric as Banach
    spaces, the norms for multipliers in $\CB(L^{r}_2(\V), \V)$ and
    $\CB(L^{c}_2(\V), \V)$ coincide too, even if their matrix amplifications
    do not. Indeed we obtain that $$\| T_m \|_{\CB(L_2^r(\V), \V)} = \| m \|_{L_2(G)} = \| T_m \|_{\CB(L_2^c(\V), \V)}.$$
    For non-hyperfinite $\V$, the space of Fourier multipliers in 
    $\CB(L_2(\V), \V)$, may be difficult to describe as an operator space. 
    Nevertheless, as a consequence of the above identities, its underlying
    Banach space is the Hilbert space $L_2(G)$.
  }
\end{remark}

\begin{remark}
  \emph{
    \label{CBRimpAmenable}
    As a consequence of the above, if $G$ is a group and
    $\S = (T_{e^{- t \, \psi}})_{t \geq 0}$ is a semigroup of Fourier
    multipliers satisfying $\mathrm{CBR}_{\Phi}^{2,\infty}$ for any
    function $\Phi$, then $G$ is amenable. To see it just notice that
    $e^{-t \psi} \in L_2(G)$ and so $e^{- 2 t \psi} \in L_1(G)$ for all $t > 0$. But
    a group is amenable iff there is a
    sequence of integrable positive type functions converging to $1$
    uniformly in compacts.  
  }
\end{remark}

\subsection{ \hskip3pt Standard assumptions}

Let $\V_+^\wedge$ denote the extended positive cone of $\V$. As it will become clear along the paper, we shall treat unbounded operators $X$ in $\V_+^{\wedge}$ as noncommutative or quantized metrics over $\V$. Note that  if $G$ is LCH and abelian, any translation-invariant metric over its dual group can be associated with the positive function $\Delta: \chi \mapsto d(\chi, e)$. The metric conditions impose that $\Delta$ is symmetric, does not vanish outside $e$ and $\Delta(\chi_1 \chi_2) \le \Delta(\chi_1) + \Delta(\chi_2)$. Here we will only require $X$ to be symmetric, i.e.: to satisfy $\sigma X = X$. Recall that the anti-automorphism $\sigma$ extends to $\V_+^\wedge$. Following the intuition relating symmetric operators in $\V_+^{\wedge}$ to metrics, we will say that $X \in \V_+^{\wedge}$ is doubling iff the function $\Phi_X(r) = \tau(\chi_{[0,r)}(X))$ is doubling. When the dependency on the operator $X$ can be understood from the context we will just write $\Phi$. In a similar fashion, we will say that $X$ satisfies the $L_p$-Hardy-Littlewood maximal property when
\begin{equation}
  \Big\| {\sup_{r \geq 0}}^{+} \Big\{ \frac{\chi_{[0,r)}(X)}{\Phi_X(r)} \star u \Big\} \Big\|_p \lesssim \| u \|_p,
  \label{HL}
  \tag{$\mathrm{H L}_p$}
\end{equation}
If we say that $X$ has the $\mathrm{HL}$ property, omitting the dependency on
$p$, we mean that the HL property is satisfied for every $1 < p \leq \infty$,
with constants depending on $p$. When the property \ref{HL} holds uniformly 
for all matrix amplifications, we will say that $X$ satisfies the \emph{completely bounded
Hardy-Littlewood maximal property} ($\mbox{CBHL}_p$ in short).
Let $\psi: G \to \R_+$ be a conditionally negative length generating a 
semigroup $\S$. We will say that $\S$ has $L_2$ Gaussian bounds with respect to
$X$ when there is some $\beta > 0$ such that
\begin{equation}
  \tau \Big\{ \chi_{[r,\infty)}(X) | \lambda( e^{-t \psi} ) |^2 \Big\}
  \lesssim
  \frac{e^{-\beta \frac{r^2}{t}}}{\Phi_X(\sqrt{t})}.
  \label{L2GB}
  \tag{$L_2\mbox{GB}$}
\end{equation}

\begin{definition}
  A triple $(\V, \S, X)$, where $\S$ is a Markov semigroup of Fourier
  multipliers generated by $\psi : G \to \R_+$ and $X \in (\V)_+^{\wedge}$,
  is said to satisfy the \emph{standard assumptions} when  
  \begin{itemize}
    \item[i)]   $X$ is symmetric and doubling.
    \item[ii)]  $\S$ has $L_2\mathrm{GB}$ with respect to $X$.
    \item[iii)] $X$ satisfies the $\mathrm{CBHL}$ property.
  \end{itemize}
  Since $\V$ is determined by $G$ and $\S$ by $\psi$ we shall often write 
  $(G,\psi, X)$ instead.
\end{definition}

\begin{remark}
\emph{If $\S$ has $L_2\mathrm{GB}$ then it admits $\mbox{CBR}_{\Phi_X}^{2,\infty}$ ultracontractivity. Namely if we take $r = 0$ in ($L_2\mathrm{GB}$), it follows from Theorem \ref{CBL2Multiplier} and Remark \ref{CBL2Norm}. 
If $X$ is in addition doubling, $\S$ has the whole range of ultracontractivity properties $\mbox{CBR}_{\Phi_X}$.}
\end{remark}

\subsubsection{{\bf Stability under Cartesian products.}}

It is interesting to note that the standard assumptions are stable under certain algebraic operations, the most trivial of them is probably the Cartesian product. Stability under crossed products also holds under natural conditions, see Remark \ref{RemCrossed} below. 

\begin{lemma}
\label{maxCartesian}
Assume that $$\S^j = (S^{j}_{\omega_j})_{\omega_j \in \Omega_j} : L_p(\M_j) \to L_p(\M_j; L_\infty(\Omega_j))$$ is completely positive for $j \in \{1,2\}$. Then $\S^{1} \otimes \S^{2}$ is also c.p. and
\begin{eqnarray*}
\lefteqn{\hskip-108pt \Big\| \S^{1} \otimes \S^{2} : L_p(\M_1 \weaktensor \M_2) \to L_p \big( \M_1 \weaktensor \M_2; L_\infty(\Omega_1) \otimes_{\min} L_\infty(\Omega_2) \big) \Big\|_{\mathrm{cb}}} \\
\hskip105pt & \lesssim & \prod_{j \in \{1,2\}} \Big\| \S^{j} : L_p(\M_j) \to L_p(\M_j; L_\infty(\Omega_j)) \Big\|_{\mathrm{cb}}. 
\end{eqnarray*}
\end{lemma}

\begin{proof}
It follows from $\S^1 \otimes \S^2 = (\S^1 \otimes \Id) \circ \S^2$ and \eqref{PostiveMax}, details are omitted.
\end{proof}

\begin{theorem}
\label{stabilityCartesian}
Let $(G_j, \psi_j, X_j)$ be triples satisfying the standard assumptions for $j=1,2$ and consider the Cartesian product $G = G_1 \times G_2$ . Then $(G, \phi, X)$ also satisfies the standard assumptions with the c.n.length $\psi(g_1,g_2) = \psi_1(g_1) + \psi_2(g_2)$ and $X \in \V_+^\wedge$ determined by the formula $X^2 = X^2_1 \otimes \1 + \1 \otimes X_2^2$.
\end{theorem}

\begin{proof}
Proving that $X$ is doubling and that the semigroup generated by $\psi$ has Gaussian bounds amount to a trivial calculation. Indeed, $\Phi_X$ is controlled from the inequalities $\chi_{\left[0,r/2   \right)}(a) \chi_{\left[0,r/2\right)}(b) \leq \chi_{[0,r)}(a + b) \leq \chi_{[0, r)}(a) \, \chi_{[0,r)}(b)$, which are valid for positive and commuting operators $a,b$. On the other hand, the $L_2\mathrm{GB}$ follow similarly from the inequality $\chi_{[r,\infty)}(a+b) \le \chi_{[r/2,\infty)}(a) + \chi_{[r/2,\infty)}(b)$. Let us now justify the $\mbox{CBHL}$ property. Let $m : L_p(\V; L_\infty \mintensor L_\infty) \to L_p(\V; L_\infty)$ be the map given by $m(x \otimes f \otimes g) = x \otimes f g$, which is c.p. By Lemma \ref{maxCartesian} 
  \[
    \mathcal{R}^1 \otimes \mathcal{R}^2 = \big( R^1_s \otimes R^2_t \big)_{s,t \geq 0} : L_p(\V) \to L_p \big( \V;L_\infty(d s) \mintensor L_\infty(d t) \big),
  \]
where $R^j_s(x) = \Phi_{X_j}(s)^{-1} \, \chi_{[0,s)}(X_j) \star x$ is c.p. As a consequence $m \circ (\mathcal{R}^1 \otimes \mathcal{R}^2)$ is also completely positive. Therefore, by the doubling property we obtain the following estimate
  \[
    \Big( \frac{\chi_{[0,r)}(X)}{\Phi_X(r)} \star x \Big)_{r \ge 0} \lesssim_{(D_{\Phi_1}, D_{\Phi_2})} \Big( \frac{\chi_{[0,r)}(X_1)}{\Phi_{X_1}(r)} \otimes \frac{\chi_{[0,r)}(X_2)}{\Phi_{X_2}(r)} \star x \Big)_{r \ge 0 }= m \circ (\mathcal{R}^1 \otimes \mathcal{R}^2) (x)
  \]
for $x \ge 0$. This is all what we need to reduce CBHL of $X$ to that of $X_1$ and $X_2$.
\end{proof}

\begin{remark} \label{RemCrossed}
\emph{Let $H$ and $G$ be LHC unimodular groups and $\theta : G \to \Aut(H)$ be a measure preserving action. Let $(H, \psi_1, X_1)$ and $(G, \psi_2, X_2)$ be triples satisfying the standard assumptions. It is possible to prove that, under certain invariance conditions on $X_1$ and $\psi_1$, the semidirect product $K = H \rtimes_\theta G$ satisfies the standard assumptions for some $X \in \mathcal{L}K_*^\wedge$ and certain c.n. length function $\psi: K \to \R_+$ built up from $X_1,X_2$ and $\psi_1,\psi_2$ respectively. Since the techniques required to prove this result are quite involved and of independent interest, we postpone its proof to a forthcoming paper were we shall explore other applications involving Bochner-Riesz summability and related topics.}
\end{remark}

\subsection{ \hskip3pt H\"ormander-Mikhlin criteria}
In this subsection we shall give a proof of Theorem B i) by means of a suitably chosen max-square decomposition. The key is to prove that, if $B_t = \lambda(m \, \eta(t \psi))$, then
\begin{eqnarray}
  B_t & = & \underbrace{B_t \, \Big(\1 + \frac{X^2}{t} \Big)^{\frac{\gamma}{2}} {\Phi(\sqrt{t})}^{\frac{1}{2}}}_{\Sigma_t} \underbrace{{\Phi(\sqrt{t})}^{-\frac{1}{2}} \Big(\1 + \frac{X^2}{t}\Big)^{-\frac{\gamma}{2}}}_{M_t}.
  \label{SMdec}
\end{eqnarray}
is a square-max decomposition for $\gamma > D_\Phi/2$. Breaking the symbol $m$ into its real and imaginary parts and using Remark \ref{switch}, we obtain a max-square decomposition by placing the smoothing factor $(\1 + {X^2}/{t})^{\gamma/2}$ on the left hand side of $B_t$. The proof of the maximal inequality consists in expressing the maximal operator as a linear combination of Hardy-Littlewood maximal operators associated to $X$ and apply \eqref{averageMax}. For the square estimate we will use the smoothness condition.

\begin{lemma}
  \label{MaxLemma}
  Assume that $F_t \in C_0(\R_+)$ is
  a family of bounded variation functions parametrized by $t > 0$. Let $dF_t$
  be its Lebesgue-Stjeltjies derivative and $| dF_t ( \lambda) |$ its absolute
  variation, then for every doubling operator $X$, we have:
  \begin{eqnarray*}
    \Big\| \Big( {\sup_{t \geqslant 0}}^{+} F_t (X) \star x \Big)  \Big\|_{L_p} & \le &
    \Big( \sup_{t > 0} \| \Phi \|_{L_1 ( | dF_t |)} \Big) \Big\|  \Big( {\sup_{r > 0}}^+
    \frac{\chi_{[0, r)} ( X)}{\Phi ( r)} \star x \Big)  \Big\|_{L_p}
  \end{eqnarray*}
\end{lemma}

\begin{proof} By integration by parts we have that
  \begin{eqnarray*}
    F_t(s)  & = & \int_{\R_+} F_t ( r) d \delta_{s}(r) \ = \ \int_{\R_+} F_t ( r) \partial \chi_{( s, \infty)}(r)\\
                  & = & - \int_{\R_+} \chi_{( s, \infty)} ( r) \partial F_t(r) \ = \ - \int_{\R_+} \frac{\chi_{[ 0, r)} ( s)}{\Phi ( r)} \Phi(r) \partial F_t(r) .
  \end{eqnarray*}
  By functional calculus, the same holds for $F_t(X)$.
  Applying \eqref{averageMax} ends the proof.
\end{proof}

According to Theorem A, the right choice for the square-max decomposition is given by $F_t(s) = |M_t|^2(s) = \Phi(\sqrt{t})^{-1} (1 + s^2 / t)^{-\gamma}$. It will suffice to pick here $\gamma > D_\Phi/2$, the condition in Theorem B i) will be justified later on. In order to prove the finiteness of the maximal bound in Theorem A, we just need to verify the condition of Lemma \ref{MaxLemma} for this concrete function. 

\begin{lemma}
  \label{MaxCon}
  For any doubling $\Phi$, we find
  \[
    \int_{\R_+} \Phi (s)  \Big| \frac{d}{ds}  \Big( 1 + \frac{s^2}{t} \Big)^{- \frac{D_{\Phi} + \varepsilon}{2}}  \Big| \, ds
    \, \lesssim_{(D_{\Phi}, \varepsilon)} \,
    \Phi(\sqrt{t}).
  \]
\end{lemma}

\begin{proof}
Changing variables $s \mapsto \sqrt{tv}$, we obtain
\begin{eqnarray*}
\int_{\R_+} \Phi (s)  \Big| \frac{d}{ds}  \Big( 1 + \frac{s^2}{t} \Big)^{- \frac{D_{\Phi} + \varepsilon}{2}}  \Big| \, ds & \sim_{(D_\Phi)} & \int_{\R_+} \Phi(s)  \Big( 1 + \frac{s^2}{t} \Big)^{-\frac{D_{\Phi} + 2 + \varepsilon}{2}} \frac{2 s}{t} ds \\ & = & \int_{\R_+} \Phi( \sqrt{t}  \sqrt{v}) \big(1 + v\big)^{- \frac{D_{\Phi} + 2 + \varepsilon}{2}} dv \\ & = & \Big( \int_0^1 + \sum_{k = 0}^{\infty} \int_{4^k}^{4^{k + 1}} \Big) \ = \ A + \sum_{k = 0}^{\infty} B_k .
\end{eqnarray*}
The monotonicity of $\Phi$ gives $A \le \Phi(\sqrt{t})$, while its doublingness yields $$\begin{array}{rll} B_k & \le & \displaystyle \Phi(\sqrt{t}) \, 2^{D_{\Phi} (k + 1)} \int_{4^k}^{4^{k + 1}} \big( 1 + v \big)^{- \frac{D_{\Phi} + 2 + \varepsilon}{2}} dv \\ & \sim_{(D_{\Phi})} & \Phi(\sqrt{t}) \, 2^{D_{\Phi} (k + 1)} 2^{- ( D_{\Phi} + \varepsilon) k} \\ [6pt] & \sim_{(D_{\Phi})} & \Phi(\sqrt{t}) \, 2^{- \varepsilon k}. \end{array}$$  
Since the sequence of $B_k$s is summable, we have proved the desired estimate.
\end{proof}

For the estimate of the square part, let us start by extending the Gaussian bounds to the complex half-plane $\mathbb{H} = \{z \in \C : \mathrm{Re}(z)>0 \}$. We need the following version of the Phragmen-Lindel\"off theorem, see \cite{Davies1995141} for the proof.

\begin{theorem}
If $F$ is analytic over $\mathbb{H}$ and satisfies
  \begin{eqnarray*}
    |F(|z| e^{i \theta})| & \lesssim & \big( |z| \, \cos \theta \big)^{- \beta}, \\
    |F(|z|)|              & \lesssim & |z|^{-\beta} \exp \big(- \alpha |z|^{-\rho} \big),
  \end{eqnarray*}
  for some $\alpha, \beta >0$ and $0 < \rho \leq 1$, then we find the following estimate $$|F(|z| e^{i \theta})| \lesssim_{(\beta)} \big( |z| \, \cos \theta \big)^{- \beta} \exp \Big( - \frac{\alpha \rho}{2} |z|^{-\rho} \cos \theta \Big).$$
\end{theorem}

\noindent We may now generalize the Gaussian $L_2$-bounds to the complex half-plane.

\begin{proposition} \label{z-Gauss}
Let $G$ be a unimodular group, $\psi : G \to \R_+$ a c.n. length and
$X \in \V_+^{\wedge}$ a doubling operator satisfying $L_2 \mathrm{GB}$. If we set
  $h_z = \lambda(e^{- z \psi})$, the following bound holds for every 
  $z \in \mathbb{H}$
  $$\tau \Big\{ \chi_{[r, \infty)} (X) | h_z |^2 \Big\} \, \lesssim \, \frac{1}{\Phi(\sqrt{\mathrm{Re} \{ z \}})} e^{- \frac{\beta}{2} \frac{r^2}{| z |} \frac{\mathrm{Re} \{ z \}}{| z |}}.$$
 \end{proposition}

\begin{proof}
Let $x$ be an element of $L_2(\V)$ with $\| x \|_2 \le 1$. Assume in addition that $x = px$ for $p = \chi_{[ r,\infty)}(X)$. Then we define $G_x$ as the following
  holomorphic function
  \[
    G_x (z) = e^{- \frac{z}{t}} \Phi(\sqrt{t}) \tau (h_z x)^2.
  \]
  Then, the estimate below holds in $\mathbb{H}$
\begin{eqnarray*}
| G_x (z) | & = & e^{- \frac{\mathrm{Re} \{ z \}}{t}} \Phi(\sqrt{t}) | \tau(h_z x) |^2 \ \le e^{- \frac{\mathrm{Re} \{ z \}}{t}} \Phi(\sqrt{t}) \tau(|h_z|^2) \\ & = & e^{- \frac{| z | \cos \theta}{t}} \Phi(\sqrt{t}) \| h_{\mathrm{Re} \{ z \}} \|_{L_2(\V)}^2 \ \lesssim \ e^{- \frac{\mathrm{Re} \{ z \}}{t}} \Phi(\sqrt{t}) / \Phi(\sqrt{\mathrm{Re} \{ z \}}).
\end{eqnarray*}
Note that the second identity above follows from Plancherel theorem and the last inequality from $L_2$GB for $r=0$. On the other hand, since $\Phi$ is doubling it satisfies $\Phi(s (1 + r)) \lesssim \Phi(s) (1 + r)^{D_\Phi}$ for every $r > 0$ and $$|G_x(z) | \lesssim e^{- \frac{\mathrm{Re} \{ z \}}{t}} \frac{\Phi(\sqrt{t})}{\Phi(\sqrt{\mathrm{Re} \{ z \}})} \lesssim e^{- \frac{\mathrm{Re} \{ z \}}{t}} \Big( 1 + \frac{\sqrt{t}}{\sqrt{\mathrm{Re} \{ z\}}} \Big)^{D_{\Phi}} \lesssim \Big( \frac{t}{| z | \cos \theta} \Big)^{\frac{D_{\Phi}}{2}},$$ by using that $e^{- s^2}  ( 1 + 1 / s)^a \lesssim (1/s)^a$ in the last inequality. We also have
\begin{eqnarray*}
| G_x(| z |) | & = & e^{- \frac{| z |}{t}} \Phi(\sqrt{t}) \big| \tau(h_{| z |} x) \big|^2 \\ & \le & e^{- \frac{| z |}{t}} \Phi(\sqrt{t}) \, \tau \big\{ p \, h_{| z |}^* h_{| z |} \big\} \ \lesssim \ e^{- \frac{| z |}{t}} \frac{\Phi(\sqrt{t})}{\Phi(\sqrt{| z |})} e^{- \beta \frac{r^2}{| z |}} \\ & \lesssim & e^{- \frac{| z |}{t}} \Big( 1 + \frac{\sqrt{t}}{\sqrt{| z |}} \Big)^{D_{\Phi}} e^{- \beta \frac{r^2}{| z |}} \ \lesssim \ \Big( \frac{t}{| z |} \Big)^{\frac{D_{\Phi}}{2}} e^{-\beta \frac{r^2}{| z |}} .
\end{eqnarray*}
  The Phargmen-Lindel{\"o}f theorem allows us to
  combine both estimates, giving
  \begin{equation*}
    | G_x(| z | e^{i \theta}) | \, \lesssim \, t^{\frac{D_{\Phi}}{2}} \big( |z| \cos \theta \big)^{- \frac{D_{\Phi}}{2} } e^{- \frac{\beta r^2}{2}  \frac{\cos\theta}{|z|}}.
  \end{equation*}
  Taking the supremum over all $x$ with $\|x\|_2 \leq 1$ and $x = p \, x$ we get
  \[
    \sup_{x} |G_x(z)| = e^{- \frac{\mathrm{Re}\{z\}}{t}} \Phi(\sqrt{t}) \tau(p \, |h_z|^2),
  \]
  Our previous estimate then yields
  \[
    e^{- \frac{\mathrm{Re}\{z\}}{t}} \Phi(\sqrt{t}) \, \tau(p \, |h_z|^2) \, \lesssim \, t^{\frac{D_{\Phi}}{2}} \big( |z| \cos \theta \big)^{- \frac{D_{\Phi}}{2} } e^{- \frac{\beta r^2}{2}  \frac{\cos\theta}{|z|}},\\
  \]
  Choosing the parameter $t \geq 0$ to be
  $t = \mathrm{Re} \{ z \}$ gives the desired estimate.
\end{proof}

\begin{lemma} \label{decay}
If $X \in \V_+^{\wedge}$ is doubling and $\psi : G \to \R_+$ has $L_2\mathrm{GB}$, then $$\tau \Big\{ \Big( \1+ \frac{X^2}{t} \Big)^{\kappa} | h_{t ( 1 - i \xi)} |^2 \Big\}^\frac{1}{2} \, \lesssim_{(\kappa)} \, \frac{1}{\Phi(\sqrt{t})^{\frac{1}{2}}} \big( 1 + |\xi| \big)^{\kappa} \quad \mbox{for all} \quad \kappa > 0.$$
\end{lemma}

\dem
Writing $z = t ( 1 - i \xi)$ in Proposition \ref{z-Gauss} gives
  \[
    \tau \Big( \chi_{[ r, \infty)} ( X)  | h_z |^2 \chi_{[ r, \infty)} ( X) \Big)
    \lesssim
    \frac{1}{\Phi(\sqrt{t})} e^{- \frac{\beta}{2} \frac{r^2}{t }\frac{1}{( 1 + | \xi |^2)}} .
  \]
Using the spectral measure $d E_X$ of $X$ and since since $(1 + s^2)^\kappa \lesssim_{(\kappa)} 1 + s^{2 \kappa}$
  \[
    \begin{array}{>{\displaystyle}r>{\displaystyle}l>{\displaystyle}l}
      \tau \Big\{ \Big( \1+ \frac{X^2}{t} \Big)^{\kappa} | h_{t ( 1 - i \xi)} |^2 \Big\} & \lesssim_{(\kappa)} & \tau \big\{  | h_{t ( 1 - i \xi)} |^2  \big\} + \tau \Big\{ | h_{t ( 1 - i \xi)} |^2  t^{-\kappa} X^{2\kappa} \Big\} \\
                                                                                                    & \lesssim            & \frac{1}{\Phi \left( \sqrt{t} \right)} + \underbrace{\tau \Big\{ |h_{t (1 - i \xi)}  |^2 \int_{\R_+} \Big( \frac{s^2}{t} \Big)^{\kappa} dE_X (s) \Big\}}_{\mathrm{A}}. 
    \end{array}
  \]
  To estimate the term $\mathrm{A}$ we use integration by parts
  \begin{eqnarray*}
    \mathrm{A} & = & \int_{\R_+} \Big( \frac{s^2}{t} \Big)^{\kappa} \tau \big\{ |h_{t ( 1 - i \xi)}|^2 d E_X (s) \big\} \\
                & = & \int_{\R_+} \Big( \frac{s^2}{t} \Big)^{\kappa} (- \partial_s) \tau \big\{ |h_{t ( 1 - i \xi)}|^2 \chi_{[s,\infty)} (X) \big\} \\
                & = & \int_{\R_+} \frac{d}{ds} \Big( \frac{s^2}{t} \Big)^{\kappa} \tau \big\{ | h_{t ( 1 - i \xi)}  |^2 \chi_{[s,\infty)} (X) \big\} \, ds.
  \end{eqnarray*}
  In the second line, by $- \partial_s \tau \{^{}  | h_{t ( 1 - i \xi)} 
  |^2 \chi_{[ s, \infty)} ( X) \}$, we mean the Lebesgue-Stjeltjes measure
  associated with the increasing function $g ( s) = - \tau \{^{}  | h_{t ( 1 -
  i \xi)}  |^2 \chi_{[ s, \infty)} ( X) \}$ and the third line is just an
  application of the integration by parts formula for Lebesgue-Stjeltjes
  integrals. A calculation gives the desired result
  \[  
    \hskip40pt \begin{array}{>{\displaystyle}c>{\displaystyle}l>{\displaystyle}l}
      \mathrm{A} & \lesssim        & \int_{\R_+} \left( \frac{2 \kappa s^{2 \kappa - 1}}{t^{\kappa}} \right) \frac{1}{\Phi(\sqrt{t})} e^{- \frac{\beta}{2} \frac{s^2}{t} \frac{1}{( 1 + |\xi|^2)}} ds\\
                 & \sim_{(\kappa)} & \frac{(1 + | \xi |^2)^{\kappa}}{\Phi(\sqrt{t})} \int_{\R_+} s^{2 \kappa - 1} e^{- \frac{\beta}{2} s^2 } ds \ \sim_{(\kappa)} \ \frac{(1 + | \xi |)^{2 \kappa}}{\Phi(\sqrt{t})}. \hskip50pt \square
    \end{array}
  \]

\begin{proposition} \label{MainL2Est}
  Let $B_t = \lambda(m(\psi) \eta_1(t \psi))$ where 
  $\eta_1(z)=\eta(z) \, e^{-z}$ for some $\eta \in \H_0^\infty$. Assume
  also that $X$ is a doubling operator satisfying $L_2 \mathrm{GB}$, then the
  following estimate holds for every $\delta > 0$ and $\kappa > 0$
  $$\tau \Big\{ \Big( \1+ \frac{X^2}{t} \Big)^{\kappa} | B_t |^2 \Big\}^{\frac12}
    \lesssim_{(\kappa, \delta)}
    \frac{1}{\Phi(\sqrt{t})^\frac{1}{2}} \big\| m(t^{-1}\cdot ) \eta(\cdot)  \big\|_{W^{2, \kappa + \frac{1 + \delta}{2}} (\R_+)}.$$
\end{proposition}

\begin{proof}
By Fourier inversion formula $$m(s) \eta_1(t s) = \underbrace{m(s) \eta(ts)}_{m_t(ts)} e^{-t s} = \Big( \frac{1}{2 \pi} \int_{\widehat{\R}} \widehat{m}_t (\xi) e^{i \xi t s} d \xi \Big) e^{- t s}.$$ Thus, by composing with $\psi$ and applying the left regular
representation $$B_t = \frac{1}{2 \pi} \int_{\widehat{\R}} \widehat{m}_t (\xi) h_{t (1 - i \xi)} d \xi.$$ Triangular inequality for the $L_2$-norm with weight $(\1 + X^2/t)$ and Lemma \ref{decay} give
$$\begin{array}{rll} \displaystyle \tau \Big\{ \Big( \1+ \frac{X^2}{t} \Big)^{\kappa} | B_t |^2 \Big\}^\frac{1}{2} & \hskip-5pt = & \hskip-8pt \displaystyle  \tau \Big\{ \Big( \1+ \frac{X^2}{t} \Big)^{\kappa} \Big| \frac{1}{2 \pi} \int_{\widehat{\R}}  \widehat{m}_t ( \xi) h_{t ( 1 - i \xi )} d \xi \Big|^2 \Big\}^{\frac{1}{2}} \\ &\hskip-5pt \le & \hskip-8pt \displaystyle  \frac{1}{2 \pi} \int_{\widehat{\R}} | \widehat{m}_t(\xi) | \, \tau \Big\{ \Big(\1+ \frac{X^2}{t} \Big)^{\kappa} \big| h_{t (1 - i \xi )} \big|^2 \Big\}^{\frac{1}{2}} d \xi \\ [8pt] & \hskip-5pt \lesssim_{(\kappa)} & \hskip-8pt \displaystyle  \frac{1}{\Phi(\sqrt{t})^{\frac{1}{2}}} \int_{\widehat{\R}} | \widehat{m}_t (\xi)| (1 + | \xi |)^{\kappa + \frac{1 + \delta}{2}}  ( 1 + | \xi |)^{- \frac{1 + \delta}{2} } d \xi = \mathrm{A}.
\end{array}$$
H\"oder's inequality in conjunction with the definition of Sobolev space then yield
\begin{eqnarray*}
\Phi(\sqrt{t})^{\frac{1}{2}} \mathrm{A} &  \le & \Big( \int_{\widehat{\R}} \big(1 + | \xi |\big)^{-(1 + \delta) } d \xi \Big)^{\frac{1}{2}} \big\| m(t^{-1} \, \cdot ) \eta(\cdot) \|_{W^{2, \kappa + \frac{1 + \delta}{2} } (\R_+)}
\end{eqnarray*}
The the integral above is dominated by $(1 + \delta^{-1})^\frac12$ and the assertion follows.  
\end{proof}

\demB
Let $B_t = \lambda(m(\psi) \eta_1(t \psi))$ with $\eta_1(s) = e^{-s} \eta(s)$ and $B_t = \Sigma_t M_t$ be the decomposition \eqref{SMdec} with $\gamma > D_\Phi/2$ . Since we are assuming $X$ to be symmetric, we have that $\sigma|M_t|^2 = |M_t|^2$ and, by Lemma \ref{MaxLemma} and Lemma \ref{MaxCon}, $M_t$ satisfies the maximal inequality of \eqref{SM}. By Proposition \ref{MainL2Est} we have that
  \[
    \begin{array}{>{\displaystyle}c>{\displaystyle}l>{\displaystyle}l}
      \sup_{t > 0} \| \Sigma_t \|_{L_2(\V)} & \lesssim_{(\gamma)} & \sup_{t > 0} \big\| m(t^{-1} \, \cdot) \eta(\cdot) \, \big\|_{W^{2, \gamma + \frac{1 + \delta}{2} } (\R_+)}.      
    \end{array}
  \]
  Therefore $B_t = \Sigma_t M_t$ is a square-max decomposition. By similar means we obtain a max-square decomposition $B_t = M_t \Sigma_t$. Since our maximal bounds trivially extend to matrix amplifications, we may apply Theorem \ref{BoundMultiplierMaximal} in conjunction with Remark \ref{CBMultiplierMax} to deduce complete bounds of our multiplier $T_{m \circ \psi}$ in both row and column Hardy spaces. Finally, arguing as in Corollary \ref{ColandRowMultiplierMax} and noticing that $m \circ \psi \equiv m(0)$ on the subgroup $G_0 = \{g \in G : \psi(g) = 0\}$, we deduce the assertion. \fin

\begin{remark}
  \emph{
    It is interesting to observe that the proof given here can be adapted to 
    the classical case. Indeed, let $S_t = e^{-t A}$ be a Markovian semigroup acting on
    $L_\infty(X,\mu)$. Assume further that the metric measure space 
    $(X, d_\Gamma, \mu)$, where $d_\Gamma$ is the gradient metric
    \cite[Definition 3.1]{Sa2009}, is doubling, i.e.:
    \[
      \esssup_{x \in X} \sup_{r  > 0} \left\{ \frac{\mu(B_x(2 r))}{\mu(B_x(r))} \right\} < \infty
    \]
    and that its integral kernel $k_t(x,y)$ has Gaussian bounds with respect to the 
    gradient distance, i.e.:
    \[
      \left\| \chi_{[r,\infty)}(d_\Gamma(x,\cdot)) \, k_t(x,\cdot) \right\|_2^2 \, \lesssim \, \frac{e^{- \beta \, \frac{r^2}{t}}}{\mu(B_x(\sqrt{t}))}.
    \]
    In that case we can apply the well known covering arguments for doubling
    spaces to prove that the Hardy-Littlewood maximal operator is of weak type $(1,1)$
    and by interpolation the $\mathrm{HL}$ inequalities hold. Since
    $(X,d_\Gamma, \mu)$ is a doubling metric measure space 
    with bounded Hardy-Littlewood maximal inequalities and Gaussian Bounds 
    we can apply the results above to reprove the classical spectral
    H\"ormander-Mikhlin theorem as stated in \cite{DuOuSi2002}. We shall
    consider this a new proof of the classical spectral H\"ormander-Mikhlin.
    Interestingly, some of the steps of the proof are parallel to that of
    \cite{DuOuSi2002} even when the main idea of our approach is to use maximal inequalities
    instead of Calder\'on-Zygmund estimates for the kernels. 
  }
\end{remark}

\subsection{ \hskip3pt The $q$-Plancherel condition}

In this subsection we shall refine our results by proving Theorem B ii). Our first task is to introduce the noncommutative form of the Plancherel condition assumed in the statement. 

\begin{definition}
  Let $(\M,\tau)$ be a noncommutative measure space and let $\S$ be a submarkovian  semigroup generated by $A$.
  We say that $\S$ satisfies the completely bounded $q$-Plancherel condition, denoted by $\mathrm{CBPlan}^{\Phi}_q$, where $\Phi$ is some
  increasing function and $q \in  (2, \infty]$, whenever $$\| {F (A)}  \|_{\CB( L_2(\V), \V)} \lesssim \frac{1}{\Phi(\sqrt{t})^{\frac{1}{2}}} \| F(t^{-1} \, \cdot ) \|_{L_q(\mathbbm{R}_+)},$$
  for every $t > 0$ and for every function $F: \R_+ \to \R_+$ with $\mathrm{supp}(F) \subset \left[0, t^{- 1}\right]$.
\end{definition}

\begin{remark}
  \emph{
    In the context of this paper $\M = \V$ for some LCH unimodular group
    $G$ endowed with its canonical trace and $\S = (T_{e^{-t \psi}})_{t \geq 0}$ is a semigroup of convolution type. In that case $F(A) = T_{F(\psi)}$ and
    by Theorem \ref{CBL2Multiplier} and Remark \ref{CBL2Norm} we have that
    \begin{eqnarray*}
    \| T_{F(\psi)} \|_{{\CB( L_2(\V), \V)}} & = & \| T_{F(\psi)} \|_{{\CB( L^r_2(\V), \V)}} \\ & = & \| T_{F(\psi)} \|_{{\CB( L^c_2(\V), \V)}} \ = \ \| F(\psi) \|_{L_2(G)}.
    \end{eqnarray*}
    Thus, the $\mathrm{CBPlan}^{\Phi}_q$ condition can be restated as a
    bound on the ${\CB( L^\dag_2(\V), \V)}$ norm, where $\dag$ is either the column or
    the row o.s.s. of $L_2(\V)$, or as a bound in the $L_2(G)$-norm of the symbol
    $F(\psi)$. Furthermore, since $\psi$ determines $\S$ we will sometimes say that 
    $\psi$ has the $\mathrm{CBPlan}_q^\Phi$.
  }
\end{remark}

For every $F$ with $\mathrm{supp}(F)\subset \left[ 0, t^{- 1} \right]$ we have
that $F(t^{-1} \, \cdot)$ is supported in $[0, 1]$. Using that 
$L_q([0,1]) \subset L_p([0,1])$, with contractive inclusion, we see that 
$\mathrm{CBPlan}_{p}^{\Phi} \Rightarrow \mathrm{CBPlan}_{q}^{\Phi}$ for 
$p \leq q$.

\begin{proposition}
  \label{SobImpCBPlan}
  Let $(G,\psi)$ be a pair formed by a LCH unimodular
  group and a c.n. length. Let $\Phi$ be a doubling
  function. If $\psi$ satisfies the utracontractivity estimates
  $\mathrm{CBR}_{\Phi}^{2, \infty}$ then it satisfies 
  $\mathrm{CBPlan}_{\infty}^{\Phi}$.
\end{proposition}

\begin{proof}
  We pick $s > 0$, to be chosen later, and notice that $$F(\psi(g)) = F(\psi(g)) e^{s \psi (g)} e^{- s \psi(g)} = G_s(\psi(g)) e^{- s \psi}$$
  where $G_s$ is a bounded function with 
  $\| G_s \|_{\infty} \leq \| F \|_{\infty} e^{s / t}$. Therefore
  \begin{eqnarray*}
  \| T_{F ( \psi)} \|_{\CB(L_2(\V), \V)} & = & \| T_{G_s(\psi)} S_s \, \|_{\CB( L_2(\V), \V)} \\ [2pt] & \le & \| T_{G_s(\psi)} \, \|_{\CB(L_2(\V))}  \| S_s \|_{\CB(L_2(\V), \V)} \\ & \lesssim & \| F \|_{\infty} e^{s / t}\Phi(\sqrt{s})^{-\frac{1}{2}}.
  \end{eqnarray*}
  Making $s = t$ and noticing that
  $\| F \|_{\infty} = \| F(t^{-1} \, \cdot) \|_{\infty}$ gives
  the desired result.
\end{proof}

The terminology of the $q$-Plancherel condition comes from the so-called spectral Plancherel measures which arise in the study of spectral
properties of infinitesimal generators of Markovian semigroups 
over some measure spaces \cite{Si1996,DuOuSi2002}. In the 
case of a semigroup of Fourier multipliers generated by a c.n. length we can
define the Plancherel measure $\mu_\psi$, as the only $\sigma$-finite measure
over $\R_+$ satisfying that for every $F \in C_c(\R_+)$
\begin{equation}
  \| T_{F(\psi)} \|_{\CB(L_2(\V),\V)} = \Big( \int_{\R_+} | F(s)|^2 \, d \mu_{\psi}(s) \Big)^{\frac{1}{2}}.
  \label{PlanMeasure}
\end{equation}
It is trivial to see that 
$d \mu_\psi(r) = \partial_r \mu(\{g \in G : \psi(g) \leq r \})$, where
$\partial_r$ represents the Lebesgue-Stjeltjes derivative of the increasing
function $g(r) \! = \! \mu(\{g \in G \! : \psi(g) \leq r \})$.

\subsubsection{{\bf Characterization of the $q$-Plancherel condition}}
By formula (\ref{PlanMeasure}) the $\CB(L_2(\V), \V)$ norm of $T_{F(\psi)}$ can be
expressed as an integral of $F$. The following lemma (whose proof is straightforward and we shall omit) allows to express the 
$\mathrm{CB Plan}_q^{\Phi}$ property as a $L_{(q/2)'}(\R_+)$ bound on
$\mu_\psi$. 

\begin{lemma}
  \label{LpMeasures}
  Let $(\Omega, \Sigma)$ be a measurable space and consider two measures $\mu$, $\nu$ on it. Assume in addition that $\mu$ is a positive measure. Then, we have the inequality 
  \begin{equation}
    \left| \int_{\Omega} f(\omega) d \nu ( \omega) \right| \leq K \| f \|_{L_p ( d \mu)}
    \label{LpMeasure}
  \end{equation}
  if and only if $\nu \ll \mu$ and $\phi = {d \nu}/{d \mu}$ satisfies $\| \phi \|_{L_{p'} (d \mu)} \leq K$. Furthermore, the optimal $K$ in \eqref{LpMeasure} is precisely $\| \phi \|_{L_{p'} (d \mu)}$. If $\nu$
  is also positive, it is enough for \eqref{LpMeasure} to hold only for positive
  functions.
\end{lemma}


\begin{proposition}
  \label{CharacterizationPlan}
  Let $G$ be a LCH unimodular group equipped with a c.n. length $\psi: G \to \R_+$. Then, this pair satisfies the $\mathrm{CBPlan}_q^{\Phi}$ property with respect to some increasing function $\Phi: \R_+ \to \R_+$ if and only if $d \mu_{\psi}(r) =
  \partial_r \mu \{ g \in G : \psi(g) \leq r \}$ fulfills the following conditions$\hskip1pt :$
  \begin{itemize}
    \item[i)] $d \mu_{\psi} \ll dm$. 
    
    \vskip5pt
    
    \item[ii)] $\displaystyle \Big\| \frac{d \mu_{\psi}}{dm} \chi_{[0, R]} \Big\|_{L_{( q/2 )'}(\R_+)} \lesssim \Phi( R^{- \frac{1}{2}})^{- 1} R^{- \frac{2}{q}}$ for every $R > 0$.
      \end{itemize}
\end{proposition}

\begin{proof}
  Let $t = 1/R$ and $G(s) = |F(s)|^2$. By (\ref{PlanMeasure}), $\mathrm{CBPlan}_q^{\Phi}$ is equivalent to
  \begin{eqnarray*}
  \int_0^R G(s) \, d \mu_\psi(s) & \lesssim & \Phi( R^{- \frac{1}{2}})^{- 1} \Big( \int_0^1 | F(t^{-1} s) |^q \, ds \Big)^\frac{2}{q} \\
  & = & \Phi( R^{- \frac{1}{2}})^{- 1} R^{- \frac{2}{q}} \Big( \int_0^R | G(s) |^\frac{q}{2} \, ds \Big)^\frac{2}{q}.
  \end{eqnarray*}
  Then, the result follows applying Lemma \ref{LpMeasures} to $(\Omega, d\nu, d\mu) = (\R_+, d\mu_\psi, dm)$.
\end{proof}

The result above uses the crucial fact that  the spectrum of the semigroup
$\S$ generated by $\psi$ can be identified with $G$. Therefore, spectral 
properties of the semigroup can be translated into geometrical properties of
$G$. It is also interesting to note that the characterization in Proposition 
\ref{CharacterizationPlan} can be expressed as a bound for the size of 
the spheres associated to the pseudo-metric
$d_\psi(g,h) = \psi(g^{-1} h)^{1/2}$.

\subsubsection{{\bf Stability under direct products}}

Consider two pairs $(G_j,\psi_j)$ of LCH unimodular groups equipped with c.n. lengths for $j=1,2$. Then it is clear that $\psi : G_1 \times G_2 \to \R_+$ given by $\psi(g,h) = \psi_1(g) + \psi_2(h)$ is also a c.n. length. Notice that
\begin{eqnarray*}
\big\| T_{F(\psi)} \big\|_{\CB(L_2(\V), \V)}^2 & = & \int_{G_1 \times G_2} \big| F(\psi_1(g) + \psi_2(h)) \big|^2 \, d\mu_{G_1}(g) \, d\mu_{G_2}(h) \\
                                         & = & \int_{\R_+} \int_{\R_+} \big| F(\xi + \zeta) \big|^2 \, d\mu_{\psi_1}(\xi) \, d\mu_{\psi_2}(\zeta) \\
                                         & = & \int_{\R_+} |F(\xi)|^2 d\, (\mu_{\psi_1} \ast \mu_{\psi_2})(\xi).
\end{eqnarray*}
Thus, the Plancherel measure is 
$\mu_{\psi} = \mu_{\psi_1} \ast \mu_{\psi_2}$ and we obtain the following result. 

\begin{theorem}
  Assume $(G_j,\psi_j)$ satisfy $\mathrm{CBPlan}_{q_j}^{\Phi_j}$ for $j=1,2$.
  Then the pair $(G_1 \times G_2, \psi)$ defined above satisfies the $\mathrm{CBPlan}_{q}^{\Phi}$ property with $\Phi = \Phi_1 \Phi_2$ and with $$q = \max \Big\{ 2, \Big( \frac{1}{q_1} + \frac{1}{q_2} \Big)^{-1} \Big\}.$$
\end{theorem}

\begin{proof}
The result is a simple consequence of Young's inequality for convolutions and we shall just sketch the argument for the (slightly more involved) case where $1/q_1 + 1/q_2 > 1/2$, so that $q=2$. According to Proposition \ref{CharacterizationPlan}, it suffices to see that $$\Big\| \frac{d\psi_1}{dm} * \frac{d\psi_2}{dm} \Big\|_{L_\infty(0,R)} \le \frac{1}{R \, \Phi_1(R^{-1/2}) \Phi_2(R^{-1/2})}.$$ The $\mathrm{CBPlan}_{q_1}^{\Phi_1}$ property of $(G_1, \phi_1)$ implies $$\Big\| \frac{d\psi_1}{dm} * \frac{d\psi_2}{dm} \Big\|_\infty \le \Big\| \frac{d\psi_1}{dm} \Big\|_{(\frac{q_1}{2})'} \Big\| \frac{d\psi_2}{dm} \Big\|_{\frac{q_1}{2}} \le \frac{1}{R^{\frac{2}{q_1}} \, \Phi_1(R^{-1/2})} \Big\| \frac{d\psi_2}{dm} \Big\|_{\frac{q_1}{2}}.$$ Now, since $1/q_1 + 1/q_2 > 1/2$ it turns out that $$\frac{1}{q_1/2} = \frac{1}{(q_2/2)'} + \frac{1}{r} \Rightarrow \Big\| \frac{d\psi_2}{dm} \Big\|_{\frac{q_1}{2}} \le R^{\frac{1}{r}} \Big\| \frac{d\psi_2}{dm} \Big\|_{(\frac{q_2}{2})'}.$$ The result follows from the characterization of $\mathrm{CBPlan}_{q_2}^{\Phi_2}$ in Proposition \ref{CharacterizationPlan}. 
\end{proof}

\begin{remark}
  \emph{
    A result along the same lines can be obtained for crossed products
    under invariance assumptions on $\psi_1$. This goes in the same spirit as Remark \ref{RemCrossed}.
  }
\end{remark}

\subsubsection{{\bf Refinement of the smoothness condition}}

Here we are going to see how we can prove the optimal smoothness order in the H\"ormander-Mikhlin condition of Theorem B ii) when $\psi$ satisfies the $\mathrm{CBPlan}_q^{\Phi}$ property. We need several preparatory lemmas.
In the next one we denote by $W_{\eta}^{p,s}(\R_+)$, where 
$\eta \in \H_0^{\infty}$, the Sobolev space given by completion 
with respect to the norm
\[
  \| f \|_{W_{\eta}^{p,s}(\R_+)} = \big\| (1-\partial_x^2)^{s/2} (\eta f) \big\|_p.
\]

\begin{lemma}
  \label{inclusionSob}
  Given $f,g : \R_+ \to \C$, the following holds:
  \begin{itemize}
    \item[i)] For every $\varepsilon > 0$
    \[
      \big\| (1 - \partial_x^2)^{s / 2} (f g) \big\|_2
      \lesssim_{(s, \varepsilon)}
      \big\| (1 - \partial_x^2)^{( s + 1 + \varepsilon) / 2} f \big\|_{\infty} \big\| (1 - \partial_x^2)^{s / 2} g \big\|_2.
    \]
    \item[ii)] If $\rho(z) = z^s e^{-z}$ and $\eta \in \H_0^{\infty}$ $$\big\| (1 - \partial^2_x)^{s/2}(\eta \rho f) \big\|_2 \lesssim_{(s,\varepsilon)} \big\| (1-\partial^2_x)^{(s + 1 + \varepsilon)/2} (\eta f) \big\|_{\infty}. $$ Equivalently, we find the embedding $W_{\eta}^{\infty, s + 1 +\varepsilon}(\R_+) \subset_{(s,\varepsilon)} W_{\eta \rho}^{2,s}(\R_+)$.
  \end{itemize}
\end{lemma}

\begin{proof}
  The second point follows immediately from the first one by noticing that 
  $\rho(z) = z^s e^{-z}$ has finite $W^{2,s}(\R_+)$ norm. We are going to 
  prove the first point for $s \in \N$ and use interpolation. Given $s \in \N$, we have
  $$\begin{array}{rll} \big\| ( 1 - \partial_x^2)^{s / 2} (fg)  \big\|_2 & \sim & \displaystyle \sum_{k = 0}^s \| \partial_x^k  (fg) \|_2 \\ & = & \displaystyle  \sum_{k = 0}^s \Big\| \sum_{j = 0}^k \binom{k}{j} (\partial_x^j f) (\partial_x^{k - j} g) \Big\|_2 \\
  & \lesssim_{(s)} & \displaystyle  \Big( \max_{0 \le j \le s} \| \partial_x^j f \|_{\infty} \Big)  \Big( \sum_{k = 0}^s  \| \partial_x^k g \|_2  \Big) \\ [12pt] & \sim & \displaystyle \Big( \max_{0 \le j \le s} \| \partial_x^j f \|_{\infty} \Big) \big\| (1 - \partial_x^2)^{s / 2} g \big\|_2. \end{array}$$
  Thus, all we have to see is that for every $j \in \{0,1,2, ...,s\}$
  \[
    \big\| \partial_x^j (1 - \partial_x^2)^{-(s + \varepsilon + 1)/2} f \big\|_\infty \lesssim_{(s,\varepsilon)} \| f \|_\infty.
  \]
  Recall that if the symbol of a Fourier multiplier is given by the Fourier 
  transform of finite measure, then it is bounded in $L_\infty(\R)$. Thus, we 
  just need to see that there is a finite measure $\mu_{j,s}$ such that
  \begin{eqnarray*}
    \widehat{\mu}_{j,s}(\xi) & = & \frac{\xi^j}{(1+|\xi|^2)^{\frac{s + \varepsilon + 1}{2}}}\\
                         & = & \mathrm{sgn}(\xi)^j \, \frac{1}{(1 + |\xi|^2)^{\frac{ s + \varepsilon - j + 1}{2}}} \, \frac{|\xi|^j}{(1+|\xi|^2)^{\frac{j}{2}}} \ = \ ( H_{[j]}(\nu_{s,j}) \ast m_{j} )^{\wedge}(\xi),
  \end{eqnarray*}
  where $H_{[j]}$ is the Hilbert transform for $j$ odd and the identity map for $j$ even. 
  By \cite[V.3, Lemma 2]{Ste1970Singular}
  $m_j$ is a finite measure.  Therefore, it is enough to see that if 
  $\widehat{\nu}_{s,j}(\xi) = 1/(1 + |\xi|^2)^{(s + \varepsilon - j + 1) /2}$, then
  $H_{[j]}(\nu_{s,j})$ is a finite measure. Applying the Hilbert transform or identity map to
  \cite[V.(26)]{Ste1970Singular} gives the desired result.
\end{proof}

\begin{lemma}
  \label{decaySupp}
  Assume $G$ is a LCH unimodular group, $\psi : G \to \R_+$ is a c.n. length
  and that they satisfy the
  $\mathrm{CBPlan}_q^{\Phi}$ property. If 
  $\eta_1$, $\eta_2 \in \H^\infty_0(\Sigma_\theta)$,
  with $\eta_1$ satisfying that there is $\gamma > 0$ such that
  $|\eta_1(z)| \lesssim e^{-\gamma \mathrm{Re}(z)}$ for all 
  $z \in \Sigma_\theta$, then the following estimate holds for all $m \in L_\infty(\R_+)$ $$\big\| \lambda \big( m(\psi) \eta_1(t \psi) \eta_2(t \psi) \big) \big\|_{L_2(\V)}
    \lesssim_{( D_\Phi, q, \gamma)} \frac{1}{\Phi(\sqrt{t})^\frac{1}{2}} \big\| m(t^{-1} \, \cdot) \eta_2(\cdot) \big\|_{L_q (\R_+)}.$$
\end{lemma}

\dem Using integration by parts we obtain 
\begin{eqnarray*}
\big\| \lambda \big( m(\psi) \eta_1(t \psi) \eta_2(t \psi) \big) \big\|_{L_2 (\V)} & = & \Big\|  \int_{\R_+} \lambda \big( m(\psi) \eta_1'(r) \eta_2(t \psi) \chi_{[0, r)} (t \psi) \big) dr \Big\|_{L_2 (\V)} \\ & \leq             & \int_{\R_+} \eta_1'(r) \big\| \lambda \big( m(\psi) \eta_2(t \psi) \chi_{[0, r)} (t \psi) \big) \big\|_{L_2 (\V)} \, dr.
\end{eqnarray*}
Nos, applying the $\mathrm{CBPlan}_q^\Phi$ property, we obtain 
  \[
    \begin{array}{>{\displaystyle}c>{\displaystyle}l>{\displaystyle}l}
      & & \hskip-80pt \big\| \lambda \big( m(\psi) \eta_1(t \psi) \eta_2(t \psi) \big) \big\|_{L_2 (\V)} \\ & \lesssim_{(q)}   & \int_{\R_+} \eta_1'(r) \frac{1}{\Phi(\sqrt{t / r})^{\frac{1}{2}}} \big\| m((r / t) \cdot) \eta_2(r \cdot) \big\|_{L_q ( [ 0, 1])} dr \\
                                                                      & =                & \Big( \int_{\R_+} \eta_1'(r) \frac{r^{- 1 / q}}{\Phi(\sqrt{t / r})^{\frac{1}{2}}} dr \Big) \big\| m(t^{-1} \, \cdot) \eta_2(\cdot) \big\|_{L_q (\R_+)}.
    \end{array}
  \]
So, we just need to estimate the integral in the right hand side term $$\int_{\R_+} \eta_1'(r) \frac{r^{- 1 / q}}{\Phi(\sqrt{t / r})^{\frac{1}{2}}} dr = \Big\{ \int_0^1 + \sum_{j = 0}^{\infty} \int_{4^j}^{4^{j + 1}} \Big\} \eta_1'(r) \frac{r^{- 1 / q}}{\Phi(\sqrt{t / r})^{\frac{1}{2}}} dr = A + \sum_{j = 0}^{\infty} B_j.$$ The first term is bounded as follows $$A \le \frac{1}{\Phi(\sqrt{t})^{\frac{1}{2}}} \int_0^1 \eta_1'(r) r^{- 1 / q} dr \lesssim_{(q)} \frac{1}{\Phi(\sqrt{t})^{\frac{1}{2}}}.$$
  For the rest of the terms, we apply the doubling condition to obtain $$B_j \le 3 \cdot 4^j \| \eta_1'\|_{L_\infty([4^j,4^{j+1}))} \frac{2^{\frac{D_{_{\Phi}}}{2} ( j + 1)}}{\Phi(\sqrt{t})^{\frac12}} =  \frac{3 \cdot 2^{\frac{D_{\Phi}}{2}}}{\Phi(\sqrt{t})^{\frac{1}{2}}} \, \| \eta_1'\|_{L_\infty([4^j,4^{j+1}))} \, 2^{\left( \frac{D_{\Phi}}{2} + 2 \right) j}.$$
  The function $\eta_1$ decreases exponentially and so does $\eta_1'$. Therefore
  $\eta_1'(z) \lesssim e^{- \gamma z}$ for $\mathrm{Re}\{z\}$ large enough.
  That allows us to sum up all the terms in the series obtaining $\sum_j B_j \lesssim \Phi(\sqrt{t})^{-\frac{1}{2}}$ up to a constant depending on $( D_\Phi,\gamma)$, as desired.
\fin

\begin{proposition}
  \label{weightedToSmoothOpt}
  Assume $G$ is a LCH unimodular group, $\psi : G \to \R_+$ is a c.n. length
  and that they satisfy the
  $\mathrm{CBPlan}_q^{\Phi}$ property. Assume in addition that $X \in \V_+^{\wedge}$ is
  doubling and admits $L_2\mathrm{GB}$. Then, we find for $\kappa, \delta, \varepsilon > 0$
  \[
    \tau \Big\{ \Big( \1 + \frac{X^2}{t} \Big)^{\kappa} | B_t |^2  \Big\}^{\frac{1}{2}}
    \lesssim_{(D_\Phi, q, \kappa, \delta, \varepsilon)}
    \frac{1}{\Phi(\sqrt{t})^{\frac{1}{2}}} \big\| m(t^{-1} \, \cdot) \eta(\cdot) \big\|_{W^{p, \kappa + \delta}(\R_+)},
  \]
  where $B_t = \lambda \left( m(\psi) \eta(t \psi) e^{-2 t \psi} (t \psi)^a \right)$, 
  $\eta$ is a $\H_0^\infty$-cut-off and  $a  = 2 \kappa / \delta + (1 + \varepsilon)/2$.
\end{proposition}

\begin{proof}
  Fix $\kappa, \delta, \varepsilon > 0$ and $a = 2 \kappa / \delta + (1+\varepsilon)/2$. We define 
  the linear, unbounded map $K_t : D \subset L_\infty(\R_+) \to L_2(\V)$ 
  by $K_t(m) = \lambda(m(t \psi) \eta(t \psi) e^{-2 t \psi} (t \psi)^a)$. 
  Using Lemma \ref{decaySupp} with $\eta_1(z) = z^a e^{-2 z}$ and
  $\eta_2(z) = \eta(z)$ gives that
  \begin{equation}
    \Big\| K_t : W_\eta^{q, 0}(\R_+) \to L_2(\V) \Big\| \lesssim_{( D_\Phi, q)} \frac{1}{\Phi(\sqrt{t})^\frac{1}{2}}.
    \label{interp1}
  \end{equation}
  Let us denote by $\phi_{t,\kappa}$ the family of weights given by $\phi_{t,\kappa}(x) = \tau \{ ( \1 + t^{-1} X^2)^{\kappa} x\}$
  and let $L_2(\V,\phi_{t,\kappa})$ be the Hilbert spaces associated to the GNS
  construction of $\phi_{t,\kappa}$. We know from Proposition \ref{MainL2Est} that 
  \[
    \left\| K_t : W_{\eta \rho}^{2,s+\frac{1+\varepsilon}{2}}(\R_+) \to L_2(\V,\phi_{t,s}) \right\|
    \lesssim_{(\kappa, \delta, \varepsilon)}
    \frac{1}{\Phi(\sqrt{t})^\frac{1}{2}},
  \]
  where $s = 2 \kappa / \delta$ and $\rho(z) = z^a e^{-z}$. Composing with the inclusion 
  $$W_{\eta}^{q,s + \frac{1+\varepsilon}{2} + 1 + \varepsilon'}(\R_+) \subset_{(s,\varepsilon')} W_{\eta \, \rho}^{2,s + \frac{1+\varepsilon}{2}}(\R_+),$$ which follows by interpolation from Lemma \ref{inclusionSob} for $q=\infty$ and the trivial inclusion for $q=2$, gives 
  \begin{equation}
    \left\| K_t : W_{\eta}^{q,s + \frac{1+\varepsilon}{2} + 1 + \varepsilon'}(\R_+) \to L_2(\V,\phi_{t,s}) \right\|
    \lesssim_{(\kappa, \delta, \varepsilon, \varepsilon')}
    \frac{1}{\Phi(\sqrt{t})^\frac{1}{2}}.
    \label{interp2}
  \end{equation}
  Notice that the spaces obtained through GNS construction 
  $L_2(\V, \phi_{t,\kappa})$ are well behaved with respect to the complex 
  interpolation method. In particular, the expected identity below holds 
  \[
    \big[ L_2\left(\V,\phi_{t,\kappa_1}\right), L_2\left(\V,\phi_{t,\kappa_2}\right) \big]_\theta = L_2 \big( \V,\phi_{t,(1-\theta)\kappa_1 + \theta  \kappa_2} \big).
  \]
  Therefore, interpolating (\ref{interp1}) and (\ref{interp2}) with 
  $\theta = \delta/2$ yields
  \[
    \Big\| K_t : W_{\eta}^{q,\kappa + \frac{\delta}{2} (\frac{1+\varepsilon}{2} + 1 + \varepsilon')}(\R_+) \to L_2(\V,\phi_{t,\theta s}) \Big\|
    \lesssim_{(D_\Phi, q, \kappa, \delta, \varepsilon, \varepsilon')}
    \frac{1}{\Phi(\sqrt{t})^\frac{1}{2}}.
  \]
  Finally, choosing $\varepsilon$ and $\varepsilon'$ such that 
  $((1 + \varepsilon)/2 + 1 + \varepsilon') \leq 2$ gives
  \[
    \Big\| K_t : W_{\eta}^{q,\kappa + \delta}(\R_+) \to L_2(\V,\phi_{t,\kappa}) \Big\|
    \lesssim_{(D_\Phi, q, \kappa, \delta)}
    \frac{1}{\Phi(\sqrt{t})^\frac{1}{2}}.
  \]
  Therefore, applying this bound to the function $m(t^{-1} \cdot)$ proves the assertion. 
\end{proof}

\demBB 
  Let $s > D_\Phi/2$. For any $\eta \in \H_0^{\infty}$ and 
  $\delta, \varepsilon > 0$ we can define $\eta_1(z) = \eta(z) e^{-2z} z^a$,
  where $a = 2 s / \delta + (1 + \varepsilon)/2$. Set 
  $B_t = \lambda(m(\psi) \eta_1(t \psi))$ and apply \eqref{SMdec}. By Proposition 
  \ref{weightedToSmoothOpt}
  \[
    \begin{array}{>{\displaystyle}r>{\displaystyle}l>{\displaystyle}l}
      \sup_{t > 0} \| \Sigma_t \|_{L_2(\V)} & \lesssim_{(D_\Phi, q, s, \delta, \varepsilon)} & \sup_{t > 0} \big\| m(t^{-1} \, \cdot) \eta(\cdot) \big\|_{W^{p, s + \delta}(\R_+)}.
    \end{array}
  \]
  Once this is settled, the argument continues as in the proof of Theorem B i). 
\fin

\subsection{ \hskip3pt An application for finite-dimensional cocycles}

Our aim is to recover the main result in \cite{JunMeiPar2014}
for the case of radial multipliers to illustrate how the Sobolev dimension
approach is, a priori, more flexible than the one used in \cite{JunMeiPar2014}.
We will start proving that c.n. lengths coming from surjective and proper
finite-dimensional cocycles satisfy the standard assumptions. Then we will 
reduce the case of general finite-dimensional cocycles to surjective
and proper ones. 

Let $b : G \to \R^n$ be a finite-dimensional cocycle. Assume that $b$ is 
surjective and proper (i.e. $b^{-1}[K]$ is a compact set for every compact $K$). 
Then the pullback of the Haar measure $b^{\ast} \mu(E) = \mu(b^{-1}[E])$ in $\R^n$ 
is translation invariant and therefore satisfies satisfies that $d \, b^{\ast} \mu(\xi) = c d \, \xi$. 
Indeed, let $\alpha: G \to \mathrm{Aut}(\R^n)$ be the action naturally associated to $b$. Given 
a Borel compact set $E \subset \R^n$ with $b^{-1}(E) = A \subset G$ and since $b(gA) = \alpha_g(b(A)) 
+ b(g)$, we conclude that $$b^*\mu(E) = \mu(A) = \mu(gA) = b^*\mu(\alpha_g(E) + b(g)).$$ Note that $\mu(A)$ 
is well-defined and finite since $b$ is continuos and proper.  Applying this identity to the $\alpha$-invariant sets 
$E = B_r(0)$ and using the subjectivity of $b$, we conclude the assertion. 
An important consequence of this fact is that $$ \| S_t \|^2_{\CB(L_2(\V), \V)} = \int_G |e^{-t \psi(g)}|^2 \, d \mu(g) = \int_{\R^n} e^{-2 t | \xi |^2} d \, (b^{\ast} \mu)(\xi) = \frac{1}{\Phi(\sqrt{t})},$$
where $\S = (S_t)_{t \geq 0}$ is the semigroup associated with $\psi(g) = \| b(g) \|^2$ and $\Phi(t) \sim t^n$. Therefore, the semigroup 
associated to any proper and surjective finite-dimensional 
cocycle satisfies the  $\mathrm{CBR}^{\Phi}$ property.
In the same way, the measure $\mu_\psi$ defined in (\ref{PlanMeasure})
can be expressed (using polar coordinates) as in terms of $b^{\ast} \mu$ and a trivial calculation
gives that $\psi$ has the $\mathrm{CB Plan_2^{\Phi}}$ property. 
We need to find a suitable $X_b \in \V^{\wedge}_+$. We shall prove that $b$
induces a natural transference map from functions $f : \R^n \to \C$ into operators
$x \in \V$ given by $$\mathcal{J}(f) = \lambda(\widehat{f} \circ b).$$ Therefore, if $\mathcal{R}$ is
a distribution in $\R^n$ such that $\widehat{\mathcal{R}}(x) = |x|$, our choice will be
$X_b \hskip-1pt = \hskip-1pt \lambda(\mathcal{R}(b))$. \hskip-2pt Before proving $X_b \in \V_+^\wedge$ we will need the following auxiliary result.

\begin{lemma}
  \label{transMult0}
  If $\varphi_j: \R^n \to \C$ are radial $L_1$-functions $$\lambda ( \varphi_1 \circ b ) \, \lambda ( \varphi_2 \circ b ) = \lambda \big( ( \varphi_1 \ast \varphi_2 ) \circ b \big)$$ for any group $G$ equipped with a proper and surjective cocycle $b: G \to \R^n$.  
  \end{lemma}

\begin{proof}
  Up to constants, we know that $d (b^{\ast} \mu) = d m$, so that
  \begin{eqnarray*}
    (\varphi_1 \circ b) * (\varphi_2 \circ b) (g) & = &  \int_G \varphi_1(b(h)) \varphi_2(b(h^{-1}g)) \, d \mu(h) \\ & = & \int_G \varphi_1(b(h)) \varphi_2(\alpha_{h^{-1}}(b(g) - b(h))) \, d \mu(h) \\
                                                             & = & \int_G \varphi_1(b(h)) \varphi_2(b(g) - b(h)) \, d \mu(h)\\
                                                             & = & \int_{\R^n} \varphi_1(\zeta) \varphi_2(b(g) - \zeta) d \, (b^{\ast} \mu) (\zeta)\\
                                                             & = & \int_{\R^n} \varphi_1(\zeta) \varphi_2(b(g) - \zeta) d \, \zeta \ = \ (\varphi_1 \ast \varphi_2)(b(g)).
  \end{eqnarray*}
Taking the left regular representation at both sides yields the assertion. \end{proof}

It is straightforward to restate Lemma \ref{transMult0} in terms of the transference operator $\mathcal{J}$. Namely, we shall be working with the following subclasses of radial functions in the Euclidean space $\R^n$
\begin{eqnarray*}
\mathcal{A} & = & \Big\{ \phi: \R^n \to \C \, \big| \ \phi \mbox{ radial}, \ \widehat{\phi} \mbox{ is a finite measure in }\R^n \Big\}, \\
\mathcal{A}_+ & = & \Big\{ \phi: \R^n \to \C \, \big| \ \phi \mbox{ radial and positive}, \ \widehat{\phi} \mbox{ is a finite measure in }\R^n \Big\}.
\end{eqnarray*}
Observe that $\phi_j  \in \mathcal{A}$ implies by Lemma \ref{transMult0} that
\begin{equation}
\mathcal{J}(\phi_1 \, \phi_2) = \mathcal{J}(\phi_1) \, \mathcal{J}(\phi_2).
\label{transMult}
\end{equation}
In fact, we will make use of the following consequences:
\begin{itemize}
\item[i)] $\mathcal{J}: \mathcal{A} \to \V$ is completely bounded.

\item[ii)] $\mathcal{J}(\mathcal{A})$ is an abelian subalgebra of $\V$.
\end{itemize}
Indeed, it follows from \eqref{transMult} that $\mathcal{J}$ is an $*$-homomorphism on $\mathcal{A}$. In particular, it is completely positive and its c.b. norm it determined by $\mathcal{J}(\1)$. The Fourier transform of $\1$ is the Dirac delta $\delta_0$ at $0$. Let us approximate 1 in the weak-$*$ topology by $h_\delta(\xi) = \exp(-\delta |\xi|^2)$ as $\delta \to 0^+$. By the weak-$*$ continuity of $\mathcal{J}$, it turns out that $$\big\| \mathcal{J}(\1) \big\|_{\V} = \lim_{\delta \to 0^+} \big\| \lambda (\widehat{h}_\delta \circ b ) \big\|_{\V}  \le \lim_{\delta \to 0^+} \big\| \widehat{h}_\delta \circ b \big\|_{L_1(\V)} = \lim_{\delta \to 0^+} \int_{\R^n} \widehat{h}_\delta(\xi) \, d\xi = 1.$$ Thus $\mathcal{J}$ is a completely positive contraction. Once this is settled, ii) follows from \eqref{transMult}. In order to define $X_b$ as an element of $\V_+^\wedge$, we need to express it as the supremum of positive operators in $\V$. We use
\[
  1 = \int_{\R_+} s \, |\xi|^2 \, e^{- s |\xi|^2} \frac{d \, s}{s},
\]
and think of $\eta_s(\xi) = |\xi|^2 \, s \, e^{- s |\xi|^2}$ as a
continuous partition of the unit. Hence $$|\xi| = \int_{\R_+} |\xi| \, \eta_s(\xi) \frac{d \, s}{s} \leadsto \phi_{\varepsilon, R}(\xi) := \int_{\varepsilon}^R |\xi| \, \eta_s(\xi) \frac{d \, s}{s} \leadsto X_b := \sup_{0 < \varepsilon \leq R < \infty} \mathcal{J}(\phi_{\varepsilon, R}).$$ This presents $X_b$ as a well-defined element of the extended positive cone $\V_+^{\wedge}$.

\begin{theorem}
\label{fdStandard}
Let $G$ be a LCH unimodular group and consider an $n$-dimensional proper and surjective cocycle $b : G \to \R^n$ equipped with the conditionally negative length $\psi(g) = \|b(g)\|^2$. Then $(G,\psi,X_b)$ satisfies the standard assumptions.
\end{theorem}

\dem
  We will start by proving the $L_2 \mathrm{GB}$. By noticing that 
  $\zeta \mapsto \chi_{[r,\infty)}(\zeta)$ is an increasing function 
  and the normality of the weight
  $x \mapsto \tau \left\{ x \, |\lambda(e^{- t \psi})|^2\right\}$ we obtain
  that
  \[
    \tau \Big\{ \chi_{[r,\infty)}(X_b) \, |\lambda(e^{- t \psi})|^2\Big\}
    =
    \sup_{0 < \varepsilon \leq R < \infty} \tau \Big\{ \chi_{[r,\infty)}(\mathcal{J}(\phi_{\varepsilon, R})) \, |\lambda(e^{- t \psi})|^2\Big\}.
  \]
  If $P$ is a polynomial, \eqref{transMult} gives $P(\mathcal{J}(\phi)) = \mathcal{J}(P(\phi))$. The function $\chi_{[r,\infty)}$ 
  may not be a polynomial but we can approximate it by analytic functions as 
  follows. Let $F$ be
  \[
    F(\zeta) = \frac{1}{2} + \frac{1}{\pi} \int_0^{\zeta} e^{- s^2} \, ds.
  \]
  We define the function $\chi_{n,r} \ge 0$ by
  \[
    \chi_{n,r}(\zeta) = \big( F(n \, (\zeta - r)) - F(-n \, r ) \big)^2.
  \]
  For $r>0$, the positive functions $\chi_{r,n}$ converge pointwise and 
  boundedly to $\chi_{[r,\infty)}$ as $n \to \infty$. Furthermore,
  $\chi_{n,r}(0) = 0$ and $\chi_{n,r}$ is a real analytic function
  with arbitrarily large convergence radius. By the analyticity it holds that
  for any radial $\phi$ in the Schwartz class
  \[
    \chi_{n,r}(\mathcal{J}(\phi)) = \mathcal{J}(\chi_{r,n}(\phi)).
  \]
  The right hand side is well-defined since $\chi_{r,n}(\phi)$ is again
  a Schwartz class function and so its Fourier transform is integrable. By
  \cite[Proposition 1.48]{Foll1995} if $\chi_{n,r}$ converges to
  $\chi_{[r,\infty)}$ pointwise and boundedly then $\chi_{n,r}(x)$ converges
  to $\chi_{[0,\infty)}(x)$ is the SOT topology for any positive $x \in \V$.
  We have that
  \begin{eqnarray*}
  \tau \Big\{ \chi_{[r,\infty)}(X_b) \, |\lambda(e^{- t \psi})|^2 \Big\} & = & \sup_{0 < \varepsilon \leq R < \infty} \tau \Big\{ \sotlim_{n \to \infty} \chi_{r,n}(\mathcal{J}(\phi_{\varepsilon, R})) \, |\lambda(e^{- t \psi})|^2 \Big\} \\
                                                                             & = & \sup_{0 < \varepsilon \leq R < \infty} \lim_{n \to \infty} \tau \Big\{ \mathcal{J}(\chi_{r,n} \circ \phi_{\varepsilon, R}) \, |\lambda(e^{- t \psi})|^2 \Big\} \\
                                                                             & \le & \lim_{n \to \infty} \sup_{0 < \varepsilon \leq R < \infty} \tau \Big\{ \mathcal{J}(\chi_{r,n} \circ \phi_{\varepsilon, R}) \, |\lambda(e^{- t \psi})|^2 \Big\}.
  \end{eqnarray*}
On the other hand, $\mathcal{J}$ is trace preserving since $$\tau \circ \mathcal{J}(\phi) = \widehat{\phi} \circ b(e) = \int_{\R^n} \phi \, dm.$$ Moreover, $\lambda(e^{-t\psi}) = \mathcal{J}(h_t)$ for the heat kernel $h_t$ in $\R^n$ and 
  \begin{eqnarray*}
  \tau \Big\{ \chi_{[r,\infty)}(X_b) \, |\lambda(e^{- t \psi})|^2 \Big\}  & = & \lim_{n \to \infty} \sup_{0 < \varepsilon \leq R < \infty} \tau \Big\{ \mathcal{J}(\chi_{r,n} \circ \phi_{\varepsilon, R}) \, |\mathcal{J}(h_t)|^2 \Big\} \\
    & = & \lim_{n \to \infty} \sup_{0 < \varepsilon \leq R < \infty} \tau \Big\{ \mathcal{J} \big( (\chi_{r,n} \circ \phi_{\varepsilon, R}) \, |h_t|^2 \big) \Big\} \\
    & = & \lim_{n \to \infty} \sup_{0 < \varepsilon \leq R < \infty} \int_{\R^n} \chi_{r,n}( \phi_{\varepsilon, R}(\xi) ) \, |h_t(\xi)|^2 d\,\xi \\
    & \leq     & \lim_{n \to \infty} \int_{\R^n} \chi_{r,n}( |\xi| ) \, |h_t(\xi)|^2 d\,\xi \\
    & =        & \int_{\R^n} \chi_{[r,\infty)}( |\xi| ) \, |h_t(\xi)|^2 d\,\xi \ \lesssim \ \frac{1}{\Phi(\sqrt{t})} e^{- \frac{r^2}{2 t}}.
  \end{eqnarray*}
The $\mathrm{CBHL}$ inequality will follow from the \emph{$L_\infty$ Gaussian lower bounds}
  \begin{equation}
    \label{lowerGaussian}
    \hskip-30pt \Big\| \Big( \chi_{[0,r)}(X_b) \, \lambda(e^{-t \psi}) \, \chi_{[0,r)}(X_b) \Big)^{-1} \Big\|^{-1}_{\V} \, \gtrsim \, \frac{e^{- \beta \frac{r^2}{t}}}{\Phi(\sqrt{t})}.
    \tag{$L_\infty \mathrm{GLB}$}
  \end{equation}
  Recall that if $x \in \M_+$ and $p$ is a projection then
  $p \| (p x p)^{-1} \|^{-1} \leq p x p$ and so we can understand
  the right hand side of \eqref{lowerGaussian} as a lower bound on $\chi_{[0,r)}(X_b) \, \lambda(e^{-t \psi}) \, \chi_{[0,r)}(X_b)$. The 
  $L_\infty$GLB allow to bound
  the noncommutative Hardy-Littlewood maximal operator by the maximal operator 
  associated with the semigroup. Indeed, since $X_b$ and 
  $\lambda( e^{-t \psi} )$ commute from \eqref{transMult} we deduce that \eqref{lowerGaussian} yield 
  $$\frac{{\chi_{[0,t)}(X_b)}}{\Phi(t)} \lesssim \chi_{[0,t)}(X_b) \, \lambda(e^{-t^2 \psi}) \, \chi_{[0,t)}(X_b) \le \lambda(e^{-t^2 \psi}).$$ This implies
  \[
    \frac{\chi_{[0,t)}(X_b)}{\Phi(t)} \star x \lesssim S_{t^2}(x),
  \]
  for every positive $x$. Now, using the maximal inequalities 
  for semigroups of \cite{JunXu2007} gives the boundedness of the 
  noncommutative Hardy-Littlewood maximal for every $1 < p < \infty$.
  The fact that $S_t \otimes \Id$ is again a Markovian semigroup
  gives the complete bounds and so the $\mathrm{CBHL}$ inequality holds.
  To prove that \eqref{lowerGaussian} holds we use that $\mathcal{J}: \mathcal{A} \to \V$ is a complete contraction.
  Justifying the calculations like in the case of upper $L_2$ Gaussian bounds and using \eqref{transMult}
  we obtain that
  \begin{eqnarray*}
  \Big\| \Big( \chi_{[0,r)}(X_b) \lambda( e^{-t \psi} ) \chi_{[0,r)}(X_b) \Big)^{-1} \Big\|_{\V}  & \le & \Big\| \left( \lambda( e^{-\frac{t}{2} \psi} )  \widetilde{\chi}_{r}(X_b) \lambda( e^{-\frac{t}{2} \psi} ) \right)^{-1} \Big\|_{\V} \\ & =  & \big\| \mathcal{J} \left( \widetilde{\chi}_{r}(|\cdot|) h^{-1}_t \right) \big\|_{\V}\\
                                                                                                         & \leq & \big\| \chi_{[0,r)}(|\cdot|) h^{-1}_t \big\|_{L_\infty(\R^n)} \\ & \lesssim & t^{\frac{n}{2}} e^{\beta \frac{r^2}{t}}. 
  \end{eqnarray*}
where $\widetilde{\chi}_{r} \in C_c^\infty(\R_+)$ is an smooth decreasing function which is identically $1$ in $[0,r)$ and supported by $[0,2r)$. Taking inverses gives us the desired inequality. \fin


\begin{corollary}
\label{specHMfdcocycles}
Given a LCH amenable unimodular group $G$, let $b: G \to \R^n$ be a finite-dimensional cocycle with associated c.n. length $\psi(g) = |b(g)|^2$. Then, given a symbol $m: \R_+ \to \C$ and $1 < p< \infty$, the following estimate holds for any $\mathcal{H}_0^\infty$ cut-off function $\eta$ and any $s > n/2$ $$\big\|T_{m \circ \psi} \big\|_{\CB(L_p(\V))} \, \lesssim_{(p)} \, \sup_{t > 0} \big\| m(t \, \cdot ) \eta(\cdot) \big\|_{W^{2,s}(\R_+)}.$$
\end{corollary}

\begin{proof}
If the cocycle $b$ is surjective and proper the result follows from Theorem \ref{TheoremB}. Indeed, in that case we know from Theorem \ref{fdStandard} that $(G, \psi, X_b)$ satisfies the standard assumptions with $\Phi(s) = s^{n}$ and Sobolev dimension $D_\Phi = n$. Moreover, the $\mathrm{CBPlan}_2^\Phi$ property also holds as we explained before Lemma \ref{transMult0}. In the general case take $G_{\rtimes} = \R^n \rtimes_{\alpha} G$ where
  $\alpha: G \to O(n)$ is the orthogonal representation that makes 
  $g \mapsto ( x \mapsto \alpha_g x + b(g) )$ an affine representation.
  The function $b_{\rtimes} : G_{\rtimes} \to \R^{n}$ given by 
  $b_{\rtimes}(\xi,g) = \xi + b(g)$ satisfies the cocycle law with cocycle action $\beta: G_\rtimes \to \R^n$ given by $\beta_{(\xi,g)} = \alpha_g$. Indeed, we have
\begin{eqnarray*}
    b_{\rtimes}( \xi + \alpha_g \zeta, g \,h ) & = & \xi + \alpha_g \zeta + b( g h ) \\
                                              & = & \xi + \alpha_g \zeta + \alpha_g b(h ) + b(g )\\
                                              & = & \beta_{(\xi,g)} ( b_{\rtimes}(\zeta,h) ) + b_\rtimes (\xi,g).
  \end{eqnarray*}
  Furthermore $b_{\rtimes}$ is clearly surjective but it may not be 
  proper. In that case, we shall take the associated affine representation 
  $\pi_{\rtimes} : G_{\rtimes} \to \R^n \rtimes O(n)$ and note that 
  the quotient representation $\pi^{\circ}_\rtimes : G_{\rtimes}^{\circ} =
  G_{\rtimes} / \ker(\pi_\rtimes) \to \R^n \rtimes O(n)$ satisfies that its
  associated cocycle $b^{\circ}_\rtimes: G_{\rtimes}^{\circ} \to \R^n$ is
  always proper (even if it is not injective). To see that, let $p_1 : \R^n
  \rtimes O(n) \to \R^n$ be the natural projection into the first component and consider a compact set $K \subset \R^n$. Then
  \[
    (b_{\rtimes}^{\circ})^{-1}[K] =
    (\pi_{\rtimes}^{\circ})^{-1}[p_1^{-1}[K]] =
    (\pi_{\rtimes}^{\circ})^{-1}[K \times O(n)]
  \]
  and the last term is compact since $K \times O(n)$ is compact and
  $\pi_\rtimes^\circ$ is a continuous group isomorphism and hence proper.
  Summing up, we have the following commutative diagram
  \begin{equation*}
    \xymatrix{
      G \ar[r]^{b} \ar@{^{(}->}[d]  & \R^n\\
      \llap{$\R^n \rtimes_\alpha G = $ }G_\rtimes \ar[ru]^{b_\rtimes} \ar@{->>}[d] & \\
      \llap{$(\R^n \rtimes_\alpha G)/\ker(\pi_\rtimes) = $ }G_\rtimes^\circ \ar@/_1pc/[ruu]^{b^\circ_\rtimes} & 
    }
  \end{equation*}
  According to Theorem \ref{fdStandard}, for the last cocycle we can use that $(G_\rtimes^\circ, \psi_\rtimes^\circ, X_{b^{\circ}_\rtimes})$ satisfy the standard
  assumptions, where $\psi_\rtimes^\circ$ is the c.n. length naturally associated to $b_\rtimes^\circ$. By Theorem \ref{TheoremB}, this implies 
  \[
    \big\| T_{m \circ \psi_\rtimes^\circ} \big\|_{\CB (L_p(\V_{\rtimes}^{\circ}) )}
    \, \lesssim_{(p)} \, \sup_{t > 0} \big\| m(t \, \cdot ) \eta(\cdot) \big\|_{W^{2,s}(\R_+)}.
  \]
  Now, using de Leeuw's type periodization 
  \cite[Theorem 8.4 iii)]{CasParPerrRic2014} we obtain the same complete
  bounds for $T_{m \circ \psi_\rtimes}$ in $L_p(\V_{\rtimes})$
  for every $1 < p < \infty$. In order to prove the assertion, 
  we just need to restrict to the subgroup 
  $\{0\} \times G \leq G_{\rtimes}$.This follows from the de Leeuw's restriction type
  result in  \cite[Theorem 8.4 i)]{CasParPerrRic2014}. 
  \end{proof}

\subsection{ \hskip3pt  Foreword}

During the exposition of the contents of Section \ref{Sect2} several  natural
questions arise.
\begin{itemize}
  \item[1.] The first question is whether all finite-dimensional
  proper cocycles, not necessarily surjective, such that their associated c.n.
  length satisfy $\mathrm{CBR_\Phi}$ have $L_2\mathrm{GB}$ for some
  $X \in \V_+^{\wedge}$. We have only been able to prove it in the easier
  case of surjetive cocycles. To that end, our intuition is that a 
  (probably nontrivial) generalization of \eqref{transMult} will be
  required.
  
\vskip3pt

  \item[2.] The second point sprouts from the annoyance of the fact that we
  have not been able to produce explicit examples of infinite-dimensional
  cocycles with $L_2\mathrm{GB}$. We are not confident about their existence. It will be of great interest for us to either construct
  infinite-dimensional cocycles having $L_2\mathrm{GB}$ or to prove that all
  c.n. lengths admitting $X$ with $L_2\mathrm{GB}$ come from finite-dimensional cocycles.
  A way of relaxing such problem is to change the family of c.n. lengths
  arising from finite-dimensional to the family of (real) analytic c.n. lengths (in order
  to make sense of analyticity we will require $G$ to be a Lie group). Note
  that every finite-dimensional cocycle $b : G \to \R^n$ over a Lie group $G$
  induces a group homomorphism of Lie groups $\pi : G \to \R^n \rtimes O(n)$.
  Such homomorphisms are automatically analytic. Therefore, the function
  $\psi : G \to \R_+$ is real analytic. It is reasonable to conjecture that
  every $\psi : G \to \R_+$ defined on a Lie group and with $L_2\mathrm{GB}$
  is analytic.    
  
\vskip3pt

  \item[3.] A possible strategy for constructing conditionally negative lengths coming from
  infinite-dimensional cocycles with $L_2\mathrm{GB}$ is to extend the
  stability results (announced in Remark \ref{RemCrossed}) for crossed products to
  non $\theta$-invariant $\psi_1 : H \to \R_+$ and $X_1 \in \mathcal{L}H_+^{\wedge}$.
  If either $G$ is amenable or $\theta : G \to \Aut(H)$ is an amenable action,
  some sort of averaging procedure may give new c.n. lengths having
  $L_2\mathrm{GB}$ if the original ones do have $L_2\mathrm{GB}$. 
  It will also be desirable to extend the stability of the standard 
  assumptions to extensions of topological groups.

\end{itemize}

\section{{\bf Non-spectral multipliers \label{Sect3}}}

\subsection{ \hskip3pt Polynomial co-growth}

As we have seen, elements in the extended positive cone
$\V^{\wedge}_+$ can be understood as quantized metrics over $\V$. Indeed, 
when $G$ is abelian, any invariant distance over its dual group is determined by the (positive unbounded) function $d(e,\chi)$ affiliated to 
$L_\infty(\widehat{G})$, since $d(\chi_1,\chi_2) = d(e,\chi_1^{-1} \chi_2)$. It may seem natural
to require $X$ to satisfy properties analogous to the triangular inequality, 
the faithfulness and the symmetry. Nevertheless, such assumptions will not 
be necessary here since we will need just \lq\lq asymptotic'' properties of $X$.
Indeed, one of our main families of examples will come from the 
unbounded multiplication symbols of invariant Laplacians over $G$. In order 
to match the classical case of $\R^n$ with the standard Laplacian, whose 
multiplication symbol is $| \xi |^2$, we will use the convention that $X$ 
behaves like $d(e,\chi)^2$. That will explain the $1/2$ exponent in some 
of the formulas. 

\begin{definition}
  \label{defCoGrowth}
  Given $X \in \V_{+}^{\wedge}$, we say that $X$ has polynomial 
  co-growth of order $D$ iff
  \[
    D = \inf \left\{ r > 0 : \big( \1 + X \big)^{-r/2} \in L_1(\V)\right\} < \infty.
  \]
\end{definition}

The definition is motivated by the fact that if we are in an abelian group and
$X$ is the unbounded positive function given by $d(e,\chi)^2$, where $d$ is a 
translation invariant metric then, defining 
$\Phi(r) = \tau(\chi_{[0,r^2)}(X) ) = \mu(B_{r}(e))$, we get $$\big\| \big(\1 + X\big)^{-D/2-\varepsilon} \big\|_1 = \int_{\R_+} \frac{1}{( 1 + r^2 )^{\frac{D}{2} + \varepsilon}} d \Phi(r) = \Big(\frac{D}{2} +\varepsilon \Big) \int_{\R_+} \frac{2 r \Phi(r)}{( 1 + r^2 )^{\frac{D}{2} +1 + \varepsilon}} d r.$$ In particular the last expression is finite whenever 
$\mu(B_r(e)) \lesssim r^D$.

\begin{remark}
  \emph{
    In the proof of Theorem \ref{TheoremC} we are only going to use that the 
    convolution operator $u \mapsto u \star (1 + X)^{- \beta}$ is completely 
    bounded on $L_p(\V)$ for $\beta > D$. Any element in $L_1(\V)$ induces such bounded 
    operator. Indeed we could have defined a similar notion of polynomial 
    co-growth alternatively as
    \[
      D = \inf \left\{ r > 0 : (1 + X)^{-r/2} \in \CB(L_1(\V)) \right\} < \infty,
    \]
    where $(\1 + X)^{-r/2}$ is identified with the operator 
    $x \mapsto (\1 + X)^{-r/2} \star x$. This condition is a priori weaker
    than co-polynomial growth although they coincide for amenable groups. We
    will stick to the original since it is a condition general enough to 
    allow us to prove Theorem \ref{TheoremC} and restrictive enough to be fully 
    characterized.
  }
\end{remark}

Now we are going to prove the existence of unbounded operators affiliated to 
$\V$  behaving like multiplication symbols for left or right invariant 
Laplacians. Recall that a submarkovian semigroup $\S$ acting on $L_\infty(G)$ is respectively called left/right invariant when $S_t \circ \lambda_g = \lambda_g \circ S_t$ or $S_t \circ \rho_g = \rho_g \circ S_t$ accordingly. 

\begin{proposition}
  \label{multSymbols}
  Let $G$ be a LCH unimodular group and consider any submarkovian semigroup $\S$ over $L_\infty(G)$. Let $A$ denote its positive generator. Then, the following properties hold:
  \begin{enumerate}
    \item[i) \label{itm:rightMult}]  If $\S$ is left invariant 
    then there is a densely defined and closable unbounded positive operator 
    $\widehat{A}$ affiliated to $\V$ such that, for all
    $f \in \mathrm{dom}(A) \subset L_2 (G)$
    \[
      \lambda(A f) = \lambda(f)  \widehat{A}.
    \]
    \item[ii) \label{itm:leftMult}] If $\S$ is right invariant 
    then there is densely defined and closable unbounded positive operator 
    $\widehat{A}$ affiliated to $\V$ such that, for all
    $f \in \mathrm{dom}(A) \subset L_2 (G)$
    \[
      \lambda(A f) = \widehat{A} \lambda(f).
    \]
  \end{enumerate}
\end{proposition}

\begin{proof}
  We start by proving $\mathrm{ii)}$. Notice that 
  $A : \mbox{dom} (A) \subset L_2(G) \rightarrow L_2 (G)$ is densely 
  defined. It is affiliated with $\V$ iff for every unitary 
  $u \in \V' = \mathcal{R}G$ we 
  have that $u A = A u$. Since $S_t$ is $\rho$ invariant and 
  we can approximate in
  the SOT topology every element in $\mathcal{R}G$ by linear 
  combinations of elements in $( \rho_{g})_{g\in G}$, we obtain that $S_t$ commutes with any element 
  $x \in \mathcal{R} G$. A function $f \in L_2(G)$ is in $\mathrm{dom}(A)$ when $$\lim_{t \to 0^+} \frac{\Id - S_t}{t} f$$ exists in $L_2(G)$ and we then have
  \[
    \lim_{t \rightarrow 0^+} \Big\| A f - \frac{\Id - S_t}{t} f \Big\|_2 = 0.
  \]
  This implies $u \, \mathrm{dom}(A) \subset \mathrm{dom}(A)$ for any $U (\mathcal{R} G)$. 
  Multiplying by $u$ we obtain $$\big\| u A f - A u f \big\|_2 \le \lim_{t \to 0^+} \Big\| u A f - \frac{\Id - S_t}{t} u f \Big\|_2 + \lim_{t \to 0^+} \Big\| \frac{\Id - S_t}{t} u f - A u f \Big\|_2 = 0$$
  for every $f \in \mathrm{dom}(A)$. This proves 
  that $A$ is affiliated 
  with $\mathcal{R} G$. Notice that $\lambda : L_2(G) \rightarrow L_2(\V)$ 
  unitarily. We will define $\widehat{A} = \lambda A \lambda^{\ast}$. By definition
  $\widehat{A}$ is an unbounded operator on $L_2(\V)$ affiliated with 
  $(\lambda \mathcal{R}G \lambda^{\ast})' = \lambda \V \lambda^{\ast}$ which is
  also equal to the von Neumann algebra $\V$ acting by left
  multiplication in the $\mbox{GNS}$ construction associated to its trace. The
  operator $\widehat{A}$ is densely defined and closable since $A$ is
  densely defined and closable. The identity of $\mathrm{ii)}$
  follows by definition. The construction for $\mathrm{i)}$ is somewhat 
  analogous. We need two trivial observations:
  \begin{itemize}
    \item[1. \label{GNSi}] The anti-automorphism $\sigma:\V \to \V$ extends to a unitary 
    operator $\sigma_2 : L_2(\V) \to L_2(\V)$ since $\tau \circ \sigma = \tau$.
    If $\pi_{r} : \V_{\op} \to \B(L_2(\V))$ and 
    $\pi_{\ell} : \V \to \B(L_2(\V))$ are the right and left GNS 
    representations, then
    $\sigma_2 \circ \pi_r(x) = \pi_{\ell}(\sigma x ) \circ \sigma_2$.

    \item[2. \label{GNSii}] The anti-automorphism $\sigma$ extends to an automorphism 
    of the extended positive cone $\V^{\wedge}_+$. We are going to denote
    such extension again by $\sigma$.
  \end{itemize}    
  Notice that, since $\pi_{\ell}[\V]' = \pi_{r}[\V]$, any element 
  in $x \in \pi_{\ell}[\V]'$ can be expressed as
  $\pi_{r}(x')$ for some $x' \in \V$. By point 1, the map that
  sends $x$ to $x'$ is given, after identifying $\V$ with its GNS 
  representation $\pi_{\ell}[\V]$, by  
  $x' = \sigma ( \sigma_2 \, x \, \sigma_2 )$.  Let $S$ be given by 
  $S = \lambda A \lambda^{\ast}$. Then $S$ is 
  affiliated with $( \lambda \V \lambda^{\ast} )' = \pi_{\ell}[\V]'$.  If we
  define $\widehat{A}$ as $\widehat{A} = \sigma ( \sigma_2 \, S \, \sigma_2 )$, 
  where $\sigma$ is the extension of point 2, we obtain $\mathrm{i)}$.
\end{proof}

\begin{remark}
  \emph{  
    Since $G$ is unimodular, the unitary
    $\iota : L_2(G) \to L_2(G)$ given by $f(g) \mapsto f(g^{-1})$
    is an isometry that intertwines $\rho_g$ and $\lambda_g$. We can characterize
    the pairs of left and right invariant operators $A_1, A_2$ whose
    left and right multiplication symbols, $\widehat{A}_1$ and $\widehat{A}_2$
    respectively, coincide. By a trivial calculation those are the operators such
    that $A_1 \iota = \iota A_2$. Indeed, using that 
    $\lambda : L_2(G) \to L_2(\V)$ satisfies 
    $\lambda \circ \iota = \sigma_2 \circ \lambda$ and that if $A$ is the 
    infinitesimal generator of a submarkovian semigroup then $A^{\top} = A$, we obtain that $$\sigma(\lambda A_1 \lambda^{\ast}) = \sigma_2 \lambda A_2 \lambda^{\ast} \sigma_2 = \lambda \iota A_2 \iota \lambda^{\ast},$$
    but the right hand side satisfies that
    $\sigma(\lambda A_1 \lambda^{\ast}) = \lambda A_1^{\top} \lambda^{\ast} = \lambda A_1 \lambda^{\ast}$.
  }
\end{remark}

Now we are going to characterize those semigroups whose infinitesimal generator has polynomial co-growth. In order to prove the characterization we will need the following two lemmas. Recall that the Fourier algebra $AG$ is defined as those $f: \G \to \C$ such that $\lambda(f) \in L_1(\V)$ with $\|f\|_{AG} = \|\lambda(f)\|_{L_1(\V)}$. We will use below the straightforward inequalities for $f \in A G$
\begin{equation}
  \label{NCHaussYoung}
    | \tau(\lambda(f)) | \leq \| f \|_{\infty} \leq \tau( | \lambda(f) | ).
\end{equation}
Indeed, both follow from the identity $\tau(\lambda_g^* \lambda(f)) = f(g)$ which is valid for $f \in A G$.

\begin{lemma} \label{DensityWLemma}
  Let $G$ be a LCH unimodular group and $\S$ a semigroup of right $($resp. left$)$ invariant operators satisfying that 
  $S_t:C_0(G) \rightarrow C_0(G)$. Let $A$ be the positive generator and assume further 
  that $\widehat{A}$ has polynomial cogrowth of order $D$. Then $W_{A}^{2,s}(G) \cap A G$ is dense inside $W_{A}^{2,s}(G)$ for every 
  $s > {D/2}+\varepsilon$.
\end{lemma}

\begin{proof}
  We will prove only the right invariant case. Notice that $A G$ is
  closed by left and right translations. The fact that $S_t : C_0(G)
  \to C_0(G)$, together with the Riesz representation theorem gives
  that for every $g \in G$ there is weak-$\ast$ continuous family of
  probability measures on $G$, $(\mu_t^{g})_{g \in G, t \geq 0}$ 
  such that
  \[
    S_t f (g) = \int_G f (h) d \mu_t^g (h).
  \]
  Applying the right invariance gives us that $d \mu_t^g ( h) = d \mu_t^e (
  hg^{- 1})$. This yields
  \begin{equation}
    S_t f(g) = \int_G \rho_g f d \mu_t^e = \iota^{\ast} \mu^e_t \ast f (g),
    \label{repConv}
  \end{equation}
  where $(\iota^{\ast} \mu^e_t)(E) = \mu^e_t(E^{-1})$. It is clear that 
  $\| S_t f - f \|_{L_2(G)} \rightarrow 0$ as $t \rightarrow 0^{+}$. 
  Recall that the same is true for $f \in W_{A}^{2, s}(G)$ in the 
  $W_{A}^{2, s} (G)$-norm for every $s > 0$. Suppose that 
  $f \in W^{2, s}_{A}(G)$, then, applying the formula (\ref{repConv}) 
  together with the polynomial co-growth, we have that
 $$S_t f = \iota^{\ast} \mu_t \ast f = \iota^{\ast} \mu_t \ast (\1 + A)^{- \frac{s}{2}}  (\1 + A)^\frac{s}{2} f = h_{t, s} \ast g,$$
  where $g = (\1 + A)^{s/2} f$. We have that $\| g \|_2 = \| f \|_{W^{s, 2}_{A}}$
  and $$\| h_{t, s} \|_2  \le \big\| ( \1 + \widehat{A} \hskip1pt )^{-s/2} \big\|_{L_2 (\V)} \le |\mu_t^e| \big\| (\1+\widehat{A} \hskip1pt )^{-s/2} \big\|_{L_2(\V)} < \infty.$$ 
  This proves that $S_t f \in A G \cap W_{A}^{2, s}(G)$. Making 
  $t \to 0^+$ completes the claim.
\end{proof}


\begin{theorem}
  \label{coGrowthSob}
  Let $G$ be a unimodular LCH group and let $\S$ be a right $($resp. left$)$ 
  invariant submarkovian semigroup over $G$. Let $A$ be its infinitesimal generator and assume further that 
  $S_t : C_0(G) \rightarrow C_0(G)$. Then, the following assertions are 
  equivalent:
  \begin{itemize}
    \item[i)]  The multiplication symbol $\widehat{A}$ 
    of $A$ has polynomial 
    co-growth of order $D$.
    \item[ii)] $\S$ satisfies the following inequality for every $\varepsilon > 0$ 
    \[
      \Big\| ( \1 + A )^{-(\frac{D}{4}+\varepsilon)} : L_2(G) \rightarrow L_\infty(G) \Big\| \lesssim_{(\varepsilon)} 1.
    \]
  \end{itemize}
\end{theorem}

\begin{proof}
  To prove $\mbox{i)} \Rightarrow \mbox{ii)}$, pick 
  $f \in A G \cap W^{2,s}(G)$ for $s = D/2 +2 \varepsilon$ and note
  \begin{eqnarray*}
      \| f \|_\infty & \leq                  & \big\| \lambda \big( (\1 + A)^{-s/2} (\1 + A)^{s/2} f \big) \big\|_{1}  \\
                     & =                     & \big\| \big( \1 + \widehat{A} \hskip1pt \big)^{-s / 2} \lambda \big( (\1 + A)^{s/2} f \big) \big\|_{1} \\
                     & \leq                     & \big\| \big(\1 + \widehat{A} \hskip1pt \big)^{-s / 2} \big\|_{2} \big\| \lambda \big( (\1 + A)^{s/2} f \big) \big\|_{2} \\
                     & =                     & \big\| (\1 + \widehat{A} \hskip1pt )^{-s} \big\|^{1/2}_{1} \| f \|_{W_{A}^{2,s}(G)} \ \lesssim_{(\varepsilon)} \ \| f \|_{W_{A}^{2,s}(G)}.
    \end{eqnarray*}
  We have used \eqref{NCHaussYoung} in the first inequality, Proposition
  \ref{multSymbols} in the first identity and the polynomial cogrowth
  in the last inequality. By the density Lemma \ref{DensityWLemma} 
  we conclude that $W_{A}^{2,s}(G)$ embeds 
  in $L_\infty(G)$ which is a rephrasal 
  of $\mathrm{ii)}$. For the implication $\mbox{ii)} \Rightarrow \mbox{i)}$ 
  we note that from \eqref{NCHaussYoung} $$\Big| \tau \Big( \big( \1 + \widehat{A} \hskip1pt \big)^{- \frac{D}{4} - \varepsilon} \lambda(f) \Big) \Big| \le \big\| \big( \1 + A \big)^{- \frac{D}{4} - 2 \varepsilon} f \big\|_{\infty} \lesssim_{( \varepsilon)} \| f \|_2.$$
  Taking the supremum over $f \in L_2(G)$ with norm 1 gives the desired result. 
\end{proof}

\begin{remark}
  \label{remarkRd}
  \emph{
    Due to Proposition \ref{RdeqSob} we obtain that the point $\mathrm{ii)}$
    is equivalent to satisfying the ultracontractivity property
    $\mathrm{R}_{D+\varepsilon}(0)$ for every $\varepsilon > 0$. Since 
    $\mathrm{R}_D(0)$ implies $\mathrm{R}_{D+\varepsilon}(0)$ for every 
    $\varepsilon > 0$, it is sufficient to prove 
    $\mathrm{R}_D(0)$ in order to have polynomial co-growth of order $D$.
  }
\end{remark}

\begin{remark} \label{LocalvsAssymptotic}
  \emph{
    Sobolev inequalities involving powers of $\1+A$ are sometimes called 
    local \cite[II.X]{VaSaCou1992} since they are tightly connected to 
    the ultracontractivity estimates for $0 < t \leq 1$ and in many contexts 
    that amounts to describing the growth of ball of small radius. Therefore 
    Theorem \ref{coGrowthSob} relates the behaviour of the large balls of $\V$
    with the behaviour of small balls in $G$. This goes along the common 
    intuition that taking group duals exchanges local and asymptotic/coarse 
    properties.
  }
\end{remark}

\demC
  Let $B_t = \lambda( m \, \eta( t \psi ) )$ and let $\widehat{A}_1$ be the multiplication symbol associated with the generator of the right invariant semigroup
  $\S_1$ which is determined by Proposition \ref{multSymbols}. Then 
  \begin{eqnarray*}
    B_t & = & \underbrace{ \big( \1 + \widehat{A}_1 \big)^{-\frac{s_1}{2}} }_{M_t} \underbrace{ \big( \1 + \widehat{A}_1 \big)^{\frac{s_1}{2}} B_t}_{\Sigma_t}
  \end{eqnarray*}
  is a max-square decomposition. By the definition of co-polynomial growth we have that
  $\sigma |M_t|^2 = (\1 + \sigma \widehat{A}_1)^{-s_1} \in L_1(\V)$ and therefore
  it is a c.b. multiplier in every $L_p(\V)$ for $1 \leq p \leq \infty$. Since
  $M_t$ does not depend on $t$, the maximal inequality (MS$_p$) is 
  satisfied trivially. By the construction of $\widehat{A}_1$ we
  have $$\sup_{t > 0} \| \Sigma_t \|_{L_2(\V)} = \sup_{t > 0} \big\| \big( \1 + \widehat{A}_1 \big)^{\frac{s_1}{2}} \lambda( m \eta( t \psi )) \big\|_{L_2(\V)} = \sup_{t > 0} \| m \eta( t \psi ) \|_{W^{2,s_1}_{A_1}(G)}.$$
  The square-max decomposition is manufactured in exactly the same way.
\fin

\subsection{ \hskip3pt Sublaplacians over polynomial-growth Lie groups}

Here we are going to work with left (resp. right) invariant 
submarkovian semigroups over $L_\infty(G)$ generated by sublaplacians.
Let $M$ be a smooth manifold, $\mathbb{X} = \{ X_1,..,X_r\}$ be a family of 
smooth vector fields and $\mu$ a $\sigma$-finite measure over $M$. Let us 
denote by $(\sigma_j(t))_{t \in (-\varepsilon_j,\varepsilon)_j}$ the one-parameter 
diffeomorphism generated by $X_j$ and assume further that $\mu$ is invariant 
under $(\sigma_j(t))_{t \in (- \varepsilon_j, \varepsilon_j)}$. Then, the semigroup
whose infinitesimal generator is given by the sublaplacian associated to 
$\mathbb{X}$
\[
  \Delta_{\mathbb{X}} = - \sum_{j=1}^r X^2_j
\]
is submarkovian. This is a consequence of the theory of 
symmetric Dirichlet forms \cite{FuOTa2010}. If $M = G$ is a Lie group, 
$\mu$ its left Haar measure and $\mathbb{X} = \{X_1, ...,X_r\}$ left
invariant vector fields. By the invariance under the one parameter subgroup
generated by $X_j$ of $\mu$ we have that $S_t = e^{-t \Delta_{\mathbb{X}}}$
is a submarkovian semigrop of left invariant operators.  The same 
construction can be performed using right invariant vector fields if 
$G$ is unimodular. Any sublaplacian carries a natural subriemannian
metric given by
\[
  d_{\mathbb{X}}(x,y) = \inf_{\gamma: [0,1] \to M} \Big\{ \Big( \int_0^1 |\gamma'(t)|^2 \, dt \Big)^{\frac12} \, \big| \ \gamma(0)=x, \gamma(1)=y, \gamma'(t) \in \mathrm{span} \hskip2pt \mathbb{X}(\gamma(t)) \Big\}.
\]
This metric coincides with the Lipschitz distance given by the gradient form,
also known as Meyer's carre de champs \cite{Me1976}. 
Observe also that, if $G$ is a connected Lie group, then its subriemannian 
distance is finite iff $\mathbb{X}$ generates the whole Lie algebra. Similarly, 
$f \in \mathrm{Ker}_p(\Delta_{\mathbb{X}})$ iff $f \in L_p(M)$
and $f(x) = f(y)$ whenever the subriemannian distance 
$d_{\mathbb{X}}(x,y)$ is finite.

The main family of illustrations of Theorem \ref{TheoremC} comes from Lie 
groups endowed with right and left invariant sublaplacians. Indeed, let 
$\mathbb{V} = \{v_1, v_2, ..., v_r\} \subset T_e G$ be a collection of, 
linearly independent, vectors generating the whole Lie algebra and 
$\mathbb{X}_1 = \{X_1, ...,X_r\}$ and $\mathbb{X}_2 = \{Y_1, ...,Y_r\}$ 
be its right and left invariant extensions respectively. Then their 
associated sublaplacians satisfy 
$\iota \Delta_{\mathbb{X}_1} = \Delta_{\mathbb{X}_2} \iota$ where we use 
$\iota f(g) = f(g^{-1})$. Hence,
it suffices to study the polynomial co-growth for
$\widehat{\Delta}_{\mathbb{X}_1}$. By Remark \ref{remarkRd} we just need
to show that $S_t = e^{- t \Delta_{\mathbb{X}}}$ has the $\mathrm{R}_D(0)$
property and by \cite[Theorem VIII.2.9]{VaSaCou1992} we known that if $G$ is
a Lie group of polynomial growth, then
\[
  \frac{e^{- \beta_1 \frac{d_{\mathbb{X}_1}(x,y)^2}{t}}}{\mu (B_e(\sqrt{t}))}
  \lesssim 
  h_t(x,y)
  \lesssim
  \frac{e^{- \beta_2 \frac{d_{\mathbb{X}_1}(x,y)^2}{t}}}{\mu (B_e(\sqrt{t}))},
\]
where $h_t$ is the heat kernel associated with $S_t$, $d_{\mathbb{X}_1}$ is the
subriemannian distance associated to $\mathbb{X}_1$ and $B_e(r)$
are the balls of radius $r$ with respect to that metric. It is a 
well known fact, see \cite{VaSaCou1992}, that
\[
  \mu (B_e(r)) \sim t^{D_0},
\]
for $t$ small. Here $D_0$ is the \emph{local dimension} associated to 
$\mathbb{X}_1$, given by
\[
  D_0 = \sum^{\infty}_{j = 0} j \dim( F_{j + 1} / F_j ),
\]
where $F_0 = \{0\}$, $F_1 = \mathbb{X}_1$ and 
$F_{j+1} = \mathrm{span}\{ F_j, [F_j, \mathbb{X}_1] \}$. As a consequence
$S_t$ has the $R_{D_0}(0)$ property and therefore 
$\widehat{\Delta}_{\mathbb{X}_1}$, and so $\widehat{\Delta}_{\mathbb{X}_2}$, have
polynomial co-growth of order $D_0$. As a corollary we obtain the 
following theorem.

\begin{theorem}
  \label{LieGroup}
  Let $G$ be a polynomial growth Lie group equipped wit a c.n. length $\psi : G \to \R_+$.
  Let $\eta \in \Ha_0$ and consider a generating set $\mathbb{X} = \{ X_1, X_2, ...X_r\}$  
  of independent right invariant vector fields. Let us write $\Delta_{\mathbb{X}}$ for its
  sublaplacian. Then, the following inequality holds for every $1 < p < \infty$ and any $s > D_0/2$ 
  $$\| T_m \|_{\CB(L_p^{\circ}(\V))} \lesssim_{(p)} \sup_{t \geq 0} \max \Big\{ \big\| \eta(t \psi) m 
  \big\|_{W_{\Delta_{\mathbb{X}}}^{2,s}(G)}, \big\| \eta(t \psi) \iota m \big\|_{W_{\Delta_{\mathbb{X}}}^{2,s}(G)} \Big\}.$$
\end{theorem}

\vskip3pt

\noindent \textbf{Acknowledgement.} Gonz\'alez-P\'erez and Parcet are partially supported by the ERC StG-256997-CZOSQP, Junge is partially supported by the NSF DMS-1201886 and all authors are also supported in part by the ICMAT Grant \lq\lq Severo Ochoa" SEV-2011-0087 (Spain).

\bibliographystyle{amsplain}
\bibliography{bibliography.bib}

\null

\vskip10pt

\hfill \noindent \textbf{Adri\'an Gonz\'alez-P\'erez} \\
\null \hfill Instituto de Ciencias Matem{\'a}ticas \\ \null \hfill
CSIC-UAM-UC3M-UCM \\ \null \hfill Universidad Aut\'onoma de Madrid \\ \null \hfill C/ Nicol\'as Cabrera 13-15.
28049, Madrid. Spain \\ \null \hfill\texttt{adrian.gonzalez@icmat.es}

\vskip2pt

\hfill \noindent \textbf{Marius Junge} \\
\null \hfill Department of Mathematics
\\ \null \hfill University of Illinois at Urbana-Champaign \\
\null \hfill 1409 W. Green St. Urbana, IL 61891. USA \\
\null \hfill\texttt{junge@math.uiuc.edu}

\vskip2pt

\enlargethispage{2cm}

\hfill \noindent \textbf{Javier Parcet} \\
\null \hfill Instituto de Ciencias Matem{\'a}ticas \\ \null \hfill
CSIC-UAM-UC3M-UCM \\ \null \hfill Consejo Superior de
Investigaciones Cient{\'\i}ficas \\ \null \hfill C/ Nicol\'as Cabrera 13-15.
28049, Madrid. Spain \\ \null \hfill\texttt{javier.parcet@icmat.es}
\end{document}